\documentclass[oneside,12pt,a4paper,english]{amsart}
\usepackage{pdflscape}
\usepackage{booktabs}
\usepackage{graphicx}
\usepackage{amsmath}
\usepackage{amssymb}
\usepackage{multirow}
\usepackage[table]{xcolor}
\usepackage[english]{babel}
\usepackage{amsmath}
\usepackage{amssymb}
\usepackage{amsthm}
\usepackage{amstext}
\usepackage{hyperref}
\usepackage{xcolor}
\usepackage{diagbox}
\usepackage{placeins}
\usepackage[most]{tcolorbox} %Colored text boxes
\usepackage{geometry}

\definecolor{Blueish}{HTML}{1c4197}

\newtcolorbox{Chaptersetting}[1][]{%
	colback=black!5,
	fonttitle=\bfseries,
	title={#1},
	coltitle=black,
	sharp corners,
	borderline west={2pt}{0pt}{Blueish},
	enhanced,
}

\DeclareMathOperator{\R}{\mathbb{R}}

\DeclareMathOperator{\Z}{\mathbb{Z}}
\DeclareMathOperator{\N}{\mathbb{N}}
\DeclareMathOperator{\C}{\mathbb{C}}

\DeclareMathOperator{\SU}{SU}
\DeclareMathOperator{\SO}{SO}
\DeclareMathOperator{\U}{U}
\DeclareMathOperator{\Uni}{U}
\DeclareMathOperator{\Sp}{Sp}
\DeclareMathOperator{\Spin}{Spin}
\renewcommand{\k}{\ensuremath{\mathfrak{k}}}
\renewcommand{\t}{\ensuremath{\mathfrak{t}}}

\renewcommand{\u}{\ensuremath{\mathfrak{u}}}

\DeclareMathOperator{\g}{\ensuremath{\mathfrak{g}}}

\renewcommand{\sp}{\ensuremath{\mathfrak{sp}}}
\DeclareMathOperator{\so}{\mathfrak{so}}
\DeclareMathOperator{\su}{\mathfrak{su}}

\newcommand{\m}{\ensuremath{\mathfrak{m}}}
\newcommand{\n}{\ensuremath{\mathfrak{n}}}

\newcommand{\z}{\ensuremath{\mathfrak{z}}}
\newcommand{\h}{\ensuremath{\mathfrak{h}}}
\newcommand{\myrestriction}{{\,\restriction\,}}
\newcommand{\mycolon}{{\,:\,}}
\DeclareMathOperator{\spann}{span}
\DeclareMathOperator{\x}{\mathfrak{X}}
\renewcommand{\c}{\ensuremath{\mathcal{C}}}
\DeclareMathOperator{\tr}{tr}
\DeclareMathOperator{\Ad}{Ad}

\newcommand{\Ca}{\operatorname{C}}
\renewcommand{\i}{\mathrm{i}}

\DeclareMathOperator{\Ric}{Ric}
\DeclareMathOperator{\vol}{vol}

\DeclareMathOperator{\id}{Id}

\DeclareMathOperator{\sgn}{sgn}

\DeclareMathOperator{\mult}{mult}

\renewcommand{\r}[1]{{#1}^{\bullet}}
\DeclareMathOperator{\End}{End}
\newcommand{\modulo}[2]{{#1}/{#2}} %das entspricht dem normalen A/B, die beiden Kommandos werden leider nicht eideutig verwendet, machen allerdings das selbe.
\newcommand{\spec}{\Sigma}
\renewcommand{\hat}[1]{\widehat{#1}}

\theoremstyle{plain}
\newtheorem{thm}{Theorem}[section]
\newtheorem*{thm*}{Theorem}
\newtheorem{lem}[thm]{Lemma}
\newtheorem*{lem*}{Lemma}
\newtheorem{prop}[thm]{Proposition}
\newtheorem{cor}[thm]{Corollary}

\theoremstyle{definition}
\newtheorem{rem}[thm]{Remark}
\newtheorem{dfn}[thm]{Definition}
\newtheorem{stand}[thm]{Standing Assumption}

\newtheorem{exa}[thm]{Example}

\usepackage{graphicx}

\usepackage{tikz}

\begin{document}

	\title[The Laplace-Beltrami spectrum on Naturally Reductive Spaces]{The Laplace-Beltrami spectrum on Naturally Reductive Homogeneous Spaces}
	
	\author{Ilka Agricola and Jonas Henkel}
	
	\address{\hspace{-7mm} 
Ilka Agricola, Jonas Henkel,
Fachbereich Mathematik und Informatik,
Philipps-Universit\"at Marburg,
Campus Lahnberge,
35032 Marburg, Germany\newline
{\normalfont\ttfamily agricola@mathematik.uni-marburg.de, henkelj@mathematik.uni-marburg.de}}


	\maketitle
    \begin{center}
   \today
    \end{center}
    
 \begin{abstract}
 We prove a formula for the spectrum of the Laplace-Beltrami operator on functions
 for compact naturally reductive homogeneous spaces in terms of eigenvalues of a generalized Casimir operator and spherical representations.  We apply this result to a large family of canonical variations of normal homogeneous metrics, thus allowing for the first time to study  how the spectrum depends on the deformation parameters of the metric. As an application, we provide a formula for the full spectrum of compact positive homogeneous $3$-$(\alpha,\delta)$-Sasaki manifolds (a family of metrics which includes, in particular, all homogeneous $3$-Sasaki manifolds).

The second part of the paper is devoted to  the detailed computation and investigation of the spectrum of this family of metrics on the Aloff-Wallach manifold $W^{1,1}=\SU(3)/S^{1}$; in particular, we provide a documented Python script that allows the explicit computation in any desired range.  We recover Urakawa's eigenvalue computation for the $\SU(3)$-normal homogeneous metric on $W^{1,1}$ as a limiting case and cover all the positively curved $\SU(3)\times\SO(3)$-normal homogeneous realizations discovered by Wilking. By doing so, we complete Urakawa's list of the first eigenvalue on compact, simply connected, normal homogeneous spaces with positive sectional curvature.
\end{abstract}

	\tableofcontents

	%Kein Einrücken
	\setlength{\parindent}{0em}

\section{Introduction}% and overview of results}
%-----------------------------------------------
Let $(M,g)$ be a closed Riemannian manifold. We consider the Laplace operator $\Delta$ acting on functions; being an elliptic self-adjoint operator, its spectrum consists of discrete eigenvalues of finite
multiplicity. 
The explicit computation of the spectrum is a challenging task that only succeeds under strong symmetry assumptions. After the classical examples of spheres, projective spaces, and flat tori (see for example \cite{BergerGM}), Riemannian symmetric   
spaces $M=G/K$ were the first large class of manifolds for which the spectrum could be computed explicitly via representation theoretic tools. In this case, the connection induced from the canonical principal fibre bundle projection $G\rightarrow G/K$ is just the Levi-Civita connection, and the Laplacian can be identified with 
the Casimir operator $\Ca$, viewed as an element of the center of the universal enveloping algebra $\U (\mathfrak{g})$. Hence, a combination of the Peter-Weyl theorem for describing 
$L^2(G/K)$ and the Freudenthal formula for the eigenvalue of $\Ca$ on irreducible representations yields the result in a very satisfying way  \cite{CahnWolf76, Ikeda78}, and basically the same approach works for compact semisimple Lie groups endowed with the negative of the Killing form just as well \cite{BeersM77} (see also \cite{Boucetta09} for a concise  account of both cases).  

Let us recall the definition of two main classes of invariant metrics $g$ on
non-symmetric homogeneous reductive spaces $G/K$ in order to formulate our results in more detail. We
fix a reductive decomposition $\g=\k\oplus\m$, and 
use the same letter $g$ to denote the $\Ad(K)$-invariant inner product on $\m$ which induces the Riemannian metric $g$. The metric $g$ is called \emph{naturally reductive} if it satisfies
\[
g(x,[y,z]_\m) + g([y,x]_\m, z) = 0 \quad \text{ for all } \ x,y,z\in \m\cong T_o G/K.
\]
This is in particular the case if the metric $g$ on $\m$ is the restriction of an 
$\Ad(G)$-invariant positive definite inner product on $\g$ and $\m=\k^\perp$ with respect to that inner product. The homogeneous space is then called \emph{normal}, and if in addition
$G$ is semisimple and the inner product in question is just the negative of the Killing form, one
often speaks of a \emph{standard} normal homogeneous space. Although the definitions seem similar, the (standard) normal homogeneous case is much easier to treat than the naturally reductive one. An adaptation of the argument given above allowed, in principle, the computation of the Laplace spectrum on standard normal homogeneous spaces, a good source for which is the
seminal book \cite{Wallach harm ana} (see also the book \cite{Hall} and \cite{Bettiol} for a recent generalisation to metrics induced from Riemannian metrics on $G$ that are left $G$- and right $K$-invariant). However, only Riemannian naturally reductive metrics give the freedom needed to cover interesting \emph{families} of metrics depending on parameters. Thus, our class of manifolds is tailor-made for investigating the dependence of the spectrum on certain deformation parameters of the metric.

\medskip
The present paper is organized in two  parts. The first part deals with the general theory,
while the second part is devoted to the $7$-dimensional Aloff-Wallach manifold $W^{1,1}=\SU(3)/S^1$, which was in fact the starting point of this work (the reasons are given below), one of the most interesting manifolds admitting metrics of strictly positive sectional curvature. While the explicit computation of the required spherical representations is highly intricate and carried out here explicitly for $W^{1,1}$, the authors and L. Cagliero have recently generalized this method in a subsequent paper \cite{AgricolaCaglieroHenkel26}. By establishing a simplified spectral branching criterion and combining a reciprocity relation with the Littlewood-Richardson rule, the full spectra for the entire classical series of 3-$(\alpha,\delta)$-Sasaki manifolds as well as Stiefel manifolds have now been completely determined.

Let us summarise the most important results of the first part. We begin with two key observations: Let $(G/K,g)$ be a reductive Riemannian homogeneous space. %\\
%\textcolor{red}{alt: }A reductive homogeneous Riemannian space $( G/K,g)$ is a Riemannian submersion $G\rightarrow G/K$  with totally geodesic fibres,
%for which there exists a concept of \emph{canonical variation} of the metric.\\
%\\
%\textcolor{red}{neu:} A reductive homogeneous Riemannian space $( G/K,g)$ is a Riemannian submersion $G\rightarrow G/K$  with totally geodesic fibres. For such submersions the Laplace-Beltrami operator commutes with the projection map: Eigenfunctions of $G$ which are constant along the right $K$-orbits descend to eigenfunctions of the quotient space $G/K$. 
If $H\subset G$ is a subgroup which commutes with $K$ and whose adjoint action preserves $g$, the Riemannian submersion $G/K\rightarrow G/(K\times H)$ has totally geodesic fibres. The deformation of the metric $g$ of the quotient space $G/K$ along the fibres $H/K$ gives rise to the notion of \textit{canonical variation} of the metric. 
It is known that the Laplace-Beltrami operator may then be split into a horizontal and vertical part, and hence its spectrum behaves particularly well for metrics that are canonical variations (Subsection \ref{subsec:submersions}). Secondly,  many known interesting families of homogeneous metrics occur as canonical variations of normal homogeneous metrics. However, the realization as an abstract canonical variation is not enough to fully \emph{compute} the spectrum, and non-trivial canonical variations of normal homogeneous spaces are rarely again normal (and it is usually difficult to decide whether they are or not), so the classical results on spectra of normal homogeneous metrics do not apply.
Thus, we  prove in a first step a generalization of the Freudenthal formula (Theorem \ref{thm eigenvalues casimir}) for \emph{any}
naturally reductive homogeneous metric $(G/K,g)$, which relies on a suitable 
generalization of the classical Casimir operator (Definition \ref{dfn:gen_Casimir}).
In a second step, we define a large class of canonical variations of normal homogeneous metrics \emph{that are naturally reductive for another presentation $G'/K'$}, thus closing the gap between the two approaches (Subsection \ref{subsec:nrhm}). We prove that the corresponding families of metrics on $G/K$ and $G'/K'$  are isometric under a careful identification of parameters (Theorem \ref{thm: isometry main result}). Theorem 
\ref{thm: spectrum main result} then states an explicit formula for the spectrum in terms of $K'$-spherical representations that makes the dependency on the deformation parameters of the metric totally explicit and reduces the computation to an intricate, but theoretically solvable problem of representation theory. As an application, we conclude a formula for the full Laplace-Beltrami spectrum of a  homogeneous
$3$-$(\alpha,\delta)$-Sasaki manifold in Corollary \ref{cor: spectrum 3 alpha delta}, and we obtain a new lower bound on its first basic eigenvalue that we compare to the recent results of \cite{Semmelmann Nagy}.  

The second part is devoted to the Aloff-Wallach manifold $W^{1,1}=\SU(3)/S^1$ solely.
Aloff and Wallach proved that this manifold admits a one-parameter family of metrics of strictly positive sectional curvature \cite{Aloff Wallach 75}. Almost 25 years passed before Wilking realized that these metrics are actually $\SU(3)\times\SO(3)$-normal homogeneous, thus closing a gap in a classification of these metrics by Berger \cite{Berger}.
Based on the construction given in the first part,
we extend this normal homogeneous realization to a naturally reductive one which gives access to the spectrum of all $3$-$(\alpha,\delta)$-Sasaki metrics on $W^{1,1}$ (Theorem \ref{thm: spectrum Aloff Wallach}). 
As a limiting case, we reproduce the classical 
computation of Urakawa for the $\SU(3)$-normal homogeneous metric \cite{Urakawa}, which corresponds to one point in our family of metrics. Moreover, we complete Urakawa's list \cite{Urakawa list} of the first eigenvalue of a compact simply connected normal homogeneous space with positive curvature in Corollary \ref{cor: list of first eigenvalue on compact simply connected normal hom pos curvature}.
As an electronic supplement, we provide a Jupyter Notebook with a Python script that allows the computation of the spectrum in any desired range, and 
Table \ref{table: first eigenvalues, reps and mult} gives a comprehensive  list of  eigenvalues with the highest weights of
the respective spherical representations and the multiplicities.
We relate the spectrum to the geometry and existing geometric structures in Subsection \ref{subsection: geometric interpretation of the spectrum} and we use the results to check eigenvalue estimates for canonical variations of $3$-Sasaki manifolds in Subsection \ref{subsection: interpretation as canonical variation}.
In particular, we provide graphs of the eigenvalues that make the underlying Riemannian submersion clearly
`visible'.
Finally, we study the spectrum for fixed volume in Subsection \ref{subsection: spectrum with constant volume}. We will show that an analogous inequality to Hersch \cite{Hersch} does not exist on $3$-$(\alpha,\delta)$-Sasaki manifolds by constructing a $1$-parameter family of metrics for which the product of the first eigenvalue with the volume is neither lower nor upper bounded.

\subsubsection*{Acknowledgements and declarations.} 
It is our pleasure to thank Louis Ioos, Emilio Lauret, Henrik Naujoks, Paul Schwahn, Uwe Semmelmann,  and Leander Stecker for very  valuable discussions on different aspects of this article and for useful references to the literature. This material is partially based upon work supported by the National Science
Foundation under Grant No. DMS-1928930, while the authors were in
residence at the Simons Laufer Mathematical Sciences Institute
(formerly MSRI) in Berkeley, California, during the Fall/2024 semester. The second author was supported by a fellowship of the German Academic Exchange Service (DAAD).
Both authors contributed equally to this work and  approved the final manuscript. The Python code and data are made freely available in a Jupyter notebook \cite{Agricola Henkel 24} prepared by the second author. Both authors declare to have no competing interests.

%
%----------------------------------------------------------------
\section{The spectrum of a naturally reductive homogeneous space}
%----------------------------------------------------------------
%
In this article, $G$ will always denote a compact, connected Lie group with Lie algebra $\g$ of dimension $n$ and $ K\subset G $ a closed subgroup with Lie algebra $\k$ of dimension $k$. 
We choose a reductive decomposition $\g=\m\oplus \k$ and assume that $G$ is equipped with an $\Ad(G)$-invariant 
\emph{semi-Riemannian} metric $g$ which induces a positive definite scalar product on $\m$
that makes $G/K$ a naturally reductive homogeneous space with respect to $\m$.

\begin{rem}
Every naturally reductive homogeneous space can be obtained that way: 
Suppose that $(\tilde{G}/\tilde{K},\tilde{g})$ is a naturally reductive 
Riemannian homogeneous space  with respect to a reductive complement $ \tilde{\m}\subset \tilde{\g}$. By a theorem of Kostant \cite[Thm.\,4]{Kostant}, 
the scalar product $ \tilde{g}|_{ \tilde{\m}} $ equals $g|_{ \tilde{\m}}$, where $g$ is an 
$\Ad(G)$-invariant semi-Riemannian metric of a subgroup $ G\subset \tilde{G} $ with Lie algebra $\g:=[\tilde{\m},\tilde{\m}]+ \tilde{\m}$. The subalgebra $\k=\g\cap\: \tilde{\k}$ is orthogonal to $\tilde{\m}$ with respect to $\tilde{g}$ and $G$ acts transitively on $ \tilde{G}/\tilde{K}$ with isotropy $K$. The spaces $ (G/K,g)$ and $ ( \tilde{G}/\tilde{K},\tilde{g}) $ are isometric.
\end{rem}
Let us now decompose  $ \g=\z\oplus_{j=1}^{l}\g_{j} $     into its center $\z$ and into simple Lie algebras $\g_{j}$.
Then it is well-known that these are mutually  orthogonal with respect to $g$ and that on each simple factor $\g_{j}$, the scalar product $g$ restricts to a non-zero multiple of the Killing form $B$: 
\begin{align*}
g |_{\g_{j}}=-r_{j}B|_{ \g_{j}}\quad \text{for $r_{j}\in \R^{*}$}.
\end{align*} 
Thus, the non-degenerate symmetric bilinear form $ g $ restricts to a positive definite $( \sgn(r_{j})>0) $ or a negative definite $ ( \sgn(r_{j})<0 ) $ bilinear form on $ \g_{j} $.

\medskip
Using the theory of $G$-invariant differential operators, we introduce the generalized Casimir element which encodes the signature of $g$ and identify the Laplacian with it.
\subsection{$G$-invariant differential operators and the generalized Casimir operator}
%-------------------------------------------------------------------------------------
%
The isometric and transitive action of $G$ on itself from the left induces the left regular representation $l_{a}(f)(b):=f(a^{-1}b)$. It is a unitary representation which restricts to the quotient $G/K$ and induces a unitary representation on $ L^{2}(G/K)$ by $(af)(bK):=f(a^{-1}bK)$. 
\begin{dfn}
A differential operator $A$ on $\mathcal{C}^{\infty}(G)$ resp. $\mathcal{C}^{\infty}(G/K)$ with analytic coefficients is called \emph{$G$-invariant} if it commutes with the left regular representation of $\c^{\infty}(G)$, resp.\,its induced representation on $\c^{\infty}(G/K)$. The set of analytic $G$-invariant differential operators is denoted by $\mathcal{D}(G)$ or $\mathcal{D}(G/K)$, respectively. The subset $\mathcal{D}_{K}(G)\subset \mathcal{D}(G)$ consists of operators that are invariant under the pullback of the right action of $K$.
\end{dfn}
%$G$-invariant differential operators commute with the $G$ action on $\c^{\infty}(G/K)\subset L^{2}(G/K)$ and 
As $G$ acts transitively, the action of a $G$-invariant operator is fully determined once
its action on functions at the neutral element $e$ is known: For $f\in \c^{\infty}(G)$ and $a\in G$, we have $A(f)(a)=A(l_{a^{-1}}(f))(e)$. The left-invariant vector fields in $ \g $ act on functions $f\in  \c^{\infty}(G)$ by the Lie group exponential:
\[
X(f)(a) :=X_{a}(f)=\frac{d}{dt}\big|_{t=0}f(a\exp(tX))
\]
and thus define $G$-invariant differential operators, $\g\subset \mathcal{D}(G)$. Lichnerowicz \cite{Lichnerowicz} considers more generally $G$-invariant connections on $G$ and defines the derivative of $f$ along $X$ using the corresponding exponential map. This corresponds here to the choice of the Ambrose-Singer connection, as in this setting the Lie group exponential map coincides with the geodesic exponential map. By the universal property, the action extends to the universal algebra and generates all $G$-invariant differential operators on $\c^{\infty}(G)$:% This enables to identify them with elements in the universal enveloping algebra. 
\begin{lem}[{\cite[p. 266]{Helgason}}]
%-------------------------------------
The homomorphism $\g\rightarrow \mathcal{D}(G)$ extends uniquely to an isomorphism $U(\g)\rightarrow \mathcal{D}(G)$.% the universal enveloping algebra of $\g$ into $D(\g)$.
\end{lem} 
Functions on $G/K$ can of course be identified with $K$-right-invariant functions on $G$,
\[ 
\c^{\infty}(G/K)\cong \c^{\infty}_{K}(G):=\{f\in \c^{\infty}(G)\mid \:f(ak)=f(a)\: \forall k\in K,\: \forall a\in G\},
\]
hence the algebra $\mathcal{D}(G/K)$ of $G$-invariant differential operators on $G/K$
is isomorphic to the algebra $\mathcal{D}_{K}(G)$ restricted to $\c^{\infty}_{K}(G)$.
Note that indeed, $A\in\mathcal{U}(\g)\cong \mathcal{D}(G) $ restricts to $\c^{\infty}(G/K)$ if and only if $A\in \mathcal{D}_{K}(G)\cong \mathcal{D}(G/K)$ which is precisely the case if $\Ad_{k}A=A$ for all $k\in K$. 

As the Laplacian commutes with isometries, it defines a $G$-invariant differential operator on $G/K$. The signature of $g|_{\k}$ enters as soon as it is identified with an element in $U(\g)$.
\begin{dfn}[Generalized Casimir operator]\label{dfn:gen_Casimir}
%---------------------------------------------------------------
Let $ x_{j} $, $ j=1,\dots,n $ be an orthonormal basis of $\g$ such that $x_{j}$, $j=1,\dots,n-k$ is an orthonormal basis of $\m$ and $x_{j}$, $j=n-k,\dots,n$ is an orthonormal basis of $\k$. Their length is given by $ \varepsilon_{j}:=g(x_{j},x_{j})\in \{\pm 1\}$, in particular $\varepsilon_{j}=1$ for $j=1,\dots,n-k$. The \emph{generalized Casimir operator} is defined by 
\[
\Ca_{g} :=\sum_{j=1}^{n}\varepsilon_{j} x_{j}^{2}\in U(\g).
\]
\end{dfn}
\begin{rem}
The classical Casimir element is usually defined for simple Lie algebras equipped with the negative of the Killing form; for standard normal homogeneous metrics, no signs 
$\varepsilon_j$ occur. A similar definition of the generalized Casimir element can be found in \cite{Yamaguchi}.
\end{rem}
As in the Riemannian case, the generalized Casimir element is $\Ad_{G}(K)$-invariant and hence restricts to the quotient $\c^{\infty}(G/K)$.
Despite the appearance of the signs $\varepsilon_j$, $\Ca_g$ may still be identified with the Laplacian:
%-------------
\begin{lem}\label{lem: casimir invariant Ad}
The generalized Casimir element $ \Ca_{g} $ is independent of the choice of orthonormal basis of $ \g. $ If $ a\in G $ then $ \Ad_{a}(\Ca_{g})=\Ca_{g} $ and hence $\Ca_{g}\in \mathcal{D}(G/K)$. As a differential operator, it yields the Laplace-Beltrami operator: $-\Ca_{g}=\Delta\in \mathcal{D}(G/K)$. 
\end{lem}

\begin{proof}
If $ \{y_{i} \mid i=1,\dots,n \}$ is a different orthonormal basis of $ \g $ with length $\varepsilon_{i}\in \{\pm 1\} $, then $ y_{i}=\sum a_{ji}x_{j}=Ax_{i} $ with $ (a_{ij})=A\in O(g)$. As  $A$ preserves the inner product, it satisfies $ (\sum_{l=1}^{n}a_{ml}a_{sl}\varepsilon_{l})_{ms}=(\varepsilon_{m}\delta_{ms})_{ms} $. Hence we get
\begin{align*}
	\sum_{i=1}^{n}\varepsilon_{i} y_{i}^{2}=\sum_{i,j,s=1}^{n}\varepsilon_{i}a_{ji}a_{si}x_{j}x_{s}=\sum_{j,s=1}^{n} \varepsilon_{j}\delta_{js}x_{j}x_{s}=\sum_{i=1}^{n}\varepsilon_{i}x_{i}^{2}=\Ca_{g}.
\end{align*}
This proves that $\Ca_{g}$ does not depend on the chosen orthonormal basis. For $ a\in G $ set $ y_{i}=\Ad_{a}x_{i} $ and obtain $\Ca_{g}=\Ad_{a}(\Ca_{g})$. This holds in particular for any $a\in K$ and which proves that $\Ca_{g}\in \mathcal{D}(G/K)$. As $\Delta$ and $-\Ca_{g}$ are $G$-invariant differential operators, we simply prove that any $f\in \c^{\infty}(G/K)$ satisfies $\Delta (f) (eK)=-\Ca_{g} (f)(eK)$ and conclude that $-\Ca_{g}=\Delta$. As semi-Riemannian submersions map horizontal geodesics to geodesics \cite[Cor. 46]{ONeill}, the map $(t_{1},\dots,t_{n-k})\mapsto \exp(\sum_{i=1}^{n-k}t_{i}x_{i})eK$ defines local geodesic coordinates, i.e. $g_{eK}=\id$ and $\Gamma_{ij}^{s}(eK)=0$ for all $i,j,s$. Thus, the Laplacian simplifies in $eK$ to
\begin{align*}
	\Delta (f)(eK) =&\left.-\sum_{j=1}^{n-k}\left.\frac{d^{2}}{dt^{2}}f(\exp(tx_{j})eK)\right|_{t=0}-0\right. \\
	=& \left. -\sum_{j=1}^{n-k} \left. \frac{d^2}{dt^2} f(\exp(t x_j) eK) \right|_{t=0} -\sum_{j=n-k+1}^{n} \varepsilon_{j}\left. \frac{d^2}{dt^{2}} f(\exp(t x_j) eK) \right|_{t=0} \right.\\
	=&-\Ca_{g}(f)(eK).
\end{align*} 
\end{proof}
If $g$ is only $\Ad_{G}(K)$-invariant, the element $\Ca_{g}$ is still $\Ad_{G}(K)$-invariant and it can still be identified with the Laplacian $\Delta$ on $G/K$. However, it turns out that $\Ca_{g}$ is not in the center of $\U(\g)$ and therefore its spectrum is more difficult to compute, see Remark \ref{rem: only Ad K invariant eigenvalues}. The left-regular representation of $G$ on $L^{2}(G/K)$ decomposes into irreducible representations which are preserved by the Laplacian. Each of those representations contribute to the spectrum.  
%\begin{rem}

%\end{rem}
\subsection{The Freudenthal formula for the generalized Casimir operator}
%------------------------------------------------------------------------
Our main goal in this subsection is to prove an analogue of the classical Freudenthal formula in the more general setting of naturally reductve homogeneous spaces.
We recall some general facts and introduce some notation.
\begin{itemize}
\item Denote by $\hat{G}$ the set of irreducible, unitary representations of the Lie group $G$. It is known that these are finite dimensional and are in $1{:}1$-correspondence to dominant integral weights, whose set is denoted by $D(G)$.

\item Let $(\varrho_{1},V_{\varrho_{1}}),(\varrho_{2},V_{\varrho_{2}})$ be unitary representations of $G$. Then $[\varrho_{1} \mycolon \varrho_{2}]$ denotes the multiplicity of $\varrho_{2}$ in $\varrho_{1}$, i.e. the dimension of $G$-equivariant homomorphisms from $V_{\varrho_{1}}$ to $V_{\varrho_{2}}$. 
\item The restriction of a representation to a subgroup $K\subset G$ is denoted by 
$\varrho \myrestriction K$ and the trivial representation is denoted by $1_{K}$.
\item A representation $(\varrho,V_{\varrho})\in \hat{G}$ is called 
\emph{$K$-spherical} if $[\varrho \myrestriction  K \mycolon 1_{K}]>0$. The subspace $V_{\varrho}^{K}\subset V_{\varrho}$ is the set of fixed vectors when $\varrho$ is restricted to $K$. The set of $K$-spherical representations is denoted by $ \hat{G}_{K} $.
\item For $(\varrho,V_{\varrho})\in \hat{G}_{K}$ the number $\mult(\varrho):=[\varrho \myrestriction K \mycolon 1_{K}]\cdot \dim(V_{\varrho})=\dim(V_{\varrho}^{K})\cdot \dim(V_{\varrho})$ describes the multiplicity of eigenfunctions being generated by $\varrho$, see Theorem \ref{thm eigenvalues casimir}.
\item By $T\subset G$, we denote a maximal torus of $G$ with Lie algebra $\t$ and the dual space $\t^{*}$ is equipped with the dual of $ g(\cdot,\cdot)\vert_{\t} $ which is extended to a non-degenerate Hermitian bilinear form on its complexification $\t^{*}_{\C}$ and denoted by $g$ as well. We choose a set of positive roots $R^{+}$ of $G$ and denote their half sum by $\rho$.
\end{itemize}

As $(\varrho,V_{\varrho})$ is unitary, we can use the scalar product to identify $(V_{\varrho}^{K})^{*}\cong (V_{\varrho}^{*})^{K}$ as $\R$-vector spaces. Any $K$-spherical representation defines matrix coefficients which are invariant under $K$ and which hence descend to functions on $G/K$:
\begin{align*}
V_{\varrho}\otimes (V_{\varrho}^{*})^{K}\rightarrow L^{2}(G/K),\quad v\otimes \varphi \mapsto f_{v\otimes \varphi}, \quad f_{v\otimes \varphi}(aK):=\varphi(\varrho(a^{-1})v) .	
\end{align*}
This assignment is a $G$-homomorphism if $G$ acts on the first entry of the tensor product $ a(v,\varphi):=(\varrho(a)v,\varphi)$ and $G$ acts on $L^{2}(G/K)$ by pullback. The Peter-Weyl theorem states that the direct sum of matrix coefficients is dense in $L^{2}(G/K)$:
\begin{thm}[{\cite[Thm. 1.3]{Takeuchi}}]\label{thm: frobenius reciprocity}
$ L^{2}(G/K) $ decomposes as a $ G $-module as follows: 
\begin{align*}
	L^{2}(G/K)\cong\overline{\oplus_{\varrho\in \hat{G}_{K}}V_{\varrho}\otimes (V_{\varrho}^{*})^{K}}.
\end{align*}% where $ f_{v\otimes\varphi}(a)=\varphi(\varrho(a^{-1}v)) $.
\end{thm}
It is crucial to understand how the Laplacian acts on the matrix coefficients. As the Laplacian commutes with isometries, it restricts to a finite dimensional, self-adjoint endomorphism on $V_{\varrho}\otimes (V_{\varrho}^{*})^{K} $. We describe its action in $eK$:
\begin{align*}
(\Delta f_{v\otimes \varphi})(eK)&=-\sum_{j=1}^{n-k}\frac{d^{2}}{dt^{2}}\varphi(\varrho(\exp(tx_{j})eK)^{-1}v)=-\sum_{j=1}^{n-k}\frac{d^{2}}{dt^{2}}\varrho(\exp(tx_{j}))(\varphi)(v)\\
&=-\sum_{j=1}^{n}\varepsilon_{j}\frac{d^{2}}{dt^{2}}\varrho(\exp(tx_{j}))(\varphi)(v)=-d\varrho(\Ca_{g})(\varphi)(v).
\end{align*}
Observe that $d\varrho(\Ca_{g})$ acts on $(V_{\varrho}^{*})^{K}$ and each eigenvector $\varphi$ of $-d\varrho(\Ca_{g})$ with eigenvalue $\eta$ yields the eigenspace $ V_{\varrho}\otimes \{\varphi\}$ of the Laplacian corresponding to the same eigenvalue $\eta$.
\begin{rem}\label{rem: only Ad K invariant eigenvalues}
If we relax the condition of $g$ from being $\Ad(G)$-invariant to being $\Ad_{G}(K)$-invariant, the eigenvalues of the Laplacian are given by the eigenvalues of $-d\varrho(\Ca_{g})$ on $(V_{\varrho}^{*})^{K}$. Lauret \cite{Lauret} was able to explicitly compute the first eigenvalue of $\SU(2)$ for all left-invariant Riemannian metrics by explicity realizing and diagonalizing $-d\varrho(\Ca_{g})$ for the first $\SU(2)$-representations. This method was later used to compute the first eigenvalue of any homogeneous metric of the CROSSes \cite{Bettiol Distance Spheres, Bettiol}. In Example \ref{exa: hom metrics spheres} we describe which of the homogeneous CROSS metrics can be realized as naturally reductive spaces.
\end{rem}
In our situation of considering a biinvariant semi-Riemannian metric, the generalized Casimir acts as a scalar which is analogous to the normal homogeneous setting:
\begin{lem}\label{lem: casimir acts as scalar}
For any $(\varrho,V_{\varrho})\in \hat{G}$ the operator $-d\varrho(\Ca_{g})$ commutes with $\varrho(a)$ for all $a\in G$ and hence acts a multiple $c_{g}(\varrho)$ of the identity.
\end{lem} 
\begin{proof}
The key ingredient is that we consider a biinvariant semi-Riemannian metric, i.e. $\Ad_{a}\in O(g)$ which was used in Lemma \ref{lem: casimir invariant Ad} to prove $\Ad_{a}(\Ca_{g})=\Ca_{g} $ for all $a\in G$.	Using the relation $ d\varrho \circ \Ad_{a}=d(\varrho\circ c_{a})=d(c_{\varrho(a)}\circ \varrho)=\Ad_{\varrho(a)}\circ d\varrho $ we have that
\begin{align*}
	d\varrho(\Ca_{g})=d\varrho\circ \Ad_{a}(\Ca_{g})=\Ad_{\varrho(a)}d\varrho(\Ca_{g}).
\end{align*}
Applying Schur's Lemma implies that $-d\varrho(\Ca_{g})$ is a multiple of the identity. 
\end{proof}
We state the main theorem of the first section which computes the constant $c_{g}(\varrho)$. This generalises \cite[Lem. 5.6.4]{Wallach harm ana} where $ g $ has no signature and coincides on $ \g_{0} $ with the Killing form.
\begin{thm}[Generalized Freudenthal formula]\label{thm eigenvalues casimir}
%---------------------------------------------------------------------------
Let $ (\varrho,V_{\varrho})\in \hat{G} $ with corresponding highest weight $ \lambda $. The operator $-d\varrho(\Ca_{g})$ acts as a multiple $c_{g}(\varrho):=g(\lambda,2\rho+\lambda)$ of the identity on $ (V_{\varrho}^{*})^{K}\subset L^{2}(G)$. The spectrum of the Laplace-Beltrami operator is given by 
\begin{align*}
	\spec(G/K,g)=\{c_{g}(\varrho)\mid \varrho\in \hat{G}_{K}\}.
\end{align*} The multiplicity associated to $(\varrho,V_{\varrho})\in \hat{G}_{K}$ is given by $\mult(\varrho)= [\varrho\myrestriction K\mycolon 1_{K}]\cdot\dim(V_{\varrho}) $.
\end{thm}

\begin{proof}
%------------
We mainly follow the strategy of the proof from \cite[Lem. 5.6.4]{Wallach harm ana}. We decompose $\g=\g_{0}\oplus\z$, where $\g_{0}$ is semisimple and $\z$ is the center of $\g$. % We generalize the proof in the case that $ \g $ is equipped with an arbitrary $ \Ad(G) $-invariant metric with signature.
Note that $ g $ is definite on $ \g_{i} $, non-degenerate on $ \z =\z\cap \t$ and hence non-degenerate on $ \t=\z\oplus_{i=1}^{l}  \g_{i}\cap\t $ and therefore defines a scalar product with signature on $\t$. We want to compute the constant $c_{g}$ from Lemma \ref{lem: casimir acts as scalar}.\\
%	 \\
%	
The main idea is to choose an orthonormal basis of $ \g $ which is obtained using the root decomposition of $ \g^{\C} $ and then evaluate $- d\varrho(\Ca_{g}) $ on the weight vector of the highest weight. 
Consider the decomposition of $ \g_{0}=\oplus_{j=1}^{l}\g_{j} $ into simple subalgebras $ \g_{j} $. Let $ y_{(1,0)},\dots,y_{(k_{0},0)} $ be a orthonormal basis of $ \mathfrak{z} $ and let $ y_{(1,j)},\dots,y_{(k_{j},j)} $ be an orthonormal basis of $ \g_{j} $. Set $ \Ca_{0}=\sum_{i=1}^{k_{0}}\varepsilon_{(i,0)} y_{(i,0)}^{2} $ and set $ \Ca_{j}=\sgn(r_{j})\sum_{i=1}^{k_{j}} y_{(i,j)}^{2} $. Then $ \Ca_{g}=\sum_{j=0}^{l}\Ca_{j} $. Because each subspace is $ \Ad(G) $-invariant, we have by the same argument as in the proof of Lemma \ref{lem: casimir invariant Ad} that $ \Ad(a)\Ca_{j}=\Ca_{j} $ for all $ a\in G $ and $ j=1,\dots,l. $ Let $\lambda'   $ be the linear form on $ \mathfrak{z} $ so that $ d\varrho(z)v=\lambda'(z)v $ for all $ z\in \mathfrak{z} $, $ v\in V $, i.e. $ d\varrho(\Ca_{0})=\sum_{j=1}^{k_{0}}\varepsilon_{(j,0)} \lambda'(y_{(j,0)})^{2}\id. $\\
Let $ R $ be the root system of $ G $ relative to $ \t $. Let $ \t_{j}=\t\cap\g_{j} $ and $ \g_{j}^{\C} $ as well as $ \t_{j}^{\C} $ denote the complexifications. The root system $ R $ decomposes into the root systems $ R=\cup_{j=0}^{l} R_{j} $ of the respective Lie algebras $ \g_{j} $ with respect to $ \t_{j} $ \cite[Thm. 7.35]{Hall}. Consider the root decomposition of $ \g_{j}^{\C}=\t_{j}^{\C}+\sum_{\alpha\in R_{j}} \g_{\alpha} $ and the $ \C $-bilinear extension of the Killing form $ B_{j}^{\C} $ of $ \g_{j}^{\C} $. Then $ B_{j}^{\C}(\g_{\alpha},\g_{\beta})=0 $ if and only if $ \alpha\neq -\beta $ \cite[Lem. 3.5.6]{Wallach harm ana}. Moreover, if $ w_{\alpha}=u_{\alpha}+iv_{\alpha} \in \g_{\alpha} $ then $ w_{\alpha}^{*}:=iv_{\alpha}-u_{\alpha}\in \g_{-\alpha} $ \cite[Prop. 7.18]{Hall} and thus $ B_{j}^{\C}(w_{\alpha},w_{\alpha})=B_{j}^{\C}(w_{\alpha}^{*},w_{\alpha}^{*})=0 $. Since $ B_{j}^{\C} $ is non degenerate on $\g_{\alpha}\oplus\g_{-\alpha}$, we can choose $ w_{\alpha},w_{\alpha}^{*} $ such that $ B_{j}^{\C}(w_{\alpha},w_{\alpha}^{*})=-B(u_{\alpha},u_{\alpha})-B(v_{\alpha},v_{\alpha})= \frac{2}{|r_{j}|}$. This is up to scaling the same choice as Wallach made. By expressing $u_{\alpha}$ and $v_{\alpha}$ in terms of $w_{\alpha},w_{\alpha}^{*}$, one concludes
\begin{align*}
	B_{j}(u_{\alpha},u_{\alpha})=B_{j}(v_{\alpha},v_{\alpha})=-\frac{1}{|r_{j}|},\quad B_{j}(u_{\alpha},v_{\alpha})=0\\
	\Rightarrow g(u_{\alpha},u_{\alpha})=g(v_{\alpha},v_{\alpha})=\sgn(r_{j}),\quad	g(u_{\alpha},v_{\alpha})=0 .
\end{align*}
Moreover, we have
\begin{align*}
	&w_{\alpha}\otimes w_{\alpha}^{*}+w_{\alpha}^{*}\otimes w_{\alpha}=-u_{\alpha}\otimes u_{\alpha}-v_{\alpha}\otimes v_{\alpha}=-(u_{\alpha})^{2}-(v_{\alpha})^{2}.
\end{align*} Taking an orthogonal basis $ h_{(1,j)},\dots, h_{(s(j),j)} \in \sqrt{-1}\t_{j}$ of length $\sgn(r_{j}) $, we write each summand of the generalized Casimir element $ \Ca_{g} $ as follows:
\begin{align*}
	-\Ca_{j}=-\sum_{i=1}^{s(j)}\sgn(r_{j}) h_{(i,j)}^{2}-\sum_{\alpha\in R_{j}^{+}}\sgn(r_{j})w_{\alpha}w_{\alpha}^{*}+\sgn(r_{j})w_{\alpha}^{*}w_{\alpha}.
\end{align*}
Using $ w_{\alpha}w_{\alpha}^{*}=[w_{\alpha},w_{\alpha}^{*}]+w_{\alpha}^{*}w_{\alpha}=\frac{H_{\alpha}}{|r_{j}|}+w_{\alpha}^{*}w_{\alpha} $ \cite[Lem. 3.5.8]{Wallach harm ana} results in
\begin{align*}
	-\Ca_{j}=-\sum_{i=1}^{s(j)}\sgn(r_{j}) h_{(i,j)}^{2}-\sum_{\alpha\in R_{j}^{+}}\sgn(r_{j})\frac{H_{\alpha}}{|r_{j}|}+2\sgn(r_{j})w_{\alpha}^{*}w_{\alpha}.
\end{align*}
The weight vector $ 0\neq v\in V_{\varrho} $ of the highest weight $ \lambda=\lambda_{0}+\lambda' $ satisfies $d\varrho(\n^{+})v=0  $, where $ \mathfrak{n}^{+}:= \sum_{\alpha\in R^{+}}\g_{\alpha}$ \cite[Lem. 4.3.7]{Wallach harm ana}. In order to compute the eigenvalue $ c_{g}(\varrho) $ of the Casimir element, we evaluate $ d\varrho(\Ca_{g})(v) $ and use that $ d\varrho(w_{\alpha})v=0 $ for every $ \alpha\in R^{+} $:
\begin{align*}		
	-d\varrho(\Ca_{g})(v)&=-d\varrho(C_{0})(v)-\sum_{j=1}^{l}d\varrho(\Ca_{l})v\\
	&=-\sum_{j=1}^{k_{0}}\varepsilon_{(j,0)} \lambda'(y_{(j,0)})^{2} v-\sum_{\substack{1 \leq j\leq l\\1\leq i\leq s(j)}}\sgn(r_{j}) \lambda_{0}(h_{(i,j)})^{2}v-\sum_{\substack{1\leq j\leq l\\\alpha\in R^{+}_{j}}}\frac{\lambda_{0}(H_{\alpha})}{r_{j}}v\\
	&=g(\lambda',\lambda')v+g(\lambda_{0},\lambda_{0})v+\sum_{\alpha\in R^{+}}g(\lambda_{0},\alpha)v=g(\lambda,2\rho+\lambda)v.\\
\end{align*}
The whole spectrum is obtained by the virtue of Theorem \ref{thm: isometry main result}.
\end{proof}

\subsection{Digression: The Laplacian for Riemannian submersions}
\label{subsec:submersions}
%---------------------------------------------------------------
%
For the family of metrics we shall consider in the next subsection, the broader view as canonical variations of a Riemannian submersion with totally geodesic fibres will turn out to be crucial. We quickly summarize all relevant results.

Consider a Riemannian submersion $\pi: M\rightarrow N$ with respect to two Riemannian manifolds $(M,g)$, $(N,h)$. It yields a vertical $\mathcal{V}=\ker d\pi$ and horizontal $\mathcal{H}=\mathcal{V}^{\perp}$ distribution of the tangent bundle $TM=\mathcal{V}\oplus_{\perp}\mathcal{H}$. The fibre of $p\in M$ over $q=f(p)\in N$ is denoted by $\mathcal{V}_{p}=f^{-1}(q) $. The Laplacian yields two differential operators on $M$: The vertical Laplacian $\Delta^{v}$ which derives functions along the fibres and the horizontal Laplacian $\Delta^{h}$ given by its difference to $ \Delta_{M}$:
\begin{align*}
\Delta^{v}f(p):=\Delta_{\mathcal{V}_{p}}(f\myrestriction \mathcal{V}_{p}),\quad \Delta^{h}:=\Delta_{M}-\Delta^{v}.
\end{align*}
It is known that all three operators $\Delta_{M},\Delta^{v},\Delta^{h}$ commute and hence yield a Hilbert basis of $L^{2}(M)$ consisting of simultaneous eigenfunctions of $\Delta_{M}$ and $\Delta^{v}$ \cite[Thm. 3.6]{Bergery_et_al}. In case of a Riemannian submersion with minimal fibres, any eigenfunction $f$ of $N$ lifts to an eigenfunction $f\circ \pi$ on $M$ if the fibres are minimal: 
\begin{thm}[{\cite[Thm. 3]{Gilkey Park}}]
%----------------------------------------
The Laplacian commutes with the projection: $\Delta_{M}\circ\nolinebreak \pi^{*}=\pi^{*}\Delta_{N}$ if and only if the fibres of $\pi$ are minimal.
\end{thm}
If the Riemannian submersion has totally geodesic fibres, the fibres are isometric. Each of those Riemannian submersions  give rise to a canonical variation of the fibre given by the deformation of $g$ along the fibres $\mathcal{V}_{p}$, $p\in M$: 
\begin{dfn}
%------------
The \emph{canonical variation} of a Riemannian submersion $\pi:M\rightarrow N$ with totally geodesic fibres is the family of metrics $g_{t_{0},t_{1}}$, $t_{0},t_{1}>0$ on $M$ such that
\begin{align*}
	g_{t_{0},t_{1}}\left|_{\mathcal{V}}\right. = t_{1}g\left|_{\mathcal{V} }\right., \quad
	g_{t_{0},t_{1}}\left|_{\mathcal{H}}\right. =t_{0} g\left|_{\mathcal{H} }\right.,\quad \mathcal{V}\perp\mathcal{H}.
\end{align*}

The Laplacian of $ (M,g_{t_{0},t_{1}}) $ is denoted by $ \Delta_{M}^{t_{0},t_{1}} $. 
\end{dfn}
The rescaling of the fibre in the canonical variation rescales the vertical Laplacian $\Delta^{v}$ while it  preserves the horizontal Laplacian $\Delta^{h}$: 
\begin{prop}[{\cite[Prop. 5.3]{Bergery_et_al}}]
%----------------------------------------------
The Laplacian $\Delta_{M}^{t_{0},t_{1}}$ of the canonical variation is given by
\begin{align*}
	\Delta_{M}^{t_{0},t_{1}}=\frac{\Delta^{v}}{t_{1}}+\frac{\Delta^{h}}{t_{0}}=\Delta_{M}^{t_{0},t_{1}}+\left(\frac{1}{t_{1}}-\frac{1}{t_{0}}\right)\Delta^{v}.
\end{align*}
\end{prop}
As the fibres are isometric we conclude:
\begin{cor}\label{cor: spectrum canonical variation}
The spectrum of the canonical variation $(M,g_{t_{0},t_{1}})$ satisfies: %is contained in the rescaled spectrum \( (M,g) \) and the difference of the spectrum of the fibre $(\mathcal{V}_{p},g|_{\mathcal{V}_{p}})$ and the spectrum of its rescalation $(\mathcal{V}_{p},tg|_{\mathcal{V}_{p}})$:
$$\spec(M,g_{t_{0},t_{1}})\subset \{t^{-1}_{0}\eta_{0}+(t^{-1}_{1}-t^{-1}_{0}) \eta_{1}\mid \eta_{0}\in \spec(M,g),\quad \eta_{1}\in \spec(\mathcal{V}_{p},g|_{\mathcal{V}_{p}})\},\quad p\in M.$$
\end{cor}

The mismatching of both sets measures the topology of the fibre bundle. If the bundle is trivial, the two sets are equal. The key problem is to determine which eigenfunctions of the fibre $\mathcal{V}_{p}$ are obtained when simultaneous eigenfunctions of $\Delta_{M}$ and $\Delta^{v}$ are restricted to it; this cannot be decided for a general canonical variation, but we shall be able to do it in the context of certain canonical variations of normal homogeneous metrics that we shall discuss in the next subsection.

For the moment being, let us state that indeed, the following homogeneous spaces are Riemannian submersions with totally geodesic fibres:
\begin{prop}[{\cite[Thm. 5.2]{Falcitelli}}]\label{prop: hom space riem subm totally geod}
Let $G$ be a compact Lie group, $K_{1}\subset K_{2}\subset G$ be subgroups and equip $G$ with a left $G$- and right $K_{2}$-invariant Riemannian metric $g$. Choose the unique metrics on $G/K_{1}$ and on $G/K_{2}$  such that \begin{align*}
	\pi_{1}:G\rightarrow G/K_{1},\quad \pi_{2}:G/K_{1}\rightarrow G/K_{2}
\end{align*} are  Riemannian submersions. Then $\pi_{1}$ and $\pi_{2}$ have totally geodesic fibres.
\end{prop}
\subsection{Naturally reductive canonical variations of normal homogeneous metrics and their spectrum}\label{subsection naturally reductive realizations}\label{subsec:nrhm}
%------------------------------------------------------------------------------------------
%
Many interesting families of homogeneous metrics occur as canonical variations of normal homogeneous spaces; to cite a few, let us mention
certain metrics of positive curvature constructed by Aloff and Wallach \cite[Thm. 2.4]{Aloff Wallach 75},
homogeneous metrics on compact rank one symmetric spaces \cite{Ziller 82 Metrics crosses, Bettiol}, including in particular odd-dimensional Berger spheres (see also Example \ref{exa: hom metrics spheres}), many homogeneous spaces carrying Killing spinors \cite{BFGK}, and more recently positive homogeneous $3$-$(\alpha,\delta)$-Sasaki manifolds \cite{Agricola Dileo Stecker 21, Goertsches Roschig Stecker}.

It is desirable to realize them as naturally reductive spaces in order to be able to compute their spectrum by means of Theorem \ref{thm: spectrum main result}. It can be tricky to decide whether a homogeneous space $(G/K,g)$ can be realized as a naturally reductive space: From Berger's classification of compact, simply connected normal homogeneous spaces \cite{Berger}, it took almost $40$ years until Wilking \cite{Wilking} realized that the positively curved Aloff-Wallach manifold $W^{1,1}$ admits a realization as a normal homogeneous space, and thus was missing from the list.
Hence, our goal in this subsection is to provide a method to realize canonical variations of normal homogeneous spaces as naturally reductive spaces and compute their spectrum. %On the level of Lie groups, D'Atri and Ziller \cite{Atri Ziller} give a similar construction to realize deformations of biinvariant metrics as naturally reductive metrics on Lie groups. 
\begin{stand}
In this section, we assume that the Lie group $G$ is equipped with a biinvariant metric $g$. By $K,H_{i}\subset G$, we denote pairwise commuting subgroups of $G$ and we assume  $H_{i}$ connected for all $i=1,\ldots,s$. Let $\k,\h_{i}\subset \g$, $i=1,\dots,s$ denote their Lie algebras, which shall be pairwise $g$-orthogonal. Choose the reductive complement $\m=\k^{\perp}\subset \g$ and define the direct sum $\h=\oplus_{i=1}^{s}\h_{i}$. Its orthogonal complement in $\m$ is denoted by $\h_{0}=\h^{\perp}\subset \m$. Note that $\h_{0}$ plays a special role as it does not define a Lie algebra. We obtain the reductive decomposition $\g=\m\oplus\k$, $\m=\h_{0}\oplus\h$.
\end{stand}
As an intermediate step, we deform the metric $g$ on $G$ along the subspaces $\k$, $\h_{i}\subset\g$ and denote the new metric  by $g_{\textbf{t}}, \ 
\textbf{t}=(t_0,t_1, \ldots, t_s)\in \R^{s+1}_{>0}$:
\begin{align*}
(G,g_{\textbf{t}}),\quad g_{\textbf{t}}:=g_{t_{0},\ldots,t_{s}}=t_{0}g|_{\h_{0}\oplus\k}+ t_{1}g|_{\h_{1}} \ldots +t_{s}g|_{\h_{s}},\quad \k\perp\h_{j},\quad \h_{j}\perp \h_{i}\: \text{for all}\: i\neq j.
\end{align*}
As $\Ad_{G}(K)$ preserves $\h_{i}$, the metric $g_{\textbf{t}}$ restricts to a well-defined homogeneous metric of the quotient $G/K$:
\begin{align}\label{eq: Einführung der Metrik g t}
(\modulo{G}{K},g_{\textbf{t}}),\quad	g_{\textbf{t}}:=g_{t_{0},\dots,t_{s}}=t_{0}g|_{\h_{0}}+\dots+t_{s}g|_{\h_{s}},\quad \h_{i}\perp \h_{j}\: \text{for all}\: i\neq j.
\end{align}
\begin{rem}
The condition that the groups $H_{i},K$ commute for each $i=1,\dots,s$ guarantees that $H_{i}\times K$ is isomorphic to a subgroup of $G$. By Proposition \ref{prop: hom space riem subm totally geod}, $$\pi_{i}:G/K\rightarrow G/(H_{i}\times K)$$ is a Riemannian submersion with totally geodesic fibres. Its vertical and horizontal spaces are, for $i=1,\ldots, s$:
\[
\mathcal{V}=\h_i,\quad \mathcal{H}=\h_0 \bigoplus_{j=1,j\neq i }^s \h_j.
\]
The canonical variation of $\pi_i$ is then exactly the metric $g_\textbf{t}$,
where $t_i$  (for fixed $i$) is to be understood as the scaling of the fibration and all other parameters as describing the metric on the horizontal space.
\end{rem}
Such deformations occur in a natural way. Consider a homogeneous space $G/K$ with reductive complement $\m$. Assume that the isotropy representation $\Ad_{G}(K)\myrestriction \m$ acts trivially on at least one vector, $[\Ad_G(K)\myrestriction \m\mycolon 1_{K}]\geq 1 $, and thus defines a global vector field on $G/K$. Using the Jacobi identity, one verifies that the subspace $\m^{K}=\{X\in \m\mid \forall k\in K:\: \Ad_{G}(k)(X)=X \}$ defines a Lie subalgebra $\m^{K}=:\h$. Decomposing $\h$ orthogonally into commuting subalgebras $\h_{1},\dots,\h_{s}\subset \h$ and deforming $g$ along them, yields a deformation of the mentioned kind. In case of a compact rank one symmetric space, it follows by definition that such deformations can at most be performed along one subalgebra, i.e. $s=1$ and hence $\h_{1}=\h$.

\begin{exa}\label{exa: hom metrics spheres}
We consider a negative multiple of the Killing form $g$ of the simple Lie group $G=\SU(n+1)$ together with its subgroups %\textcolor{red}{the inclusions are given by the obvious matrix blocks embeddings}
\begin{align*}
	K=\left\{\begin{pmatrix}
		1&0\\
		0&A 
	\end{pmatrix}\mid A\in \SU(n)\right\},\quad H=\left\{\begin{pmatrix}
		e^{it}&0\\
		0& e^{\frac{-it}{n}}\id_{n}
	\end{pmatrix}\mid e^{it}\in S^{1}\right\}, 
\end{align*} i.e. $G/K\cong S^{2n+1}$ and $G/(H\times K)\cong \mathbb{CP}^{n-1} $. The deformation $g_{t_{0},t_{1}}$ yields the canonical variation of the Hopf fibration
\begin{align*}
	S^{1}\rightarrow	S^{2n+1}\rightarrow  \mathbb{CP}^{n}.
\end{align*} 
%	One refers to $(S^{2n+1},g_{t_{0},t_{1}})$ as Berger spheres.
Another example is given by the group $G=\Sp(n+1)$ equipped with a negative multiple of the Killing form $g$ and the subgroups
\begin{align*}
	K=\left\{\begin{pmatrix}
		1&0\\
		0& A
	\end{pmatrix}\mid A\in \Sp(n)\right\},\quad H=\left\{\begin{pmatrix}
		B&0\\
		0&1
	\end{pmatrix}\mid B\in \Sp(1)\right\}.
\end{align*}
The deformation $g_{t_{0},t_{1}}$ is the canonical variation of the quaternionic Hopf fibration
\begin{align*}
	S^{3}\rightarrow	S^{4n+3}\rightarrow  \mathbb{HP}^{n}
\end{align*} 
which yields at the same time the positive $3$-$(\alpha,\delta)$-Sasaki metrics on $S^{4n+3}$, see Subsection \nolinebreak \ref{subsection: The Spectrum of Sasaki Manifolds}. Note that this canonical variation descends to a one-parameter family of metrics on $\mathbb{CP}^{2n+1}=S^{4n+3}/\U(1)$, as $\U(1)\subset \Sp(1)$ commutes with $\Sp(n)$. Moreover, both canonical variations descend to the real projective spaces $\mathbb{RP}^{2n+1}$ and $\mathbb{RP}^{4n+3}$  as both fibre groups $S^{1}$ and $\Sp(1)$  commute with $\Z_{2}$. Observe that the 
canonical variation of the octonionic Hopf fibration
\begin{align*}
	S^{7} \rightarrow S^{8+7}\rightarrow \operatorname{Ca}\mathbb{P}^{2} %\text{\textcolor{red}{wie genau geht das eig? $\Spin(9)/\Spin(7)$}}
\end{align*}
and its quotient $\mathbb{RP}^{15}$ are not of this type, as the isotropy representation of $\Spin(7)$ on $S^{15}=\Spin(9)/\Spin(7)$ does not contain the trivial representation \cite{Ziller 82 Metrics crosses}.
Furthermore, one checks that a different rescaling of the Killing vector fields (if more than one) can not be obtained by our construction for similar reasons. All in all, the metrics we just described are the only homogeneous CROSS metrics: They have been classified by Ziller \cite{Ziller 82 Metrics crosses} and their spectra are studied in \cite{Bettiol Distance Spheres, Bettiol}.
%
% projection and hence descend to the quotient.
\end{exa}
%\begin{exa}
%	Let $k,n\in \N$, $n>k$ and $H_{1},\dots,H_{s}\subset O(k)$ be commuting %subgroups. The Stiefel manifold $V_{k}^{n}(\R^{n})=O(n)/O(n-k)$ equipped %with $g(X,Y)=\tr(XY)$ can be deformed along $H_{1},\dots,H_{s}$. 
%\end{exa}
%
%
Our next goal is to  give an explicit 
naturally reductive realization $(G'/K',h_{\textbf{r}})$ of the deformed normal homogeneous space $(\modulo{G}{K},g_{\textbf{t}})$.
We introduce the group $H:=H_{1}\times \dots\times H_{s}$ and the multiplication $\tilde{\iota}:H\rightarrow G$, $\tilde{\iota}(h_{1},\dots,h_{s})=h_{1}\dots h_{s}$, which is an injective  group homomorphism as the groups $H_{i}$ commute pairwise. Using additionally that $H_{i}$ and $K$ commute, we construct the following groups: \begin{align*}
G':=G\times H, \quad \Delta H:=\{(\tilde{\iota}(h),h)\mid h\in H\},\:\: 
K':=\{(k\tilde{\iota}(h),h)\mid k\in K,\: h\in H\}\subset G'.
\end{align*}
To make the final construction a bit more transparent, we make an 
intermediate step. Define an
$\Ad(G') $-invariant metric $h_\textbf{r}$ on $\modulo{G'}{\Delta H}$
depending on $s+1$ real parameters 
$\textbf{r}:=(r_{0},\dots,r_{s})\in (\R^{*})^s$ by
\begin{align*}
(\modulo{G'}{\Delta H},h_\textbf{r}),\quad h_{{\textbf{r}}}:=h_{r_{0},\dots,r_{s}}=r_{0}g|_{\g}+r_{1}g|_{\h_{1}}+\dots+r_{s}g|_{\h_{s}},\quad
\g\perp \h_{i} \ \text{ and }\ \h_{i}\perp \h_{j}\ \forall i\neq j.
\end{align*}
This makes $\modulo{G'}{\Delta H}$ a semi-Riemannian homogeneous space;
it is essentially a copy of $G$, of course.

In order to define an orthonormal basis of the Lie algebra $\g\oplus\h$ of $G'$, we need to define indices.
%The map $d\iota: \h\rightarrow\g$ denotes the inclusion. 
We set $m_{i}=\dim \h_{i}$, $\sum_{i=1}^{s}m_{i}=m$ and introduce sets of indices corresponding to $\h_{0}$: $I_{0}=\{m+1,\dots,n\}$, to $\h_{1}$: $I_{1}=\{1,\dots,m_{1}\}$ and for $j>1$ to $\h_{j}$: ${I}_{j}=\{\sum_{i=1}^{j-1}m_{i},\dots, \sum_{i=1}^{j}m_{i}\}$, their union is denoted by $I=\cup_{j=1}^{s}I_{j}$. We define with respect to $g$ orthonormal bases of $\h_{j}$, $j\geq 0$ and $\g$:
\begin{align*}
\{e_{i}\mid i\in {I}_{j}\}\subset \h_{j},\quad \{E_{i}\mid i=1,\dots,n\}\subset \g,\quad d\iota e_{i}=E_{i}.
\end{align*} 
%of $\h_{j}$ for $j>0$ and orthonormal basis $\{E_{i}\mid i=1,\dots,n\}$ of $\g$ such that
%\begin{align*}
%\forall \: i\in I:  \:\:d\iota e_{i}=E_{i} ,\quad \spann_{\R}\{E_{m+1},\dots,E_{n}\}=d\iota(\h)^{\perp}\subset \g.
%\end{align*} 
Their lengths in $(\g\oplus \h,h_{\textbf{r}})$ and in $(\g,g_{\textbf{t}})$ are given by ($ i\in I_{j} $, $ j\geq 0 $):
\begin{align*}
g_{\textbf{t}}(E_{i},E_{i})=t_{j},\quad h_{\textbf{r}}(E_{i},E_{i})=r_{0},\quad \text{ for }j\geq 1:\: h_{\textbf{r}}(e_{i},e_{i})=r_{j}
\end{align*}
%\begin{align*}
%h_{\textbf{r}}(e_{i},e_{i})=r_{j},\quad h_{\textbf{r}}(E_{i},E_{i})=h_{\textbf{r}}(E_{k},E_{k})=r_{0},\quad g_{\textbf{t}}(E_{i},E_{i})=t_{j},\quad g_{\textbf{t}}(E_{k},E_{k})=t_{0}
%\end{align*}
and all vectors are pairwise orthogonal.
A basis for the isotropy algebra is given by
%	=\left\{\left(\begin{bmatrix}
%		A&0\\
%		0&-\tr(A)
%	\end{bmatrix},A+\s^{1}\:\vline\: A\in \u(2)\right)\right\}\\
\begin{align*}
%\su(3)\oplus\so(3)& \cong \su(3)\oplus \spann_{\R}\{e_{1},e_{2},e_{3}\}\\
\Delta\h=(d\iota,\id)(\h)=\spann_{\R}\left\{E_{i}+e_{i}\mid i\in I\right\}.
\end{align*}
The following example of the Aloff-Wallach space inspired our constructions and will be discussed in detail in Section \ref{section: Aloff-Wallach manifold}. In fact, we shall explain in Subsection \ref{subsection: The Spectrum of Sasaki Manifolds} that it can even be generalised to
all positive homogeneous $3$-$(\alpha,\delta)$-Sasaki manifolds.
% \textcolor{red}{Verbindungssätze, warum funktioniert das, was ist subtiler Punkt dabei?}
\begin{exa}\label{exa: incomplete aloff wallach}
Consider the Lie group $G=\SU(3) $ and its commuting subgroups $K=S^{1},H=\SU(2)\subset \SU(3)$, where $S^{1}$ is embedded diagonally and $\SU(2)$ embedded in the upper left $2\times2$ block. They yield the subgroup $(S^{1})'=\{(SA,A)\mid S\in S^{1},A\in U(2)\}\subset \SU(3)\times \SU(2)$. The metric $g(X,Y)=-1/2\tr(XY)$ is the rescaled Killing form of $\SU(3)$. The resulting Riemannian homogeneous space is the Aloff-Wallach space $W^{1,1}$ equipped with the canonical variation of a $3$-Sasaki manifold:
\begin{align*}
	\left(\modulo{\SU(3)}{S^{1}},g_{t_{0},t_{1}}\right),\quad g_{t_{0},t_{1}}=t_{0}g|_{\su(2)^{\perp}}+t_{1}g|_{\su(2)},\quad t_{0},t_{1}>0.
\end{align*}
The semi-Riemannian homogeneous space associated to the deformation $g_{\textbf{t}}$ is
\begin{align*}
	\left(\modulo{\SU(3)\times \SU(2)}{\Delta \SU(2)},h_{r_{0},r_{1}}\right),\quad h_{r_{0},r_{1}}=r_{0}g|_{\g}+r_{1}g|_{\su(2)},\quad  r_{0},r_{1}\geq 0 \text{ or } r_{1}<-r_{0}<0.
\end{align*}
\end{exa}
The Reeb vector fields of a homogeneous Sasaki manifold are per definition invariant under the isotropy representation. Hence, the unique connected subgroup which realizes those commutes with the isotropy group. By extension of the transitively acting group we are able to give naturally reductive realizations of the canonical variation $g_{\textbf{t}}$, see Theorem \nolinebreak\ref{thm: isometry main result}. In the next step, we need to ensure that $h_{\textbf{r}}$ defines
a Riemannian naturally reductive metric on
$G'/\Delta H$ and hence also on the quotient $G'/K'$:
\begin{lem}\label{lem: pos definite bilinear form}
%-------------------------------------------------
%The Lie algebra $ (\so(3)\oplus \su(3),h_{r_{1},r_{2}}) $ is isometric to $ (\su(2)\oplus \su(3),4r_{1}\left.B\right\vert_{\su(2)} +r_{2}\left.B\right\vert_{\su(3)})$. 

The bilinear form $h_{\textbf{r}} $ is non-degenerate on $ \Delta\h$ if and only if $r_{j}\neq -r_{0}$ for all $j=1,\dots,s$ and 
\begin{align*}
	\Delta\h^{\perp}=\spann_{\R}\{ h_{\textbf{r}}(e_{j},e_{j})E_{j}-r_{0}e_{j},E_{i}\mid j\in I,\: i\in I_{0}\}
\end{align*} is a reductive complement. $h_{\textbf{r}}$ is positive definite on $\Delta\h^{\perp}$ if and only if $r_{0}>0$ and $ r_{j} \in \left]-\infty, -r_{0}\right[ \cup \left]0, \infty\right[$
for $j=1,\dots,s$.
\end{lem}
\begin{proof}
%--------------
The basis vectors $\{E_{i}+e_{i}\mid i\in I\}$ of $\Delta\h$ are pairwise orthogonal and have for arbitrary $j=1,..,s$ and $i\in I_{j}$ the squared length $ h_{\textbf{r}}(E_{i}+e_{i},E_{i}+e_{i})=r_{0}+r_{j}$. This proves that $h_{\textbf{r}}$ is non-degenerate on $\Delta\h$ if and only if $r_{j}\neq -r_{0}$ for all $j=1,\dots,s$. Let us compute the orthogonal complement $\Delta\h^{\perp}$ and check for the conditions when $h_{\textbf{r}}$ is positive definite on it: The vectors $E_{i}$ with $ i\in I_{0}$ are orthogonal to $\Delta\h$, so it is enough to compute when two vectors $\sum_{i=1}^{m}a_{i}E_{i}+b_{i}e_{i}\in \g\oplus \h$ and $\sum_{i=1}^{m}c_{i}(E_{i}+e_{i})\in \Delta\h $ are orthogonal:
\begin{align*}
	&h_{\textbf{r}}(\sum_{i=1}^{m}a_{i}E_{i}+b_{i}e_{i},\sum_{i=1}^{m}c_{i}(E_{i}+e_{i}))=\sum_{i=1}^{m}c_{i}(a_{i}h_{\textbf{r}}(E_{i},E_{i})+b_{i}h_{\textbf{r}}(e_{i},e_{i}))\\
	&=\sum_{i=1}^{m}c_{i}(a_{i}r_{0}+b_{i}h_{\textbf{r}}(e_{i},e_{i}))\overset{!}{=}0\quad  \forall c_{i}\quad \Leftrightarrow\quad a_{i}r_{0}={-b_{i}h_{\textbf{r}}(e_{i},e_{i})}.
\end{align*}
Additionally to $E_{i}$ with $i\in I_{0}$ we find the basis vectors $\{ h_{\textbf{r}}(e_{i},e_{i})E_{i}-r_{0}e_{i}\mid i\in I\}$ of $\Delta\h^{\perp}$ which are pairwise orthogonal. Since the squared length of each $E_{i}$ with $i \in I_{0}$ is $r_{0}$, it is necessary that $r_{0} > 0$ for $h_{\textbf{r}}$ to be positive definite. For $i\in I_{j}$ and $j=1,\dots,s$ we compute the squared lengths 
\begin{align*}
	h_{\textbf{r}}(r_{j}E_{j}-r_{0}e_{i},r_{j}E_{j}-r_{0}e_{i})=r_{j}^{2}r_{0}+r_{0}^{2}r_{j}=r_{j}r_{0}(r_{j}+r_{0})
\end{align*}  
which are positive if and only if $r_{0}>0$ and $ r_{j} \in \left]-\infty, -r_{0}\right[ \cup \left]0, \infty\right[$
for $j=1,\dots,s$.
\end{proof}
The enlargement of the group $G$ to $G\times H$ allows to realize the deformed normal homogeneous space as a naturally reductive one. This is a key tool to compute the spectrum of the deformed normal homogeneous space $(G/K,g_{\textbf{t}})$ defined in \eqref{eq: Einführung der Metrik g t}:  
\begin{thm}\label{thm: isometry main result}
%-----------------------------------------------
The isometric action of $G'=G\times H$ on $G$ defined by $(a,h)b=ab\iota(h^{-1})$ induces a $G$-equivariant isometry \begin{align*}
	\Psi:\quad 	(G,g_{\textbf{t}})\rightarrow (\modulo{G'}{\Delta H},h_{\textbf{r}}),\quad a\mapsto (a,e)\Delta H.
\end{align*} 
With respect to $j=1,\dots,s$, the matching of the coefficients is given by $t_{0}=r_{0}$ and $\frac{1}{r_{j}}=\frac{1}{t_{j}}-\frac{1}{t_{0}}$. All coefficients $t_{0}\in \R^{+}$, $t_{j}\in \R^{+} \setminus\{t_{0}\}$ and $r_{0}\in \R^{+}$, $ r_{j} \in \left]-\infty, -r_{0}\right[ \cup \left]0, \infty\right[$ are covered. It restricts to a $G$-equivariant isometry of the quotients 
\begin{align*}
	\Psi:\quad 	(\modulo{G}{K},g_{\textbf{t}})\rightarrow (\modulo{G'}{K'},h_{\textbf{r}})
\end{align*}
and yields a family of naturally reductive realizations $({G'}/{K'},h_{\textbf{r}})$.
\end{thm}
\begin{proof} As $H_{1},\dots,H_{s}\subset G$ commute, the action of $ G'=G\times H $ on $G$ defined by $(g_{1},h)g_{2}=g_{1}g_{2}\tilde{\iota}(h)^{-1}$ is well-defined and it is transitive. The isotropy group $ G_{e}$ is given by $\{(\tilde{\iota}(h),h)\mid h\in H\}=\Delta H$. Hence, $\Psi$ is a diffeomorphism which is $G$-equivariant considering the obvious action of $G$ on $G'/\Delta H$ by the left on the first factor.
The tangent space of $G'/\Delta H$ in the origin $\Psi(e)=(e,e)\Delta H$ is given by $$ \Delta\h^{\perp}=\{h_{\textbf{r}}(e_{i},e_{i})E_{i}-r_{0}e_{i},E_{k}\mid i\in I,\: k\in I_{0}\}. $$ The differential $ d\Psi_{e} $ is the projection $ \g\rightarrow \Delta\h^{\perp}\subset \g\oplus \h $ which is an orthogonal map. For all possible indices $ j\geq 1 $, $i\in I_{j}$, $ k\in I_{0}$, the vectors $ E_{i}/\sqrt{t_{i}} $, $ E_{k}/\sqrt{ t_{0}} $ yield an orthonormal basis of $(\g,g_{\textbf{t}})$ and $(r_{i}E_{i}-r_{0}e_{i})/(\sqrt{r_{i}^{2}r_{0}+r_{i}r_{0}^{2}}),$  $ E_{k}/\sqrt{r_{0}} $ yield an orthonormal basis of $ (\Delta \h^{\perp},h_{r_{1},r_{2}}) $. We see that $ d\Psi:\g\rightarrow \Delta\h^{\perp} $ is an isometry for the wished relations of $ (t_{0},\dots,t_{s})$ and $(r_{0},\dots,r_{s}) $:
\begin{align*}
	&1=h_{\textbf{r}}\left(\frac{E_{i}}{\sqrt{t_{i}}},\frac{r_{i}E_{i}}{\sqrt{r_{i}^{2}r_{0}+r_{i}r_{0}^{2}}}\right)^{2}=\frac{r_{i}^{2}r_{0}^{2}}{{t_{i}(r_{i}^{2}r_{0}+r_{i}r_{0}^{2})}}\quad \Leftrightarrow\quad t_{i}=\frac{r_{i}r_{0}}{r_{i}+r_{0}},\\
	&1=h_{\textbf{r}}\left(\frac{E_{k}}{\sqrt{t_{0}}},\frac{E_{k}}{\sqrt{r_{0}}}\right)^{2}=\frac{r_{0}^{2}}{t_{0}r_{0}},\quad \Leftrightarrow \quad t_{0}=r_{0}.
	%	&1=h_{r_{1},r_{2}}\left(\frac{E_{0}}{\sqrt{3t_{1}}}, \frac{12r_{1}E_{0}-3r_{2}e_{0}}{6\sqrt{3(4r_{1}^{2}r_{2}+r_{1}r_{2}^{2})}}\right)^{2}=\frac{4r_{1}^{2}r_{2}^{2}}{t_{1}(4r_{1}^{2}r_{2}+r_{1}r_{2}^{2})}\quad \Leftrightarrow\quad t_{1}=\frac{4r_{1}r_{2}}{4r_{1}+r_{2}}.
\end{align*}
\end{proof}
In case of Lie groups, such realizations were constructed by D'Atri and Ziller \cite[Thm. 1]{Atri Ziller}. They additionally proved that if $G$ is simple, every left-invariant, naturally reductive metric on $G$ is given by deformations of the biinvariant metric along commuting subalgebras \cite[Thm. 3]{Atri Ziller}.

Being able to realize the deformed homogeneous space as a naturally reductive homogeneous space, Theorem \ref{thm eigenvalues casimir} allows to compute its spectrum.
% \textcolor{red}{Bedingung. an positive Wurzeln und deren max Tori erläutern}. We denote by $R^{+}$ the positive roots of $G$, $R^{+}_{i}$ the positive roots of $H_{i}$ and $\rho=\rho_{0}=\frac{1}{2}\sum_{\lambda\in R^{+}}\lambda$ as well as $\rho_{i}=\frac{1}{2}\sum_{\lambda\in R^{+}_{i}}\lambda$.
\begin{rem}\label{rem: GxH action}
The extension of the transitively acting group $G$ to $G\times H$ leads to the following: As $H$ commutes with $K$, it defines for any $(\varrho_{0},V_{\varrho_{0}})\in \hat{G}_{K}$ a representation on $(V_{\varrho_{0}}^{*})^{K}$ denoted by $\varrho^{H,K}_{0}:=(\varrho_{0}^*\myrestriction H)\myrestriction (V_{\varrho_{0}}^{*})^{K}$. The decomposition of $(V_{\varrho_{0}}^{*})^{K}$ into irreducible $H$-modules correspond precisely to $K'$-spherical representations of $G\times\nolinebreak H$:% We will see in the following proof that $\otimes_{i=1}^{s}\varrho_{i}\in \hat{H}$ occurs in $\varrho^{H,K}_{0}$, i.e. $[\varrho^{H,K}_{0}\mycolon\nolinebreak\otimes_{i=1}^{s}\varrho_{i}]>0$ if and only if $\otimes_{i=0}^{s}\varrho_{i}\in \hat{G'}_{K'}$. 
%	In the next theorem we will see the equivalence for any $\varrho_{0}\in \hat{G}$: \textcolor{red}{achtung def $\Phi_{\varrho}$ wird später relevant!}
%	\begin{align*}
	%		[{\varrho^{K}_{0}}_{H}:\otimes_{i=1}^{s}\varrho_{i}]>0\quad\Leftrightarrow\quad \otimes_{i=1}^{s}\varrho_{i}\in \Phi_{\varrho_{0}}:=\{\otimes_{i=1}^{s}\varrho_{i}\in \hat{H}\:\mid\: \otimes_{i=0}^{s}\varrho_{i}\in \hat{G'}_{K'}\}.
	%	\end{align*}
	\end{rem}
	
	\begin{thm}\label{thm: spectrum main result}
%--------------------------------------------
An irreducible representation $\otimes_{j=0}^{s}\varrho_{j}\in \hat{G'}=\hat{G\times H} $ is $K'$-spherical, i.e. $[\otimes_{j=0}^{s}\varrho_{j}\myrestriction K'\mycolon 1_{K'}]>0$ if and only if $\varrho_{0}\in \hat{G}$ satisfies $[\varrho^{H,K}_{0}\mycolon\otimes_{i=1}^{s}\varrho_{i}]>0$ and $\varrho_{0}$ is $K$-spherical, i.e. $[\varrho_{0}\myrestriction K\mycolon 1_{K}]>0$. The spectrum of $(G/K,g_{\textbf{t}})$ is given by
\begin{align*}
	\spec(G/K,g_{\mathbf{t}})=\left\{\frac{c_{g}(\varrho_{0})-\sum_{i=1}^{s}c_{g\myrestriction\h_{i}}(\varrho_{i})}{t_{0}}+\sum_{i=1}^{s}\frac{c_{g\myrestriction\h_{i}}(\varrho_{i})}{t_{i}}\: \Bigg|  \otimes_{i=0}^{s}\varrho_{i}\in \hat{G'}_{K'}\right\}.
\end{align*}
The scalar $c_{g\myrestriction\h_{i}}(\varrho_{i})=g_{ \h_{i}}(\lambda_{i},2\rho_{i}+\lambda_{i})$ of the $H_{i}$-representation $(\varrho_{i},V_{\varrho_{i}})$ is the eigenvalue of $-d\varrho_{i}\Ca_{g\myrestriction \h_{i}}$, where $\Ca_{g\myrestriction \h_{i}}$ is the generalized Casimir element of $\h_{i}$, see Theorem $\ref{thm eigenvalues casimir}$.
	The multiplicity associated to $\otimes_{i=0}^{s}\varrho_{i}$ is given by $$\mult (\otimes_{i=0}^{s}\varrho_{i}):=\dim(V_{\otimes_{i=0}^{s}\varrho_{i}})\cdot[\otimes_{i=0}^{s}\varrho_{i}\myrestriction K'\mycolon 1_{K' }]=\dim(V_{\otimes_{i=0}^{s}\varrho_{i}})\cdot[\varrho^{H,K}_{0}\mycolon\otimes_{i=1}^{s}\varrho_{i}].$$
\end{thm}

\begin{proof}
	Applying the Peter-Weyl theorem on the different realizations leads to the irreducible decompositions of $L^{2}(G/K)\cong L^{2}(G'/K')$ as a $G$-module. By Theorem \ref{thm: isometry main result}, the diffeomorphism 
	\begin{align*}
		\Psi: L^{2}(G/K)\rightarrow L^{2}(G'/K') 
	\end{align*}
	is $G$-equivariant. The only irreducible $G$-representations occurring in $L^{2}(G/K)$ are $K$-spherical and therefore, the restriction of the $G'$-module $L^{2}(G'/K')$ to $G$ yields also precisely $K$-spherical representations. For any such representation $\varrho_{0}\in \hat{G}_{K}$ we denote by $\Phi_{\varrho_{0}}=\{\otimes_{i=1}^{s}\varrho_{i}\in \hat{H}\:\mid\: \otimes_{i=0}^{s}\varrho_{i}\in \hat{G'}_{K'}\}$ the set of irreducible $H$-representations which extend $\varrho_{0}$ to a $K'$-spherical representation. Then the following $G$-modules are equivalent:
	\begin{align*}
		V_{\varrho_{0}}\otimes (V_{\varrho_{0}}^{*})^{K}\overset{\Psi}{\longrightarrow} \bigoplus_{\otimes_{j=1}^{s}\varrho_{j}\in \Phi_{\varrho_{0}}}		V_{\otimes_{j=0}^{s}\varrho_{j}}\otimes (V_{\otimes_{j=0}^{s}\varrho_{j}}^{*})^{K'}. 
	\end{align*}
	On the right side, the action of $G$ on each factor $	V_{\otimes_{j=0}^{s}\varrho_{j}}$ extends to an $G\times H$-action. By the diagonal embedding of $H$, the action of $(e,h)\in G\times H$ corresponds to the action of $(h^{-1},e)\in G\times H$. Applying $\Psi^{-1}$ and identifying $v\otimes\varphi=f_{v\otimes\varphi}$, leads to an action of $H$ on $(V_{\varrho_{0}}^{*})^{K}$. As $H$ commutes with $K$, its action on $(V_{\varrho_{0}}^{*})^{K}$ is indeed well-defined. Hence we obtain that the given equivalence of $G$-modules extends to an equivalence as $G\times H$-modules, where $G\times H$ acts on $V_{\varrho_{0}}\otimes(V_{\varrho_{0}}^{*})^{K}$ by $ (a,h)v\otimes \varphi\mapsto \varrho_{0}(a)v\otimes\varrho_{0}^{*}(h)\varphi$. %\textcolor{red}{das folgt aus der Identifikation mit den Funktionen und muss noch angemerkt werden!} on the left-hand side and $(a,h)(v,\varphi)\mapsto (\otimes_{j=0}^{s}\varrho_{s}(a,h)v,\varphi)$ on the right-hand side. 
	We conclude that $\otimes_{i=0}^{s}\varrho_{i}$ is $K'$-spherical if and only if $\varrho_{0}$ is $K$-spherical and $[\varrho^{H,K}_{0}:\otimes_{i=1}^{s}\varrho_{i}]>0$. Each representation $\otimes_{i=0}^{s}\varrho_{i}\in \hat{G'}_{K'}$ occurs with multiplicity $$\mult (\otimes_{i=0}^{s}\varrho_{i})=\dim(V_{\otimes_{i=0}^{s}\varrho_{0}})\cdot[\otimes_{i=1}^{s}\varrho_{i}:1_{H}]=\dim(V_{\otimes_{i=0}^{s}\varrho_{i}})\cdot[\varrho^{H,K}_{0}\mycolon\otimes_{i=1}^{s}\varrho_{i}].$$
	We can apply Theorem $\ref{thm: isometry main result}$ and use that the spaces $(G/K,g_{\textbf{t}}) $ and $(\modulo{G'}{K'},h_{\textbf{r}})$ are isometric for the matching coefficients $t_{0}=r_{0}$ and $t_{i}=\frac{r_{i}r_{0}}{r_{i}+r_{0}}$. As the latter one is naturally reductive, we can use Theorem \ref{thm eigenvalues casimir} to compute its spectrum.
	\begin{align*}
		\spec(G'/K',h_{\mathbf{r}})=\left\{\frac{c_{g\myrestriction \g}(\varrho_{0})}{{r_{0}}}+\sum_{i=1}^{s}\frac{c_{g\myrestriction \h_{i}}(\varrho_{i})}{{r_{i}}} \Bigg| \otimes_{i=0}^{s}\varrho_{i}\in \hat{G'}_{K'}\right\}	
	\end{align*}
	Setting $t_{0}=r_{0}$ and $\frac{1}{r_{i}}=\frac{t_{0}-t_{i}}{t_{0}t_{i}}=\frac{1}{t_{i}}-\frac{1}{t_{0}}$ yields the wished result.
\end{proof}
%\begin{samepage}
\begin{rem}
	\begin{enumerate}
		\item[]
		\item	The deformation of the metric $g$ along the subspaces $\h_{i}$ leads to subtracting the eigenvalues of the subgroups $H_{i}$ from the spectrum and then adding them rescaled. The spectrum satisfies the inclusion stated in Corollary \ref{cor: spectrum canonical variation} with respect to the canonical variations of  $\pi_{i}:G/K\rightarrow G/(H_{i}\times K)$. The main problem which representations contribute to the spectrum is solved by computing the $K'$-spherical representations. This can be done using branching rules.   %The summand of the left side yields the spectrum of the horizontal Laplacian, where the right side yields the spectrum of the vertical Laplacian being associated to the Riemannian submersion $G/K\rightarrow G/(H\times K)$, where the whole group $H=H_{1}\times\dots\times H_{s}$ is modded out.
		\item The spectrum of $(G/(H_{i}\times K),g_{\textbf{t}})$ is obtained by the same formula as above, but only those $K'$-spherical representations $\otimes_{i=0}^{s}\varrho_{i}$ contribute which act trivially when restricted to $H_{i}$, i.e. $\varrho_{i}=1_{H_{i}}$.
		\item For any representation $(\varrho,V_{\varrho})\in \hat{G}_{K}$ the operator $-d\varrho( \Ca_{g})\in \End((V_{\varrho}^{*})^{K})$ may have several eigenvalues. The trick of enlarging the group $G$ to $G\times H$ allows to compute the eigenvalues on each irreducible $H$-submodule of $(V_{\varrho}^{*})^{K}$ which correspond precisely to those $H$-representations, which extend $\varrho$ to an $K'$-spherical representation. This is explained in more detail in the example of the Aloff-Wallach manifold in Subsection \ref{subsection: geometric interpretation of the spectrum}.
		\item As the spectrum counted with multiplicities depends continuously on the parameters $t_{0},\dots,t_{s}$ \cite{Bando Urakawa}, the above formula also holds if some parameters $t_{i}$ coincide with $t_{0}$ and the multiplicities add up.
		\item As already mentioned before, in the case of a compact simple Lie group, all naturally reductive metrics are obtained by our construction, see \cite{Atri Ziller}. In this case Theorem \ref{thm: spectrum main result} can be applied, where $H$ is a semisimple subgroup of $G$ and $K=\{e\}$. We recover the formula stated in \cite[p. 994]{Gordon Sutton}: As $K=\{e\}$, $K'=\Delta H$ and $\hat{G}_{K}=\hat{G}$, any representation $\varrho_{1}\in \hat{H}$ which occurs in a given representation $\varrho_{0}\in \hat{G}$ yields a $\Delta H$-spherical representation $\varrho_{0}\otimes\varrho_{1}\in \hat{G\times H}_{\Delta H}$ which contributes to the spectrum. To illustrate this, let us quickly sketch the deformation $(\SU(2),g_{t_{0},t_{1}})$ of the Killing form $g$ of $\SU(2)$ along the subalgebra $\h=\spann_{\R}\{e_{1}:=\text{diag}[i,-i]\}$. Choose $\h$ to be the maximal torus of $\su(2)$ and denote by $\nu_{1}$ the dual of $e_{1}$. It is a fundamental weight of $\su(2)$, i.e. $n_{1}\nu_{1}$ parametrizes for $n_{1}\in \N_{0}$ the highest weight of $\varrho_{0}(n_{1})\in \hat{\SU(2)}$. Given $\varrho_{0}(n_{1})$, the occuring $S^{1}$-weights are $\nu_{1}(n_{1}-2n_{2})$ for $0\leq n_{2}\leq n_{1}$, all with multiplicity one. Theorem \ref{thm: spectrum main result} immediately yields
		\begin{align*}
			\spec(\SU(2),g_{t_{0},t_{1}})=\left\{\frac{n_{1}(2n_{1}+1)-(n_{1}-2n_{2})^{2}}{t_{0}}+\frac{(n_{1}-2n_{2})^{2}}{t_{1}}\:\vline\:n_{1},n_{2}\in \N_{0}\mycolon\: n_{1}\geq n_{2}\right\}.
			%\spec(\SO(3),g_{t_{0},t_{1}})=\left\{\frac{n_{1}(2n_{1}+1)-(n_{1}-2n_{2})^{2}}{t_{0}}+\frac{(n_{1}-2n_{2})^{2}}{t_{1}}\:\vline\:n_{1},n_{2}\in \N_{0}\mycolon\: n_{1}\geq n_{2}\right\}
		\end{align*}
		As irreducible $\SO(3)$-representations are precisely given by $\varrho(n_{1})$ for even $ n_{1}$, the same formula holds for $\spec((\SO(3),g_{t_{0},t_{1}}))$ if $n_{1}$ is restricted to even numbers. This reproduces \cite[Prop. 3.9]{Lauret}.
		%with $\N_{0}$ and the one of $S^{1}$ by $\Z$ and consider $\varrho_{0}(m)\in\hat{SU(2)}$.
	\end{enumerate}
\end{rem}
%\end{samepage}
%
In the next subsection we will use this theorem to compute the spectrum of a family of $3$-contact manifolds that admit a lot of geometric structures, namely to $3$-$(\alpha,\delta)$-Sasaki manifolds.
\subsection{The spectrum of homogeneous $3$-$(\alpha,\delta)$-Sasaki manifolds}\label{subsection: The Spectrum of Sasaki Manifolds}
%----------------------------------------------------------------------------
%
%A almost $3$-contact metric manifold  is defined by its Reeb vector fields $\xi_{i}$, their metrical duals $\eta_{i}=g(\xi_{i},\cdot)$ and the  corresponding fundamental $2$-forms $\Phi_{i}$ ($i=1,2,3$). 
In $2020$, Agricola and Dileo \cite{Agricola Dileo 20} defined new classes of almost $3$-contact metric manifolds which admit metric connections with skew torsion. Generalizing the class of $3$-$\alpha$-Sasaki manifolds and the Heisenberg groups, they discovered the class of $3$-$(\alpha,\delta)$-Sasaki manifolds which are almost $3$-contact manifolds $(M,\varphi_{i},\xi_{i},\eta_{i},g)$ of dimension $4n+3$ defined by the equation:
\begin{align*}
	d\eta_{i}=2\alpha\Phi_{i}+2(\alpha-\delta)\eta_{j}\wedge\eta_{k}.
\end{align*}
%Here, $\eta_{i}:=g(\xi_{i},\cdot)$ are the metrical dual forms of the Reeb vector fields $\xi_{i}$ and $\Phi$ is the fundamental $2$-form defined by the 
Being parameterized by two real parameters $\alpha,\delta$, $\alpha\neq 0$ the class contains two classes of Einstein manifolds: The first class is given by $\alpha=\delta$, corresponding to the whole class of homogeneous $3$-$\alpha$-Sasaki manifolds and the second class is given by $\delta=(2n+3)\alpha$, which is for $n=1$ a class of nearly parallel $G_{2}$-manifolds. Boyer et al. \cite[Thm. C]{Boyer Galicki Mann} proved that every $3$-Sasaki manifold of dimension $4n+3$ is the total space of a locally defined Riemannian submersion over a quaternionic Kähler orbifold with scalar curvature $16n(n+2)$ \cite[Thm. A]{Boyer Galicki Mann}. From this, they classified homogeneous $3$-Sasaki manifolds which are either real projective spaces or simply connected and in $1\mycolon1$-correspondence to complex, simple Lie algebras. Using the theory of their root systems, Goertsches et al. gave a new self-contained proof of this classification \cite{Goertsches Roschig Stecker}. Similar to the work of Boyer et al., Agricola et al. \cite{Agricola Dileo Stecker 21} proved that every $3$-$(\alpha,\delta)$-Sasaki manifold admits a locally defined Riemannian submersion over a quaternionic Kähler manifold of scalar curvature $16n(n+2)\alpha\delta$. Using the work of \cite{Draper}, the authors gave a uniform description of homogeneous $3$-$(\alpha,\delta)$-Sasaki manifolds introducing the concept of generalized $3$-Sasaki data. From this they classified the positive ones, i.e. $\alpha\delta>0$: These are canonical variations of homogeneous $3$-Sasaki manifolds associated to the mentioned Riemannian submersion discovered by Boyer et al.. Additionally, the authors provide a general construction of negative, i.e. $\alpha\delta<0 $ homogeneous $3$-$(\alpha,\delta)$-Sasaki manifolds over non-symmetric Alekseevsky spaces. Negative $3$-$(\alpha,\delta)$-Sasaki manifolds are typically non-compact and therefore not of our interest. We describe positive $3$-$(\alpha,\delta)$-Sasaki manifolds using the canonical variation over Wolf spaces and compute the spectrum on those. They necessarily have fibres  isometric to either $\SU(2)$ or $\SO(3)$.

\medskip	
Let $G$ be compact and simple, $K\subset G$ such that $G/K$ is a homogeneous $3$-Sasaki manifold of dimension $4n+3$. 
Then it is known that the $3$-Sasaki metric is not standard normal homogeneous, but rather 	given by the following rescaling of the Killing form $	g=-\frac{1}{4(n+2)}B$ along the subalgebra $\h=\su(2)\subset \g$:
\begin{align*}
	g_{1/2,1}=\frac{1}{2}g|_{\su(2)^{\perp}}+g|_{\su(2)}.	
\end{align*}
In \cite[Thm. 3.1.1.]{Agricola Dileo Stecker 21}, it is proved that
every positive (i.\,e.\,$\alpha\delta>0$) homogeneous $3$-$(\alpha,\delta)$-Sasaki metric is given by
\begin{align*}
	g_{t_{0},t_{1}}=\frac{1}{2\alpha\delta}g|_{\su(2)^{\perp}}+\frac{1}{\delta^{2}}g|_{\su(2)},
\end{align*}
where the correspondence of parameters is established by setting
$t_{0}=\frac{1}{2\alpha\delta}$ and $t_{1}=\frac{1}{\delta^{2}}$.	
Clearly, $\alpha=\delta=1$ yields the $3$-Sasaki metric.
One of the key results of \cite{Agricola Dileo 20} is the existence of 
a canonical connection with parallel skew torsion on  $3$-$(\alpha,\delta)$-Sasaki manifolds (actually, existence is established even on a much larger class of $3$-contact manifolds, but the torsion is the not necessarily parallel) and the case  $\delta=2\alpha$ singles out
\emph{parallel} $3$-$(\alpha,\delta)$-Sasaki manifolds, i.\,e. those for which all structure tensors are parallel under the canonical connection.
In the homogeneous case, this coincides with the standard normal homogeneous metric.
The subgroup $H$ with Lie algebra $\h$ is either $\SU(2)$ or $\SO(3)$.

The three Reeb vector fields are Killing, generate $\su(2)$ and 
correspond to  trivial summands in the isotropy representation.
Thus, we can  apply Theorem \ref{thm: isometry main result}:\\\\

\begin{cor}
	%------------
	The $3$-$(\alpha,\delta)$-Sasaki manifold $(G/K,g_{t_{0},t_{1}})$ is
	\begin{enumerate}
		\item 
		$G\times H$-normal homogeneous for $t_{1}<t_{0}$, i.e. $2\alpha<\delta$; this includes the nearly parallel $G_{2}$-metrics $\delta=5\alpha$ in case $n=1$.
		\item $G\times H$-naturally reductive for $t_{0}<t_{1}$, i.e. $\delta<2\alpha$ (but not normal homogeneous); this includes the $3$-$\alpha$-Sasaki metrics $\alpha=\delta$.
		\item $G$-normal homogeneous if it is  parallel, i.\,e. $t_{0}=t_{1} \ (\Leftrightarrow 2\alpha=\delta)$.
	\end{enumerate}
\end{cor} 
We can apply Theorem \ref{thm: spectrum main result} 
to derive a formula for its spectrum: 
\begin{cor}\label{cor: spectrum 3 alpha delta} 
	%-------------------------------------------
	The spectrum of the $3$-$(\alpha,\delta)$-Sasaki manifold $(G/K,g_{t_{0},t_{1}})$ with $t_{0}=(2\alpha\delta)^{-1}$, $t_{1}=\delta^{-2}$ is given by
	\begin{align*}
		\spec(G/K,g_{t_{0},t_{1}})=&\{2\alpha\delta (c_{g}(\varrho_{0})-c_{g\myrestriction \su(2)}(\varrho_{1}))+\delta^{2}c_{g\myrestriction \su(2)}(\varrho_{1}) \:\vline\:\varrho_{0}\otimes \varrho_{1 }\in \hat{G\times H}_{K'}\}.	
	\end{align*}
\end{cor}
% Considering a $3$-$(\alpha,\delta)$-Sasaki manifold, isospectral implies isometric when the parameters $\alpha,\delta$ vary. We will see this in detail in the example of the Aloff-Wallach manifold $W^{1,1}$.
%Considering a fixed $3$-$(\alpha,\delta)$-Sasaki manifold,  %Assume that $s=1$. In the set of deformations $\{(G/K,g_{t_{0},t_{1}})\mid t_{0},t_{1}>0\}$, isospectral implies isometric. 
%%
%\begin{prop}
%	Let $G$ be compact and simple, $K\subset G$ such that $G/K$ is a positive $3$-$(\alpha,\delta)$-Sasaki manifold.
%	Choose simple roots of $\g$ which include the maximal root in the corresponding fundamental chamber.
%	If a representation $\varrho\in \hat{G}$ is $K'$-spherical, its highest weight is greater than or equal to the maximal root or greater than or equal to $0$. If the canonical submersion has $\SU(2)$ as fibre, it may also be greater than or equal to the half of the maximal root.   Hierfür: Theorem \cite[thm. 10.1.]{Hall}
%\end{prop}\begin{proof}
%\textcolor{red}{Das ist bisher nur eine Vermutung.}
%\end{proof}
Recently, Nagy and Semmelmann gave an estimate of the first basic eigenvalue $\eta_{1}^{B}\geq\nolinebreak 8(n+\nolinebreak1)$ on $3$-Sasaki manifolds \cite[Cor. 3.10]{Semmelmann Nagy} which yields the estimate $\eta_{1}^{B}\geq\nolinebreak 2\alpha\delta\cdot 4(n+\nolinebreak1)$ on homogeneous positive $3$-$(\alpha,\delta)$-Sasaki manifolds. Using the explicit computations of the Casimir eigenvalues $c_{g}(\varrho) $ for each simple Lie group by Yamaguchi \cite{Yamaguchi}, we present improvements of this estimate for the homogeneous case. We do not attain improvements if $G$ is $\SO(2m+1)$, $\Sp(m)$, $G_{2}$ or $F_{4}$.
\begin{cor}\label{cor: estimates basic eigenvalue}
	Let $(G/K, g_{t_{0}, t_{1}})$ be a positive $3$-$(\alpha, \delta)$-Sasaki manifold of dimension $4n + 3$, where $G$ is a compact simple Lie group. Depending on the type of the simple Lie algebra $\g^{\C}$, the first basic eigenvalue $\eta_{1}^{B}$ 
	has the lower bounds given in Table \ref{Table: lower bounds}. %The horizontal Laplacian $\Delta^{h}$ satisfies on each $\su(2)$-representation space $\varrho(z)$, $z\in \N_{0}$ the lower bounds given in Table \ref{Table: lower bounds2}.
\end{cor}

\begin{center}
	\renewcommand{\arraystretch}{1.3}
	\begin{table}[h!]
		%\rowcolors{2}{gray!5}{white}
		%\begin{tabular}{m{6cm} | m{8cm}}
		%	\rowcolor{gray!5}
		\begin{tabular}{cc}
			\toprule
			$G$ & \text{Lower bound on $\eta_{1}^{B}$} \\
			\midrule
			$\SU(m)$, $\SO(2m)$ & $\eta_{1}^{B} \geq 2\alpha\delta \cdot 4(n + 2)$ \\
			%	\textcolor{red}{diese und SP(m) muss raus!} $\SO(m)$ & $\eta_{1}^{B} \geq 2\alpha\delta \cdot 2\frac{n^{2}+2n}{n-1}$ \\
			$\SO(2m+1)$, $\Sp(m)$ & $\eta_{1}^{B} \geq 2\alpha\delta \cdot 4(n+1)$ \\
			$E_{6}$ & $\eta_{1}^{B} \geq2\alpha\delta\cdot48$\\
			$E_{7}$ & $\eta_{1}^{B} \geq2\alpha\delta\cdot72$\\
			$E_{8}$ & $ \eta_{1}^{B} \geq2\alpha\delta\cdot120$\\
			${F_4}$ & $\eta_{1}^{B} \geq 2\alpha\delta \cdot 32$ \\
			${G_2}$ & $\eta_{1}^{B} \geq 2\alpha\delta \cdot 12$ \\
			\bottomrule
		\end{tabular}
		\caption{Lower bounds for the first basic eigenvalue $\eta_{1}^{B}$ on positive $3$-$(\alpha, \delta)$-Sasaki manifolds.}
		\label{Table: lower bounds}
	\end{table}
\end{center}	
\begin{proof}
	If $G/K$ is a $3$-$(\alpha,\delta)$-Sasaki manifold, the maximal torus $\t\subset \g$ splits orthogonally into $\t_{\k}\oplus\t_{\su(2)}$, where $\t_{\k}$ and $\t_{\su(2)}$ denote the maximal torus of $\k$ and $\su(2)$ respectively. We can thus think about any weight of $\t$ as the sum of a weight of $\t_{\k}$ and a weight of $\t_{\su(2)}$. Let us fix simple roots such that the positive root of $\t_{\su(2)}$, which is the maximal root, is contained in the fundamental Weyl chamber, see \cite{Goertsches Roschig Stecker}. Now, assume that a $K$-spherical representation $\varrho\in \hat{G}_{K}$ produces a basic eigenvalue, i.e. $\varrho\otimes1_{\su(2)} $ is $K'$-spherical. This means that the restriction of $\varrho$ to $K$ acts on at least one vector $v\in V_{\varrho}$ trivially
	and that the dual representation of $\varrho$ restricted to $\su(2)$ acts on $v$ trivially too, see Theorem \ref{thm: spectrum main result}. As any dual representation of $\su(2)$ is self-dual, this implies that the weight corresponding to $v$ is trivial. This means that either $\varrho$ is the trivial representation, or that its highest weight is higher or equal to a root which is contained in the fundamental Weyl chamber. In the cases of $\SU(m)$, $\SO(2m)$, $E_{6}$, $E_{7}$, $E_{8}$ the root system is simply laced, i.e. any two roots have the same length and the Weyl group acts transitively on the set of roots. As there is only one Weyl group orbit, we can apply \cite[Prop. 3]{Moody Patera} and conclude that only one root is contained in a fundamental Weyl chamber which can be identified with the selected root of $\t_{\su(2)}$. With respect to the Killing form, this (maximal) root has Casimir eigenvalue $1$, see \cite{Yamaguchi}. This translates to the $3$-$(\alpha,\delta)$-Sasaki metric as $2\alpha\delta \cdot 4(n + 2)$. As the Casimir eigenvalue increases if a positive root is added, we obtain the estimate $c_{\g}(\varrho)\geq 2\alpha\delta\cdot 4(n+2)$. 
	%		 We will prove that this estimate is sharp. To do so, it is already enough to prove that $\varrho \in \hat{G}$ which has highest weight $\mu$ is $K$-spherical.  by computing the multiplicities $m(\mu_{1})$, $m(0)$ of the weights $\mu_{1}$ and $0$ using the modified Freudenthal formula described in \cite{Moody Patera}. It is a recursive formula which makes the computation of the multiplicity of the highest weight possible by first computing the multiplicities of the occuring lower weights. In the case of the basic eigenvalue, there are only two weights in the fundamental Weyl chamber namely the trivial weight $0$ and the maximal root $\mu_{1}$ which makes the application in this case particularly easy. As the multiplicity of the weights is invariant of the Weyl group action, we may change the positive roots such that the maximal root is actually a root of $\t_{\k}$. We will again compute the relation of the multiplicity of the trivial weight $0$ and the one of the maximal root $\mu$. The comparison of both relations will tell us how often the trivial weight $0$ corresponds to the trivial representation of $\t_{\k}$.  The group $W^{T}$ generated by the Weyl group $W$ and $-\id$ acts in in the cases $A_{l}$, $D_{l}$, $E_{6}$, $E_{7}$, $E_{8}$ transitively on the set of roots, i.e. its orbit equals the set of roots. The unique root of each $W^{T}$ orbit which is in the fundamental Weyl chamber \cite[Prop. 3]{Moody Patera} is the maximal root $\mu_{1}$.  In the case of $\SU(m)$ we obtain the    \\
	%		 \\
	In the missing cases, there is besides the maximal root an other positive root in the fundamental Weyl chamber. The minimum of their Casimir eigenvalues is less than the lower estimate $2\alpha\delta\cdot 4(n+1)$ provided by \cite{Semmelmann Nagy}.
\end{proof}
\begin{rem}\label{rem: horizontal eigenvalue estimates for each z}
	\begin{itemize}
		\item[]
		\item The estimate $2\alpha\delta\cdot 4(n+1)$ for the first basic eigenvalue is not achieved if $M$ is not a sphere, see \cite[Cor. 3.10]{Semmelmann Nagy}. If $G=\SU(m)$ and $n=1$, i.e. $G/K$ is the Aloff-Wallach manifold $W^{1,1}$, our improved estimate $2\alpha\delta \cdot 4(n+2) =2\alpha\delta \cdot 12$ is achieved, see Theorem \ref{thm: spectrum Aloff Wallach}.
		\item Let $\mu_{1}$ denote the maximal root in the fundamental Weyl chamber which corresponds to the $\su(2)$ root in $\g$ and let $z\in \N_{0}$. If one fixes an $\SU(2)$-representation $\varrho_{1}(z)$ with highest weight $\frac{z}{2}\mu_{1}$ (for even $z$ this can be identified with an $\SO(3)$-representation) and investigates the horizontal Laplacian on the function space on which this representation is realized, one can proceed analogously to the last proof to obtain eigenvalue estimates.  The Casimir eigenvalue of $\frac{z}{2}\mu_{1}$ is in all cases given by $2\alpha\delta\cdot4 (z^{2}/4+z(n+1)/2) $. Subtracting the fibre eigenvalue $2\alpha\delta\cdot z(z/2 +2)/2$ yields the lower bound $2\alpha\delta\cdot zn$ in all cases for the horizontal Laplacian. This estimate equals the estimate stated in \cite[Prop. 3.4]{Semmelmann Nagy}.% and the estimate is achieved for any $z\in 2\N_{0}$ for $G\neq \Sp(n)$ and $z\in \N_{0}$ for $G=\Sp(n)$
	\end{itemize}		
\end{rem}

The main problem for deducing from Corollary \ref{cor: spectrum 3 alpha delta} the full spectrum is to compute the spherical representations explicitly, which can be quite challenging.
The second part of this paper is solely devoted to the 
$7$-dimensional Aloff-Wallach manifold $W^{1,1}$. 
We remark that building upon the geometric framework developed here, a systematic representation-theoretic resolution of this branching problem for the entire classical series of homogeneous 3-$(\alpha,\delta)$-Sasaki manifolds (Types A, B, C, and D) has recently been achieved in \cite{AgricolaCaglieroHenkel26}.

%
%--------------------------------------------------------------------------------------------
\section{\texorpdfstring{The spectrum of the Aloff-Wallach manifold $W^{1,1}=\SU(3)/S^{1}$}{The spectrum of the Aloff-Wallach manifold W1,1 = SU(3)/S1}}\label{section: Aloff-Wallach manifold}
%-------------------------------------------------------------------------------------------
%
The Aloff-Wallach manifolds $W^{k,l}=\SU(3)/S^{1}_{k,l}$ depend on the coprime embedding parameters $k,l\in \Z$ of $S^{1}_{k,l}$ into $\SU(3)$:
\begin{align*}
	S^1_{k,l}\rightarrow \SU(3), \quad 
	z\mapsto
	\begin{pmatrix}
		z^{k} & 0 & 0 \\
		0 & z^{l} & 0 \\
		0 & 0 & z^{-k-l}
	\end{pmatrix}.
\end{align*} They are a well-known class of manifolds and realize a $1$-parameter family of metrics with strictly positive sectional curvature \cite{Aloff Wallach 75}. Since their discovery by Aloff and Wallach in 1975, they attracted a lot of attention {\cite[Exa. 4.5]{BFGK}, \cite[Prop. 9.17]{Boyer Galicki Mann}}. Each Aloff-Wallach manifold $W^{k,l}$ has precisely two homogeneous Einstein metrics \cite{Nikonorov 04,Page Pope 84,Wang 82} and admits proper nearly parallel $G_{2}$-structures \cite{Friedrich et al. 97}. Those make them very interesting for superstring theory \cite{Agricola 06} and theoretical physics: The corresponding $G_{2}$-instantons have been classified in \cite{Ball Oliveira 19}. In the special case $k=l=1$, they additionally admit a $3$-$(\alpha,\delta)$-Sasaki structure \cite[Exa. 3.3.1]{Agricola Dileo Stecker 21}. The fact that they are also $7$-dimensional is crucial for their role in $M$-theory and the Ad\nolinebreak S/CFT correspondence \cite{Billo et al. 01}.\\

Urakawa \cite{Urakawa} computed the spectrum on any of those spaces with respect to the Killing form of $\SU(3)$ which is not contained in the $1$-parameter family of positive metrics. Moreover, basing on the classification of Berger \cite{Berger}, Urakawa \cite{Urakawa list} gave a list of the first eigenvalue of any simply connected, normal homogeneous space with positive curvature. This list is incomplete as Wilking \cite{Wilking} found a gap in Berger's classification: He was able to realize the positively curved metrics of the Aloff-Wallach manifold $W^{1,1}$ as a $1$-parameter family of $\SU(3)\times\SO(3)$-normal homogeneous metrics and he named this new realization $V_{3}$. Using the methods developed in Subsections \ref{subsection naturally reductive realizations}, we extend this normal homogeneous realization to a naturally reductive one which gives access to the spectrum of all $3$-$(\alpha,\delta)$-Sasaki metrics on $W^{1,1}$, see Theorem \ref{thm: spectrum Aloff Wallach}. By doing so, we complete the list of the first eigenvalue of a compact simply connected normal homogeneous space with positive curvature in Corollary \ref{cor: list of first eigenvalue on compact simply connected normal hom pos curvature}. Moreover, Urakawa's $\SU(3)$-normal homogeneous computation of the spectrum (this metric is not of positive curvature) is recovered as a limiting case. We relate the spectrum to the geometry and existing geometric structures in Subsection \ref{subsection: geometric interpretation of the spectrum}, use the results to check eigenvalue estimates for canonical variations of $3$-Sasaki manifolds in Subsection \ref{subsection: interpretation as canonical variation} and study the spectrum for fixed volume in Subsection \ref{subsection: spectrum with constant volume}. We will prove that an analogous inequality to Hersch \cite{Hersch} does not exist on $3$-$(\alpha,\delta)$-Sasaki manifolds, by constructing a $1$-parameter family of metrics for which the product of the first eigenvalue with the volume is neither lower nor upper bounded:
\begin{align*}
	\lim_{t_{1}\rightarrow \infty}	\eta_{1}(g_{t_{1}^{-3/4},t_{1}})\cdot \vol(W^{1,1},g_{t_{1}^{-3/4},t_{1}})^{2/7}=\infty,\quad \lim_{t_{1}\rightarrow 0}	\eta_{1}(g_{t_{1}^{-3/4},t_{1}})\cdot \vol(W^{1,1},g_{t_{1}^{-3/4},t_{1}})^{2/7}=0.
\end{align*}

In order to describe the Aloff-Wallach manifold $W^{1,1}$ more closely, we denote by $$\iota:\nolinebreak\U(2)\rightarrow\nolinebreak\SU(3)$$ the upper left embedding. In this case $S^{1}_{1,1}$ corresponds to the embedded center of $\U(2)$ into $\SU(3)$ and we use the short notation $S^{1}$.
Considering the rescaled Killing form $g(XY)=-\frac{1}{2}\tr(XY)$ the $3$-$(\alpha,\delta)$-Sasaki metric $g_{t_{0},t_{1}}$ of $W^{1,1}=\SU(3)/S^{1}$ is given by
\begin{align*}
	g_{t_{0},t_{1}}=t_{0}g\vert_{\su(2)^{\perp}}+t_{1}g\vert_{\su(2)}	=\frac{1}{2\alpha\delta}g\vert_{\su(2)^{\perp}}+\frac{1}{\delta^{2}}g\vert_{\su(2)}.	
\end{align*}
Continuing Example \ref{exa: incomplete aloff wallach}, we introduce the subgroup ${S^{1}}'$of $\SU(3)\times \SU(2)$ $${S^{1}}'=\{(S\iota(A),A)\mid S\in {S^{1}},\: A\in \SU(2)\}$$ and the following subgroup of $\SU(3)\times \SO(3)$:
\begin{align*}
	\Uni^{\bullet}(2)=(\iota,\pi)(\Uni(2)),\quad \pi\colon  \Uni(2)\rightarrow \modulo{\Uni(2)}{S^{1}} \cong \SO(3).
\end{align*} 
We obtain the space % $V_{3}$ introduced by Wilking \cite{Wilking}:
\begin{align*}
	V_{3}:=\modulo{(\SU(3)\times \SU(2))}{{S^{1}}'}\cong	\modulo{(\SU(3)\times \SO(3))}{\Uni^{\bullet} (2)},
\end{align*}
which is equipped with the restricted semi-Riemannian metrics of $\SU(3)\times \SO(3)$
\begin{align*}
	h_{r_{0},r_{1}}=r_{0}g|_{\su(3)}+r_{1}g|_{\su(2)},\quad r_{0}>0,\quad  r_{1} \in \left]-\infty, -r_{0}\right[ \cup \left]0, \infty\right[.
\end{align*}
By Lemma \ref{lem: pos definite bilinear form}, $h_{r_{0},r_{1}}$ restricts to a positive definite scalar product on $\u^{\bullet}(2)^{\perp}\subset \su(3)\times\so(3)$ which extends to a Riemannian metric on $V_{3}$.
The space $(V_{3},h_{r_{0},r_{1}})$ was introduced by Wilking for positive $r_{0},r_{1}$. It denotes the missing third exceptional space in Berger's classification:
\begin{thm} \cite[Thm. 2.4, Lem 3.1]{Aloff Wallach 75}\cite{Wilking}
	%-------------------------------------------------------------------
	The $\SU(3)\times\SO(3)$-normal homogeneous space $(V_{3},h_{r_{0},r_{1}})$, $r_{0},r_{1}>0$ is isometric to the non $\SU(3)$-normal homogeneous space $(W^{1,1},g_{t_{0},t_{1}})$ with $0<t_{1}<t_{0}$. In particular, it has positive sectional curvature. All metrics are covered. 
\end{thm}
For the explicit parameter matching for the metrics, we refer to the formulas stated in Corollary \ref{cor: naturally reductive realization aloff wallach}.
In Theorem \ref{thm: isometry main result} we proved abstractly that this $\SU(3)\times\SO(3)$-normal homogeneous realization can be extended to a naturally reductive one:
%Identifying $\su(2)\cong\so(3)$, the space $V_{3}$ is a naturally reductive realization with respect to the reductive complement $\u^{\bullet}(2)$ and with respect to the metric \textcolor{red}{andere Notation der Metrikeinschränkung?}\marginpar{remark auskommentiert $\su(2)$ vs $\so(3)$}

\begin{cor}\label{cor: naturally reductive realization aloff wallach}
	The naturally reductive space $(V_{3},{h}_{r_{0},r_{1}})$ is isometric to $(W^{1,1},{g}_{t_{0},t_{1}})$ for $t_{0}=r_{0}$ and $t_{1}=\frac{r_{0}r_{1}}{r_{0}+r_{1}}>0 $, i.e. $r_1 = \frac{t_1 t_0}{ t_0 - t_1}$. All metrics on $(W^{1,1},{g}_{t_{0},t_{1}})$ with $t_{0}>0$, $t_{1}\in \R^{+}\setminus\{t_{0}\}$ and all metrics on $(V_{3},{h}_{r_{1},r_{2}})$, i.e. $r_{0}>0$ and $ r_{1} \in \left]-\infty, -r_{0}\right[ \cup \left]0, \infty\right[$ are covered.
\end{cor}
The naturally reductive realizations $(V_{3},h_{r_{0},r_{1}})$ include the nearly parallel $G_{2}$ metric:
	
	\begin{rem}
		The Killing form of $\SU(3)\times \SU(2)$ corresponds to $r_{0}=12$, $r_{1}=8$ which translates to $t_{0}=12=\frac{1}{2\alpha\delta}$ and $t_{1}=4.8=\frac{1}{\delta^{2}}$. As this implies $\delta=5\alpha$, this yields the nearly parallel $G_{2}$ Einstein metric on $W^{1,1}$ (see \cite{Agricola Dileo 20} for details).
	\end{rem}
	
	The fundamentals of representation theory of $\SU(3)\times\SO(3)$ are required to describe the spectrum.
	\subsection{\texorpdfstring{Representation theory of $\SU(3)\times \SO(3)$}{Representation theory of SU(3)x SO(3)}} We introduce an orthogonal basis for $ \su(3)\times\u(2) $ which is except for $||E_{0}||^{2}=3$ normalized with respect to $g_{1,1}$:  %The Lie algebra $\so(3)\times \su(3)$ is isomorphic to $\su(2)\times \su(3)$. 
	\begin{align*}
		&	e_{0}=\begin{bmatrix}
			i&0\\
			0&i
		\end{bmatrix},\quad e_{1}=\begin{bmatrix}
			i&0\\
			0&-i
		\end{bmatrix},\quad e_{2}=\begin{bmatrix}
			0&-1\\
			1&0
		\end{bmatrix},\quad e_{3}=\begin{bmatrix}
			0&i\\
			i&0
		\end{bmatrix},\\
		&E_{j}=d\iota e_{j},\: j=0,\dots,3, \quad\text{where } d\iota(A)=\begin{bmatrix}
			A&0\\
			0&-\tr(A)
		\end{bmatrix},\\
		&E_{4}=A_{13},\quad E_{5}=\tilde{A}_{13},\quad E_{6}=A_{23},\quad E_{7}=\tilde{A}_{23},\\
		&\text{where}\quad  A_{ij}=(\delta_{ij}-\delta_{ji})_{i,j},\quad  \tilde{A}_{ij}=(i\delta_{ij}+i\delta_{ji})_{i,j}.
	\end{align*}
	The element $E_{0}$ is not normalized as it yields a basis vector of the maximal torus and needs to be compatible with the natural projections $\lambda_{i}$:
	\begin{stand} With respect to the chosen basis we define
		\begin{enumerate}
			\item Maximal tori of $\su(3)$, $\su(2)$ and $\su(3)\oplus\su(2)$:
			\begin{align*}
				\t_{1}=\spann_{\R}\{E_{1},{E}_{0}\}\subset \su(3),\quad \t_{2}=\spann_{\R}\{e_{1}\}\subset \su(2),\quad  \t=\t_{1}\oplus \t_{2}\subset \su(3)\oplus \su(2) 
			\end{align*}
			and the maximal torus of the subalgebra $\u^{\bullet}(2)$: $$\t_{\u^{\bullet}(2)}=\spann_{\R}\left\{{E}_{0},E_{1}+e_{1}\right\}\subset \u^{\bullet}(2).$$
			\item With respect to $ A=a_{0}{E}_{0}+a_{1}E_{1}+b_{1}e_{1}\in \t $ and for $ i=1,2,3 $, $ j=1,2 $ the projections onto the diagonal entries $ \lambda_{i},\nu_{j},\mu_{j}\in i\t^{*}  $ by \begin{align*}
				\lambda_{1}(A)=i(a_{1}+ a_{0}),\quad \lambda_{2}(A)=i(a_{0}-a_{1}),\quad \lambda_{3}(A)=-2ia_{0},\quad \nu_{j}(A)=(-1)^{j+1}ib_{1}.
			\end{align*}
			%$ \lambda_{j}((a_{ij})_{ij},(b_{jl})_{jl})=ia_{jj} $ and $ \nu_{j}((a_{ij})_{ij},(b_{jl})_{jl})=ib_{jj} $. 
			Note that $ \{\lambda_{1},\lambda_{2},\nu_{1}\} $ yields a basis for $ i\t^{*}$ and define $ \mu_{j}=2\nu_{j}$.  %We write $ \mu_{1}=2\nu_{1} $ and $ \mu_{2}=2\nu_{2}. $
		\end{enumerate}
	\end{stand}
	We will need the scalar products of the dual vectors in several situations.
	\begin{lem}
		For $ i,j\in \{1,2\},\: i\neq j $ the length of the dual vectors is given by% $ i\t^{*}\subset i(\su(2)\oplus\su(3))^{*} $ is given by
		\begin{align*}
			h_{1,1}(\lambda_{j},\mu_{1}) &= 0, & h_{1,1}(\mu_{1},\mu_{1}) &=4,\\
			h_{1,1}(\lambda_{j},\lambda_{j}) &= \frac{4}{3}, & h_{1,1}(\lambda_{i},\lambda_{j}) &= -\frac{2}{3}.
		\end{align*}
	\end{lem}
	\begin{proof}
		The vectors $X_{\lambda}$ satisfying $\lambda=ih_{1,1}(X_{\lambda},\cdot)$ are given by:
		\begin{align*}
			X_{\lambda_{j}}=\frac{{E}_{0}}{3}+(-1)^{j+1}E_{1},\quad X_{\mu_{1}}=2e_{1}.
		\end{align*}
		Taking the respective scalar products yields the result.
		%	The dual vectors are supposed to satisfy \begin{align*}
			%		\lambda_{j}(a_{0}{E}_{0}+a_{1}E_{1}+b_{1}e_{1})&=ih_{1,1}(X_{\lambda_{j}},a_{0}{E}_{0}+a_{1}E_{1}+b_{1}e_{1})=i(a_{0}+(-1)^{j+1}a_{1})\\
			%		\mu_{1}(a_{0}{E}_{0}+a_{1}E_{1}+b_{1}e_{1})&=ih_{1,1}(X_{\mu_{1}},a_{0}{E}_{0}+a_{1}E_{1}+b_{1}e_{1})=2ib_{1}
			%	\end{align*}
		%	as the vectors are orthogonal with length 
		%	\begin{align*}
			%		h_{1,1}({E}_{0},{E}_{0})=3,\quad h_{1,1}(E_{1},E_{1})=1,\quad h_{1,1}(e_{1},e_{1})=1	
			%	\end{align*}
		%	We get: \begin{align*}
			%		X_{\lambda_{j}}&=\frac{{E}_{0}}{3}+(-1)^{j+1}E_{1},\quad X_{\mu_{1}}=2e_{1}\\
			%		h_{1,1}(\lambda_{j},\lambda_{j})&=\frac{1}{3}+1=\frac{4}{3},\quad h_{1,1}(\lambda_{j},\lambda_{i})=\frac{1}{3}-1=-\frac{2}{3}\quad h_{1,1}(\mu_{1},\mu_{1})= 4
			%	\end{align*}
	\end{proof}
	
	% We denote by $\t$ the Lie algebra of the diagonal maximal torus in $ \SU(3)\times \SO(3) $ and $\t_{\u^{\bullet}(2)}$ the one of the diagonal maximal torus in $\U^{\bullet}(2)$. They are given by:
	%%\begin{align*}
	%%	\t=\left\{\left(\begin{bmatrix}
		%	%		ia_{1}&0&0\\
		%	%		0&ia_{2}&0\\
		%	%		0&0&-i(a_{1}+a_{2})
		%	%	\end{bmatrix},\begin{bmatrix}
		%	%	ib&0\\
		%	%	0&-ib
		%	%\end{bmatrix}\right)\:\vline\: a_{1},a_{2},b\in \R\right\}.
		%	%\end{align*}
		%	\begin{align*}
			%		\t=\spann_{\R}\left\{{E}_{0},E_{1},e_{1}\right\},\quad 	.
			%	\end{align*}
		The root systems of these Lie algebras is a standard result in Lie theory:% In order to compute the spectrum, we give the root system with respect to the chosen torus $\t$.
		\begin{prop}[{\cite[p. 981]{Urakawa},\cite{Hall}}]\label{prop: D(SU(3)xSO(3))}% The root system of $ \su(3) $ is given by $ R_{1}=\left\{\pm(\lambda_{1}-\lambda_{2}),\pm(\lambda_{1}-\lambda_{3}),\pm (\lambda_{2}-\lambda_{3})\right\} $ and the root system of $ \su(2) $ is given by $ R_{2}=\{\mu_{1},\mu_{2}\} $ . Hence T
			The roots of $ \su(3)\oplus \su(2) $ are		
			\begin{align*}
				R=\{\pm(\lambda_{1}-\lambda_{2}),\pm (\lambda_{1}-\lambda_{3}),\pm(\lambda_{2}-\lambda_{3}),\mu_{1},\mu_{2}\}.
			\end{align*} We choose the following basis and the corresponding positive roots: \begin{align*}
				\Delta=	\{\lambda_{1}-\lambda_{2},\lambda_{2}-\lambda_{3},\mu_{1}\},\quad R^{+}=\{\lambda_{1}-\lambda_{2},\lambda_{1}-\lambda_{3},\lambda_{2}-\lambda_{3},\mu_{1}\}.
			\end{align*}  The half of the sum of positive roots is given by
			\begin{align*}
				\rho=\lambda_{1}-\lambda_{3}+\frac{\mu_{1}}{2}=2\lambda_{1}+\lambda_{2}+\nu_{1}.
			\end{align*}
			The dominant integral weights
			\begin{align*}
				D(\SU(3)\times \SO(3))=\{z_{1}\lambda_{1}+z_{2}\lambda_{2}+{z_{3}}\mu_{1}\mid z_{i}\in \N_{0}: z_{1}\geq z_{2}\geq 0 ,\: z_{3}\geq 0\}
			\end{align*}
			parametrize the unitary irreducible representations of $\SU(3)\times \SO(3)$ by
			\begin{align*}
				\varrho(z_{1},z_{2},z_{3})=\varrho_{0}(z_{1},z_{2})\otimes \varrho_{1}(z_{3}), \text{ where } z_{1}\geq z_{2}.
			\end{align*}
			Their dimension can be computed by
			\begin{align*}
				\dim(V_{\varrho(z_{1},z_{2},z_{3})})=\dim(V_{\varrho_{0}(z_{1},z_{2})})\cdot\dim(V_{\varrho_{1}(z_{3})})=\frac{(z_{1}-z_{2}+1)(z_{1}+2)(z_{2}+1)}{2}\cdot (2z_{3}+1).
			\end{align*}
			%		\begin{align*}
				%			\dim(V_{\varrho(z_{1},z_{2},z_{3})})=\frac{(z_{1}-z_{2}+1)(z_{1}+2)(z_{2}+1)(2z_{3}+1)}{2}.
				%		\end{align*}
			%If we replace $\mu_{1}$ with $\nu_{1}=\mu_{1}/2$, we obtain the dominant integral elements of$\SU(3)\times\SU(2)$.		
		\end{prop}
		With this background, the results can be specified.

		\subsection{Formulas for the spectrum of $(W^{1,1},g_{t_{0},t_{1}})$}
		%---------------------------------------------------------------------------------
		%
		Since the highest weight of any unitary irreducible representation of ${\SU(3)\times\SO(3)}$ is parametrized by $z_{1}\lambda_{1}+z_{2}\lambda_{2}+z_{3}\mu_{1}$ with $z_{1},z_{2},z_{3}\in \N_{0}$ and $z_{1}\geq z_{2}$, we use the shorthand notation
		$$m(z_{1},z_{2},z_{3}):=[\varrho(z_{1},z_{2},z_{3})\myrestriction \U^{\bullet}(2)\mycolon1_{\U^{\bullet}(2)}].$$
		Furthermore, we define the following purely combinatorial partition function
		\begin{align*}
			\wp(a_{1},a_{2})&=\# \{(m_{1},m_{2},m_{3})\in \N_{0}^{3}\mid a_{1}=3(m_{1}+m_{2}),\quad a_{2}=m_{3}-m_{1}-2m_{2}\}\\
			&=\begin{cases}
				0,\quad \text{if } \frac{a_{1}}{3}\notin \N_{0}\text{ or }a_{2}\notin \Z\\
				\max(0, \frac{a_{1}}{3}+1+\min(0,a_{2}+\frac{a_{1}}{3})),\quad \text{else}.
			\end{cases}
		\end{align*}
		%We are now in a position to give the spectrum of the Aloff-Wallach space $W^{1,1}$:
		We split the description of the spectrum $\spec(\SU(3)/S^{1}_{1,1},g_{t_{0},t_{1}})$ in three parts, which will be proved separately in the following pages.
		\begin{thm}\label{thm: spectrum Aloff Wallach}
			
			\begin{enumerate}
				\item[]
				\item \label{item: spectrum formula}The spectrum $\spec(\SU(3)/S^{1},g_{t_{0},t_{1}})$ is given for $t_{0}=(2\alpha\delta)^{-1}$, $t_{1}=\delta^{-2}$ by all real numbers
				\begin{align*}
					\eta(z_{1},z_{2},z_{3})=\:2\alpha\delta \left(\frac{4(z_{1}^{2}+z_{2}^{2}-z_{1}(z_{2}-3))}{3}-4z_{3}(z_{3}+1)\right)+\delta^{2}4z_{3}(z_{3}+1)
				\end{align*}
				for which $m(z_{1},z_{2},z_{3})>0$.	The multiplicity associated to $\varrho(z_{1},z_{2},z_{3})$ is given by
				\begin{align*}
					\mult(z_{1},z_{2},z_{3})=	m(z_{1},z_{2},z_{3})\cdot\frac{(z_{1}-z_{2}+1)(z_{1}+2)(z_{2}+1)(2z_{3}+1)}{2}.
				\end{align*}
				\item \label{item: m formula}$m(z_{1},z_{2},z_{3})$ can be computed by
				\begin{align*}
					m(z_{1},z_{2},z_{3})=&	\wp({z_{1}+z_{2},-z_{3}-z_{1}-2})+\wp({z_{1}+z_{2},z_{3}-z_{2}})\\
					&-\wp({z_{1}+z_{2},-z_{3}-z_{2}-1})			-\wp({z_{1}+z_{2},z_{3}-z_{1}-1})\\
					&	+ \wp({z_{1}-2z_{2}-3,z_{2}+1-z_{3}})	
					+\wp({z_{1}-2z_{2}-3,z_{2}-z_{1}+z_{3}})\\
					& -\wp({z_{1}-2z_{2}-3,z_{2}+2+z_{3}}).
				\end{align*} 
				The representation $\varrho(2,1,0)$ produces the first basic eigenvalue $\eta_{1}^{B}=2\alpha\delta\cdot 12$ which is at the same time the first eigenvalue if $g_{t_{0},t_{1}}$ has positive sectional curvature. Otherwise $\varrho(2,1,1)$ yields the first eigenvalue: \[
				\eta_{1} = 
				\begin{cases}
					2\alpha\delta\cdot 4 +\delta^{2}\cdot 8,\quad   &2\alpha\delta \geq \delta^{2}, \\
					2\alpha\delta\cdot 12	,\quad  &\delta^2> 2\alpha\delta.
				\end{cases}
				\]
				In the class of $3$-$(\alpha,\delta)$-metrics, isospectral implies isometric as $\eta(g_{0.5,t_{1}})=8(1+\nolinebreak t_{1}^{-1})$ is always the first non-constant eigenvalue.
				\item\label{item: S1 spherical} $\varrho_{0}(z_{1},z_{2})\in \hat{\SU(3)}_{S^{1}_{1,1}}$ if and only if $z_{1}+z_{2}\equiv 0\pmod{3}$, hence these are precisely the tuples $(z_{1},z_{2})$ which occur in the spectrum.
			\end{enumerate}% which satisfy $z_{1}+z_{2}\equiv 0 \pmod{3}$. 
			%	\marginpar{auskommentiert $SU(2)$ vs $SO(3)$}
			%	\begin{align*}
				%	\spec(\SU(3)/S^{1}_{1,1},g_{t_{0},t_{1}})=&\left\{{2\alpha\delta (\frac{4(z_{1}^{2}+z_{2}^{2}-z_{1}(z_{2}-3))}{3}-z_{3}(z_{3}+2))}+\delta^{2}z_{3}(z_{3}+2)\right.\\
				%	& \:\: \left. \vline\: z_{1}\lambda_{1}+z_{2}\lambda_{2}+z_{3}\nu_{1}\in D(\SU(3)\times \SU(2),(S^{1}_{1,1})')\right\}\\
				%	=&\left\{{2\alpha\delta (\frac{4(z_{1}^{2}+z_{2}^{2}-z_{1}(z_{2}-3))}{3}-4z_{3}(z_{3}+1))}+\delta^{2}4z_{3}(z_{3}+1)\right.\\
				%	& \:\: \left. \vline\: z_{1}\lambda_{1}+z_{2}\lambda_{2}+z_{3}\mu_{1}=:\lambda\in D(\SU(3)\times \SO(3),U^{\bullet}(2))\right\}	
				%\end{align*}
			\end{thm}
			For an evaluation of this Theorem, please refer to Subsections \ref{subsection: geometric interpretation of the spectrum} -- \ref{subsection: spectrum with constant volume}. For a table of $\U^{\bullet}(2)$-spherical representations, multiplicities and eigenvalues see the electronic supplement \cite{Agricola Henkel 24} or Table \ref{table: first eigenvalues, reps and mult}.
			\subsubsection{Proof of \eqref{item: spectrum formula} in Theorem $\ref{thm: spectrum Aloff Wallach}$}
			
			We can apply Corollary \ref{cor: spectrum 3 alpha delta} to compute the spectrum:
			\begin{align*}
				&\spec(\SU(3)/S^{1}_{1,1},g_{t_{0},t_{1}})\\
				&=\left\{2\alpha\delta (c_{g}(\varrho_{0})-c_{g\myrestriction \su(2)}(\varrho_{1}))+\delta^{2}c_{g\myrestriction \su(2)}(\varrho_{1})\right.\:\: \left. \vline\:\varrho_{0}\otimes \varrho_{1}\in \hat{\SU(3)\times \SO(3)}_{\U^{\bullet}(2)}\right\}.
			\end{align*}
			Denote by $\lambda(z_{1},z_{2})=z_{1}\lambda_{1}+z_{2}\lambda_{2}$ the highest weight of $\varrho_{0}$ and $\mu(z_{3})=z_{3}\mu_{1}$ the highest weight of $\varrho_{1}$. The half of sum of positive roots of $\su(3)$ and $\su(2)$ are given by  $\rho_{0}=2\lambda_{1}+\lambda_{2}$ and $\rho_{1}=\nu_{1}$, respectively. By the Freudenthal formula, see Theorem \ref{thm eigenvalues casimir}  we get that
			\begin{align*}
				c_{g}(\varrho_{0}(z_{1},z_{2}))&=g_{1,1}(\lambda(z_{1},z_{2})+2\rho_{0},\lambda(z_{1},z_{2}))=g_{1,1}(z_{1}\lambda_{1}+z_{2}\lambda_{2}+4\lambda_{1}+2\lambda_{2},z_{1}\lambda_{1}+z_{2}\lambda_{2})\\
				=&g_{1,1}((z_{1}+4)\lambda_{1}+(z_{2}+2)\lambda_{2},z_{1}\lambda_{1}+z_{2}\lambda_{2})\\
				&=\frac{4(z_{1}(z_{1}+4)+(z_{2}(z_{2}+2)))}{3}-\frac{2(z_{2}(z_{1}+4)+z_{1}(z_{2}+2))}{3}\\
				&=\frac{4(z_{1}^{2}+z_{2}^{2}-z_{1}(z_{2}-3))}{3}
			\end{align*}
			as well as 
			\begin{align*}
				c_{g}(\varrho_{1}(z_{3}))=	g_{1,1}(\mu(z_{3})+2\rho_{1},\mu(z_{3}))=4z_{3}(z_{3}+1).
				%=g_{1,1}((z_{3}+1)\mu_{1},z_{3}\mu_{1})
			\end{align*}
			We obtain the desired formula for the spectrum.\nopagebreak		
			\hfill	\qed
			\subsection{Computation of spherical representations}\label{subsec: spherical representations}
			%-------------------------------------------------------------------------------
			%
			The computation of spherical representations is the heart of the spectrum computation. It is done using branching rules which are described in \cite[Sec 8.2]{Goodman Wallach}. We denote the maximal torus of $\u^{\bullet}(2)$ by $\t_{\u^{\bullet}(2)}$. The restriction of a torus element $\lambda\in i\t^{*} $ to $\t_{\u^{\bullet}(2)}$ is most essential when it comes to branching. For this purpose, we introduce the restriction notation  %of  and  its dual $i\t^{*}_{\u^{\bullet}(2)}$
			%  We denote the restriction 
			\begin{align*}
				\r{\lambda}:=\lambda\myrestriction \t_{\u^{\bullet}(2)},\quad 	 R^{+}_\bullet:=\{\r{\lambda}\mid \lambda\in R^{+}\}.
			\end{align*}
			Be aware that the set of restricted positive roots $R^{+}_\bullet $ does not coincide with the set of positive roots $ R^{+}_{\u^{\bullet}(2)} $ of $\u^{\bullet}(2)$, but in order to apply the branching rules, the comparison between both is essential. First we compute the torus $ i\t_{\u^{\bullet}(2)}^{*} $ and the corresponding system of positive roots $R^{+}_{\u^{\bullet}(2)}$. 
			\begin{lem}\label{lem computation i t_k*}
				%There is an element $ Y\in \t_{\u^{\bullet}(2)} $ such that $-i\lambda(Y)>0 $ for all $ \lambda\in R^{+} $.	
				The dual Lie algebra $ i\t_{\u^{\bullet}(2)}^{*} $ of the maximal torus of $ K=\U^{\bullet}(2) $ is given by
				\begin{align*}
					i\t_{\u^{\bullet}(2)}^{*}=\left\{\r{\lambda}_{1}a_{1}+a_{2}\r{\mu}_{1}\:\vline\: a_{1},a_{2}\in \R\right\}
				\end{align*}
				and the corresponding set of positive roots is $ R^{+}_{\u^{\bullet}(2)}=\{\r{\mu}_{1}\}\subset R^{+}_\bullet $. 
				% $ i\t^{*}_{\u^{\bullet}(2)} $ is given by
			\end{lem}
			\begin{proof}
				%	The first condition is the regularity assumption \cite[p. 370]{Goodman Wallach}. 
				%	%The regularity assumption states that  This is equivalent to the existence of $ X_{\alpha}\in \t_{\u^{\bullet}(2)} $ such that $ \lambda(X_{\alpha})\neq 0 $ for all $ \lambda\in R^{+} $. 
				%	Recall that we write $ X_{\lambda}=\#(-i\lambda) $ and $ \lambda=i\flat(X_{\lambda}) $ for $ \lambda\in i\t^{*} $. By equation \eqref{eq: D(SU(3)xSO(3))} in the proof of Proposition \ref{prop: D(SU(3)xSO(3))}, we see that $ \lambda(Y)=h_{1,1}(i\flat(Y),\r{\lambda})>0 $ is true if $ i\flat(Y) $ is an element in
				%	\begin{align*}
					%		i\t_{\u^{\bullet}(2)}^{*}\cap	\{z_{1}\lambda_{1}+z_{2}\lambda_{2}+z_{3}\mu_{1}\mid z_{1}> z_{2}> 0 ,\: z_{3}>0\}.%\subset 	D(SU(3)\times SO(3)).
					%	\end{align*}\textcolor{red}{diese Menge ist nicht wohldef.}
				The maximal torus $ \t_{\u^{\bullet}(2)} $ of $ \u^{\bullet}(2) $ is given by:
				\begin{align*}
					\t_{\u^{\bullet}(2)}=\left\{a_{1}E_{1}+a_{0}{E}_{0}+a_{1}e_{1}\:\vline\: a_{1},a_{0}\in \R\right\}.
				\end{align*}
				For any $\lambda\in i\t^{*}$ we denote by $X_{\lambda}$ its dual with respect to $ih_{1,1}(\cdot,\cdot)$, i.e.
				\begin{align*}
					X_{\lambda_{1}}=E_{1}+\frac{1}{3}{E}_{0},\quad  X_{\lambda_{2}}=\frac{1}{3}{E}_{0}-E_{1},\quad 
					X_{\mu_{1}}=2e_{1}.
				\end{align*}
				Solving for $E_{1}$, $E_{0}$ and $e_{1}$ and restricting to $\t_{\u^{\bullet}(2)} $ yields
				\begin{align*}
					E_{1}=\frac{1}{2}(X_{\lambda_{1}}-X_{\lambda_{2}}),\quad {E}_{0}=\frac{3}{2}(X_{\lambda_{1}}+X_{\lambda_{2}}),\quad e_{1}=\frac{X_{\mu_{1}}}{2}.
				\end{align*}
				Hence we get:
				\begin{align*}
					\t_{\u^{\bullet}(2)}&=\left\{X_{\lambda_{1}}(a_{1}+3a_{0})+X_{\lambda_{2}}(3a_{0}-a_{1})+a_{1}X_{\mu_{1}}\:\vline\: a_{1},a_{0}\in \R\right\}\\
					i\t_{\u^{\bullet}(2)}^{*}&=\left\{\r{\lambda}_{1}(a_{1}+3a_{0})+\r{\lambda}_{2}(3a_{0}-a_{1})+a_{1}\r{\mu}_{1}\:\vline\: a_{1},a_{0}\in \R\right\}.
				\end{align*}
				%	We see that the regularity condition is fulfilled for each element in $ i\t_{k}^{*} $ with $ 3a_{0}>a_{1}>0 $. 	
				For $ i=1,2,3 $, $ j=1,2 $ and $ A=a_{1}E_{1}+a_{0}{E}_{0}+a_{1}e_{1}\in \t_{\u^{\bullet}(2)} $ the restricted dual vectors $ \r{\lambda}_{i},,\r{\mu}_{j}\in i\t^{*}_{\u^{\bullet}(2)}  $ are given by \begin{align*}
					\r{\lambda}_{1}(A)&=i(a_{1}+ a_{0}),\quad \r{\lambda}_{2}(A)=i(a_{0}-a_{1}),\quad \r{\lambda}_{3}(A)=-2ia_{0}\\
					\r{\mu}_{j}(A)&=(-1)^{j+1}2ia_{1}.
					%(2\r{\lambda}_{1}+\r{\lambda}_{2})(A)&=i(a_{1}+3a_{0}),\quad (\r{\lambda}_{1}+2\r{\lambda}_{2})(A)=i(2a_{0}-a_{1}) ,\quad \r{\mu}_{1}(A)=2a_{1}
				\end{align*}
				We see that $ \r{\lambda}_{1}-\r{\mu}_{1}=\r{\lambda}_{2} $ as well as $  \r{\mu}_{1}-2\r{\lambda}_{1}=\r{\lambda}_{3}  $, so we can rewrite $ i\t_{\u^{\bullet}(2)}^{*} $:
				\begin{align*}
					i\t_{\u^{\bullet}(2)}^{*}&=\left\{\r{\lambda}_{1}(a_{1}+3a_{0})+\r{\lambda}_{2}(3a_{0}-a_{1})+a_{1}\r{\mu}_{1}\:\vline\: a_{1},a_{0}\in \R\right\}\\
					&=\left\{\r{\lambda}_{1}(a_{1}+3a_{0})+(\r{\lambda}_{1}-\r{\mu}_{1})(3a_{0}-a_{1})+a_{1}\r{\mu}_{1}\:\vline\: a_{1},a_{0}\in \R\right\}\\
					%i\t_{\u^{\bullet}(2)}^{*}&=\left\{6\r{\lambda}_{1}a_{0}+\r{\mu}_{1}(a_{1}-(3a_{0}-a_{1}))\:\vline\: a_{1},a_{0}\in \R\right\}\\
					%	i\t_{\u^{\bullet}(2)}^{*}&=\left\{6\r{\lambda}_{1}a_{0}+\r{\mu}_{1}(a_{1}+a_{1}-3a_{0})\:\vline\: a_{1},a_{0}\in \R\right\}\\
					%i\t_{\u^{\bullet}(2)}^{*}&=\left\{6\r{\lambda}_{1}a_{0}+\r{\mu}_{1}(a_{1}2-3a_{0})\:\vline\: a_{1},a_{0}\in \R\right\}\\
					&=\left\{a_{1}\r{\lambda}_{1}+a_{0}\r{\mu}_{1}\:\vline\: a_{1},a_{0}\in \R\right\}.
				\end{align*}
				In order to compute the roots of $\u^{\bullet}(2)$, the $\u(2)$ commutator relations are relevant: $[e_{i},e_{j}]=-2\varepsilon_{ijk}e_{k}$, $[E_{i},E_{j}]=-2\varepsilon_{ijk}E_{k}$, $[e_{i},e_{0}]=[E_{i},E_{0}]=0$ for all $i,j,k=1,2,3$. 
				For understanding the adjoint action of $\t_{\u^{\bullet}(2)}$ on $\u^{\bullet}(2)_{\C}=\u^{\bullet}(2)\oplus i\u^{\bullet}(2)$, define 
				\begin{align*}
					x+iy=\sum_{j=0}^{3}(x_{j}+iy_{j})E_{j}+\sum_{j=1}^{3}(x_{j}+iy_{j})e_{j}\in  \u^{\bullet}(2)_{\C},\quad a=a_{1}E_{1}+a_{0}{E}_{0},a_{1}e_{1}\in \t_{\u^{\bullet}(2)}.  
				\end{align*} The condition $[a,x+iy]=\lambda(a)(x+iy)$ for a $\lambda\in i\t^{*}_{\u^{\bullet}(2)}$ necessarily leads to
				\begin{align*}
					[a,x+iy]&=2a_{1}((x_{3}+iy_{3})(E_{2}+e_{2})-(x_{2}+iy_{2})(E_{3}+e_{3}))\\
					&=2ia_{1}((y_{3}-ix_{3})(E_{2}+e_{2})+(-y_{2}+ix_{2})(E_{3}+e_{3}))=\lambda(a)(x+iy).
				\end{align*} 
				This is equivalent to $x_{1}=y_{1}= x_{0}=y_{0}=0 $, $ x_{2}=y_{3},\: y_{2}=-x_{3},\ x_{3}=-y_{2},\: y_{3}=x_{2} $ and $ \lambda(a)=\r{\mu}_{1}(a) $.
			\end{proof}
			The partition function is most important to investigate the branching of a representation.
			\begin{dfn}
				For $ \beta\in R^{+}_\bullet $ we set $ B_{\beta}=\{\alpha\in R^{+}\mid \r{\alpha}=\beta\} $ and define 
				\begin{align*}
					\Gamma_{0}=\{\beta\mid \beta\in R^{+}_{\u^{\bullet}(2)}\:\text{and}\: |B_{\beta}|>1\},\quad \Gamma_{1}=R^{+}_\bullet\setminus R^{+}_{\u^{\bullet}(2)},\quad \Gamma=\Gamma_{0}\cup \Gamma_{1}.
				\end{align*}
				The \emph{multiplicity} $ m_{\beta} $ for $ \beta\in \Gamma $ is defined by
				\begin{align*}
					m_{\beta}=\begin{cases}
						|B_{\beta}|,\quad\beta\notin R^{+}_{\u^{\bullet}(2)}\\
						|B_{\beta}|-1,\quad\beta\in R^{+}_{\u^{\bullet}(2)}
					\end{cases}.
				\end{align*} With respect to $\Gamma$, the \emph{partition function} $ \wp :i\t_{\u^{\bullet}(2)}^{*} \rightarrow \N_{0}$ evaluates $ \xi $ to the number of ways of writing
				\begin{align*}
					\xi=\sum_{\beta\in \Gamma}c_{\beta}\beta,\quad c_{\beta}\in \N_{0}
				\end{align*}
				where each $ \beta $ that occurs is counted with multiplicity $ m_{\beta}. $
			\end{dfn}
			Hence, we need a precise description of the indexset $\Gamma.$
			\begin{lem}
				The set $ \Gamma $ is given by $\Gamma=R^{+}_\bullet=\{3\r{\lambda}_{1} -\r{\mu}_{1},3\r{\lambda}_{1}-2\r{\mu}_{1},\r{\mu}_{1}\}.$	Moreover we have $ \r\rho=\r{\lambda}_{1}-\r{\lambda}_{3}+\r{\nu}_{1} $, $ m_{\beta}=1 $ for all $\beta\in \Gamma.$
			\end{lem}
			\begin{proof}
				The positive roots are given by	$	R^{+}=\{\lambda_{1}-\lambda_{2},\lambda_{1}-\lambda_{3},\lambda_{2}-\lambda_{3},\mu_{1}\}.$
				For $ i=1,2,3 $, $ j=1,2 $ and $ A=a_{1}(E_{1}+e_{1})+a_{0}{E}_{0}\in \t_{\u^{\bullet}(2)} $ the restricted vectors  are given by \begin{align*}
					\r{\lambda}_{1}(A)=i(a_{1}+ a_{0}),\quad \r{\lambda}_{2}(A)=i(a_{0}-a_{1}),\quad \r{\lambda}_{3}(A)=-2ia_{0},\quad \r{\mu}_{j}(A)=(-1)^{j+1}2ia_{1}.
				\end{align*}
				Note that $ \r{\lambda}_{1}-\r{\lambda}_{2}=\r{\mu}_{1} $, $  \r{\lambda}_{1}-\r{\lambda}_{3}=3\r{\lambda}_{1} -\r{\mu}_{1}$ and $ \r{\lambda}_{2}-\r{\lambda}_{3}= 3\r{\lambda}_{1}-2\r{\mu}_{1}$. Hence we have $\Gamma_{0}=\{\r{\mu}_{1}\}$ and $ \Gamma_1=\{3\r{\lambda}_{1} -\r{\mu}_{1},3\r{\lambda}_{1}-2\r{\mu}_{1}\}$.
			\end{proof}
			From this point, the partition function can easily be computed:
			\begin{lem}\label{lem: partition function calculated}
				%	In terms of the restricted dual elements $ \r{\lambda}_{1},\r{\nu}_{1} $, the dual Lie algebra $ i\t_{\u^{\bullet}(2)}^{*} $ of the maximal torus of $ K=\U^{\bullet}(2) $ is given by
				%	\begin{align*}
					%		i\t_{\u^{\bullet}(2)}^{*}=\left\{\r{\lambda}_{1}a_{1}+a_{2}\r{\mu}_{1}\:\vline\: a_{1},a_{2}\in \R\right\}.
					%	\end{align*}	
				For $ a_{1}\r{\lambda}_{1}+a_{2}\r{\mu}_{1} \in i\t_{\u^{\bullet}(2)}^{*} $ the partition function $ \wp :i\t_{\u^{\bullet}(2)}^{*} \rightarrow \N_{0}$ is given by 
				\begin{align*}
					\wp(a_{1},a_{2})&=\# \{(m_{1},m_{2},m_{3})\in \N_{0}^{3}\mid a_{1}=3(m_{1}+m_{2}),\quad a_{2}=m_{3}-m_{1}-2m_{2}\}\\
					&=\begin{cases}
						0,\quad \text{if } \frac{a_{1}}{3}\notin \N_{0}\text{ or }a_{2}\notin \Z\\
						\max\{0, \frac{a_{1}}{3}+1+\min(0,a_{2}+\frac{a_{1}}{3})\},\quad \text{else}.
					\end{cases}
				\end{align*}
				In particular, it vanishes if $3a_{2}<-2a_{1}$.
				%
				%	\begin{align*}		\wp(a_{1},a_{2})
					%		\begin{cases}
						%			%0 ,\textit{ if $ n_{1}\notin 3\N_{0} $}\\
						%			=\frac{a_{1}}{3},\textit{ if $ a_{1}\in 3\N $ and $ a_{2}\in \N_{0}$}\\
						%			\in \:	 ]0,\frac{a_{1}}{3}],\textit{ if $ a_{1}\in 3\N $ and $ a_{2}\in \Z\cap [-2a_{1}/3,\infty)$}\\
						%			=1,\textit{ if $ a_{1}=0 $ and $ a_{2}\in \N_{0} $}\\
						%			=0,\textit{ else.}
						%			%2,\textit{ if $ n_{1}\in 3\N_{0} $ and $ n_{2}\geq 2n_{1} $.}
						%		\end{cases}
					%	\end{align*}  
			\end{lem}
			\begin{proof}
				%(Note that the linear combination of a set of vectors is unique for every element in its span, precisely when the set of vectors is linearly independent.)	
				Recall that $ \r{\lambda}_{1}-\r{\lambda}_{2}=\r{\mu}_{1} $ and assume that $ a_{1}\r{\lambda}_{1}+a_{2}\r{\mu}_{1}\in i\t_{\u^{\bullet}(2)}^{*} $ for $ a_{1},a_{2}\in \R $. The number of ways to write this element as a natural number linear combination over $\Gamma$ is precisely given by the number of $ m_{1},m_{2},m_{3}\in \N_{0} $ satisfying 
				\begin{align*}
					a_{1}\r{\lambda}_{1}+a_{2}\r{\mu}_{1}&=m_{1}(3\r{\lambda}_{1} -\r{\mu}_{1})+m_{2}(3\r{\lambda}_{1} -2\r{\mu}_{1})+m_{3}\r{\mu}_{1}\\
					&=\r{\lambda}_{1}(3m_{1}+3m_{2})+\r{\mu}_{1}(m_{3}-m_{1}-2m_{2}).
				\end{align*}
				In particular, we obtain that $\wp(a_{1},a_{2})=0 $ if $3a_{1}\notin\N_{0}$ or $a_{2}\notin\Z.$ We continue to compute
				\begin{align*}
					\wp(a_{1},a_{2})&=\#\{(m_{1},m_{2},m_{3})\in \N_{0}^{3}\mid \frac{a_{1}}{3}=m_{1}+m_{2},a_{2}=m_{3}-m_{1}-2m_{2}\}\\
					&=\#\{m_{3}\in \N_{0},m_{2}\in \{0,\ldots,\frac{a_{1}}{3}\}\mid a_{2}=m_{3}-m_{2}-\frac{a_{1}}{3}\}\\
					&=\#\{m_{2}\in \{0,\ldots, \frac{a_{1}}{3}\}\mid m_{2}+\frac{a_{1}}{3}+a_{2}\geq 0\}\\
					&=\begin{cases}
						0,\quad \text{if } \frac{a_{1}}{3}\notin \N_{0}\text{ or }a_{2}\notin \Z\\
						\max(0, \frac{a_{1}}{3}+1+\min(0,a_{2}+\frac{a_{1}}{3})),\quad \text{else}.
					\end{cases}
				\end{align*}
			\end{proof}
	
			With the computed partition function, we can start the branching and find the $\U^{\bullet}(2)$-spherical representations. 
			\begin{proof}[Proof of \eqref{item: m formula} in Theorem $\ref{thm: spectrum Aloff Wallach}$ ]
				%------------------------
				Let $\varrho(z_{1},z_{2},z_{3})\in \hat{\SU(3)\times \SO(3)}$ be an irreducible representation with corresponding highest weight $\lambda$. The multiplicity of the trivial representation when restricted to $\U^{\bullet}(2)$ can be computed by the formula
				\begin{align*}
					m(z_{1,}z_{2},z_{3})=\sum_{w\in W}\sgn(w)\wp(\r{w(\lambda+\rho)}-\r{\rho}),\quad \text{\cite[\text{Thm. 8.2.1}]{Goodman Wallach}}.
				\end{align*}
				Recall that $ \rho=2\lambda_{1}+\lambda_{2}+\nu_{1} $, $ \r{\lambda}_{1}-\r{\lambda}_{2}=\r{\mu}_{1} $ and introduce the variables $$p_{1}:=z_{1}+2,\quad p_{2}:=z_{2}+1,\quad p_{3}:=-z_{3},\quad p_{4}:=z_{3}+1.$$ Introducing the Weyl group
				\begin{align*}
					W=\{\sigma^{\pm}\mid \sigma\in S_{3}\},\quad	S_{3}=\{(1,2),(2,3),(1,3),\id,(1,2,3),(1,3,2)\}
				\end{align*} 
				we start to compute
				\begin{align*}
					m(z_{1},z_{2},z_{3})&=\sum_{w\in W}\sgn(w)\wp({w\r{({{\lambda_{1}}}p_{1}+{\lambda_{2}}p_{2}+(2z_{3}+1){\nu_{1}})}}-\r{\rho})\\
					%	=&\sum_{w\in W}\sgn(w)\wp({p_{1}w\r{({\lambda_{1}})}+p_{2}w\r{({\lambda_{2}})}+(2z_{3}+1)w\r{({\nu_{1}})}-2\r{\lambda}_{1}-\r{\lambda}_{2}-\r{\nu}_{1}})\\
					%	=&\sum_{\sigma\in S_{3}}\sgn(\sigma)\wp({p_{1}\r{\lambda}_{\sigma(1)}+p_{2}\r{\lambda}_{\sigma(2)}+2z_{3}\r{\nu}_{1}-3\r{\lambda}_{1}+2\r{\nu}_{1}})\\
					%	&-\sum_{\sigma\in S_{3}}\sgn(\sigma)\wp({p_{1}\r{\lambda}_{\sigma(1)}+p_{2}\r{\lambda}_{\sigma(2)}-(2z_{3}+2)\r{\nu}_{1}-3\r{\lambda}_{1}+2\r{\nu}_{1}})\\
					=&\sum_{j=3}^{4}(-1)^{j}\sum_{\sigma\in S_{3}}\sgn(\sigma)\wp({p_{1}\r{\lambda}_{\sigma(1)}+p_{2}\r{\lambda}_{\sigma(2)}-3\r{\lambda}_{1}+p_{j}\r{\mu}_{1}})\\
					%		=&(-1)^{j}\sum_{j=3}^{4}-\wp({p_{1}\r{\lambda}_{2}+p_{2}\r{\lambda}_{1}-3\r{\lambda}_{1}+p_{j}\r{\mu}_{1}}) -\wp({p_{1}\r{\lambda}_{1}+p_{2}\r{\lambda}_{3}-3\r{\lambda}_{1}+p_{j}\r{\mu}_{1}})\\
					%		&-\wp({p_{1}\r{\lambda}_{3}+p_{2}\r{\lambda}_{2}-3\r{\lambda}_{1}+p_{j}\r{\mu}_{1}})+\wp({p_{1}\r{\lambda}_{1}+p_{2}\r{\lambda}_{2}-3\r{\lambda}_{1}+p_{j}\r{\mu}_{1}})\\
					%		&+\wp({p_{1}\r{\lambda}_{2}+p_{2}\r{\lambda}_{3}-3\r{\lambda}_{1}+p_{j}\r{\mu}_{1}})+\wp({p_{1}\r{\lambda}_{3}+p_{2}\r{\lambda}_{1}-3\r{\lambda}_{1}+p_{j}\r{\mu}_{1}})\\
					=&\sum_{j=3}^{4}(-1)^{j}\left(-\wp({(p_{2}-3)\r{\lambda}_{1}+p_{1}\r{\lambda}_{2}+p_{j}\r{\mu}_{1}}) -\wp({(p_{1}-3)\r{\lambda}_{1}+p_{2}\r{\lambda}_{3}+p_{j}\r{\mu}_{1}})\right.\\
					&\left.-\wp({p_{1}\r{\lambda}_{3}+p_{2}\r{\lambda}_{2}-3\r{\lambda}_{1}+p_{j}\r{\mu}_{1}})+\wp({(p_{1}-3)\r{\lambda}_{1}+p_{2}\r{\lambda}_{2}+p_{j}\r{\mu}_{1}})\right.\\
					&\left.+\wp({p_{1}\r{\lambda}_{2}+p_{2}\r{\lambda}_{3}-3\r{\lambda}_{1}+p_{j}\r{\mu}_{1}})+\wp({p_{1}\r{\lambda}_{3}+(p_{2}-3)\r{\lambda}_{1}+p_{j}\r{\mu}_{1}})\right).
				\end{align*}
				Now we use the equations $ {\r{\lambda}_{1}}-{\r{\mu}_{1}}=\r{\lambda}_{2} $, $  {\r{\mu}_{1}}-2{\r{\lambda}_{1}}=\r{\lambda}_{3}  $:
				\begin{align*}
					&\sum_{j=3}^{4}(-1)^{j}\left(-\wp({(p_{1}+p_{2}-3)\r{\lambda}_{1}+(p_{j}-p_{1})\r{\mu}_{1}}) -\wp({(p_{1}-3-2p_{2})\r{\lambda}_{1}+(p_{2}+p_{j})\r{\mu}_{1}})\right.\\
					&\left.-\wp({(p_{2}-3-2p_{1})\r{\lambda}_{1}+(p_{1}+p_{j}-p_{2})\r{\mu}_{1}})+\wp({(p_{1}+p_{2}-3)\r{\lambda}_{1}+(p_{j}-p_{2})\r{\mu}_{1}})\right.\\
					&\left.+\wp({(p_{1}-3-2p_{2})\r{\lambda}_{1}+(p_{2}+p_{j}-p_{1})\r{\mu}_{1}})+\wp({(p_{2}-3-2p_{1})\r{\lambda}_{1}+(p_{1}+p_{j})\r{\mu}_{1}})\right).
				\end{align*}
				Expanding the sum and plugging in $ p_{1}=z_{1}+2, p_{2}=z_{2}+1,p_{3}=-z_{3},p_{4}=z_{3}+1 $ yields:
				%	\begin{align*}
					%m(z_{1},z_{2},z_{3})	=&\sum_{j=3}^{4}(-1)^{j}\left(-\wp({(z_{1}+z_{2})\r{\lambda}_{1}+(p_{j}-(z_{1}+2))\r{\mu}_{1}})\right.\\
					%		& -\wp({(z_{1}-2z_{2}-3)\r{\lambda}_{1}+(z_{2}+1+p_{j})\r{\mu}_{1}})\\
					%		&-\wp({(z_{2}-2z_{1}-6)\r{\lambda}_{1}+(2(z_{1}+2)+p_{j}-(z_{2}+1))\r{\mu}_{1}})\\
					%		&+\wp({(z_{1}+z_{2})\r{\lambda}_{1}+(p_{j}-(z_{2}+1))\r{\mu}_{1}})\\
					%		&+\wp({(z_{1}-2z_{2}-3)\r{\lambda}_{1}+(z_{2}-1-z_{1}+p_{j})\r{\mu}_{1}})\\
					%		&\left.+\wp({(z_{2}-2z_{1}-6)\r{\lambda}_{1}+(z_{1}+2+p_{j})\r{\mu}_{1}})\right)\\	
					%	\end{align*}
				%According to Lemma \ref{lem: partition function calculated} the $ 3^{rd} $ and the $ 6^{th} $ summand vanishes. :
				\begin{align*}
					m(z_{1},z_{2},z_{3})=&\wp({(z_{1}+z_{2})\r{\lambda}_{1}+(-z_{3}-z_{1}-2)\r{\mu}_{1}})\\
					&+ \wp({(z_{1}-2z_{2}-3)\r{\lambda}_{1}+(z_{2}+1-z_{3})\r{\mu}_{1}})\\
					&+\wp({(z_{2}-2z_{1}-6)\r{\lambda}_{1}+(-z_{3}-z_{2}+z_{1}+1)\r{\mu}_{1}})\\
					&-\wp({(z_{1}+z_{2})\r{\lambda}_{1}+(-z_{3}-z_{2}-1)\r{\mu}_{1}})\\
					&-\wp({(z_{1}-2z_{2}-3)\r{\lambda}_{1}+(z_{2}-1-z_{1}-z_{3})\r{\mu}_{1}})\\
					&-\wp({(z_{2}-2z_{1}-6)\r{\lambda}_{1}+(z_{1}+2-z_{3})\r{\mu}_{1}})\\
					&-\wp({(z_{1}+z_{2})\r{\lambda}_{1}+(z_{3}-z_{1}-1)\r{\mu}_{1}})\\
					& -\wp({(z_{1}-2z_{2}-3)\r{\lambda}_{1}+(z_{2}+2+z_{3})\r{\mu}_{1}})\\
					&-\wp({(z_{2}-2z_{1}-6)\r{\lambda}_{1}+(z_{1}-z_{2}+2+z_{3})\r{\mu}_{1}})\\
					&+\wp({(z_{1}+z_{2})\r{\lambda}_{1}+(z_{3}-z_{2})\r{\mu}_{1}})\\
					&+\wp({(z_{1}-2z_{2}-3)\r{\lambda}_{1}+(z_{2}-z_{1}+z_{3})\r{\mu}_{1}})\\
					&+\wp({(z_{2}-2z_{1}-6)\r{\lambda}_{1}+(z_{1}+z_{3}+3)\r{\mu}_{1}}).
				\end{align*}
				Lemma \ref{lem: partition function calculated} implies that five of the mentioned summands vanish: For the summands $ 3 $, $ 6 $, $ 9 $, $ 12 $, the first entry in $ \wp $ is $ z_{2}-2z_{1}-6 $ which is always negative. For summand $ 5 $ note that the first entry is $ a_{1}:=z_{1}-2z_{2}-z_{3} $ and the second entry is $ a_{2}:=z_{2}-1-z_{1}-z_{3} $. One quickly verifies $ a_{2}<-2a_{1}/3 $ which proves that summand $ 5 $ vanishes. This results in
				\begin{align*}
					m(z_{1},z_{2},z_{3})	=&\wp({(z_{1}+z_{2})\r{\lambda}_{1}+(-z_{3}-z_{1}-2)\r{\mu}_{1}})\\
					&+ \wp({(z_{1}-2z_{2}-3)\r{\lambda}_{1}+(z_{2}+1-z_{3})\r{\mu}_{1}})\\
					&-\wp({(z_{1}+z_{2})\r{\lambda}_{1}+(-z_{3}-z_{2}-1)\r{\mu}_{1}})\\
					&-\wp({(z_{1}+z_{2})\r{\lambda}_{1}+(z_{3}-z_{1}-1)\r{\mu}_{1}})\\
					& -\wp({(z_{1}-2z_{2}-3)\r{\lambda}_{1}+(z_{2}+2+z_{3})\r{\mu}_{1}})\\
					&+\wp({(z_{1}+z_{2})\r{\lambda}_{1}+(z_{3}-z_{2})\r{\mu}_{1}})\\
					&+\wp({(z_{1}-2z_{2}-3)\r{\lambda}_{1}+(z_{2}-z_{1}+z_{3})\r{\mu}_{1}}).
				\end{align*}
				%A short computation shows that if the first summand does not vanish then the second will not vanish either.
				%\textcolor{red}{Überflüssig?:}
				%	Note that $ z_{1}-2z_{2}-3 =z_{1}+z_{2}-3z_{2}-3$ is divisible by $ 3 $ if and only if $ z_{1}+z_{2} $ is divisible by $ 3 $. Hence we get that $\varrho_{\lambda} $ is not $\U^{\bullet}(2)$-spherical if $ z_{1}+z_{2} $ is not divisible by $ 3 $. 
				As the partition function is complicated to compute, the calculation to find which representations $\varrho(z_{1},z_{2},z_{3})$ are $\U^{\bullet}(2)$-spherical, i.e. satisfy $m(z_{1},z_{2},z_{3})\neq 0$, is performed using a suitable Python script \cite{Agricola Henkel 24}. 
			\end{proof}
			Next we investigate the condition to obtain $S^{1}$-spherical representations more closely.
			\begin{proof}[Proof of \eqref{item: S1 spherical} in Theorem $\ref{thm: spectrum Aloff Wallach}$ ]
				%------------------------------------
				In this part of the proof we change the simple roots to $\lambda_{1}-\lambda_{2}$ and $-(\lambda_{1}-\lambda_{3})=-(2\lambda_{1}+\lambda_{2})$. By doing so, we introduce $n_{1},n_{2}\in \N_{0}$ and switch the parametrization of the highest weights to $n_{1}\lambda_{1}-n_{2}\lambda_{2}$. It translates to the old parametrization for the matching of coefficients $z_{1}=n_{1}+n_{2}$ and $z_{2}=n_{2}$. We will prove that $\varrho(n_{1},n_{2})$ is $S^{1}$-spherical if and only if $n_{1}+2n_{2}\equiv 0 \pmod 3$ which is equivalent to $z_{1}+z_{2}\equiv 0\pmod 3$.
				The weights of an irreducible representation are precisely contained in the convex hull of the Weyl group orbit of its highest weight. The weight vectors are given by adding integer linear combinations of roots to the highest weight \cite[Thm. 10.1]{Hall}. As the maximal torus $\t_{\su(3)}=\{E_{0},E_{1}\}$ of $\SU(3)$ contains the Lie algebra of $S^{1}$, $S^{1}$-spherical representation are precisely those which contain a weight vector $\lambda=n_{1}\lambda_{1}-n_{2}\lambda_{2}$ satisfying $\lambda(E_{0})=0 $. This condition for the weight is equivalent to $n_{1}=n_{2}$ which clearly satisfies $n_{1}+2n_{2}\equiv 0\pmod 3$. The relation is preserved under adding or subtracting the positive simple roots $\lambda_{1}-\lambda_{2}$, $-(\lambda_{1}+2\lambda_{2})$ and hence also applies to the highest weight vector.
				
			\end{proof}
			We were able to identify certain representations $\varrho(z_1,z_2,z_3)$ that are always  $\U^{\bullet}(2)$-spherical
			representations. They are of the form $\varrho(2z_2,z_2,z_3)$ for $0\leq z_3\leq z_2$. The case $z_2=z_3$ is of particular interest. A further discussion of this result will be given in Section \ref{NaSe}.
			\begin{cor}\label{cor: first eigenvalue horizontal}\begin{enumerate}
					\item[] 
					\item
					Fix $z_{2}\in \N_{0}$. Then $	\mult(2 z_{2},z_{2},z_{3}) = 0$ for any $z_3>z_2$ while
					the representations $\varrho(2z_{2},z_{2},z_{3})$ are $\U^{\bullet}(2)$-spherical with  multiplicity 
					\begin{align*}
						\mult(2 z_{2},z_{2},z_{3})=\frac{(z_{2}+1)(2z_{2}+2)(z_{2}+1)(2z_{3}+1)}{2}
						\text{ for } \ z_3=0,\ldots, z_2.
					\end{align*}
					The corresponding eigenvalues are given by
					\begin{align*}
						\eta(2z_{2},z_{2},z_{3})=\frac{4z_{2}}{t_{0}}+\frac{2z_{3}(2z_{3}+2)}{t_{1}}.
					\end{align*} 
					\item
					Any $\U^{\bullet}(2)$-spherical representation 
					$\varrho(z_{1},z_{2},z_{3})$ satisfies the eigenvalue inequality  
					\[
					\eta(z_{1},z_{2},z_{3})\geq \eta(2z_{3},z_{3},z_{3}).
					\] 
					In particular:  For a fixed representation $\varrho(z_{3})\in \hat{\SO(3)}$ of the fibre, the representation $\varrho(2z_{3},z_{3})$ yields the first eigenvalue $\eta_{1}^{\mathcal{H}}=4z_{3}/t_{0}$ of the horizontal Laplacian $\Delta_{\mathcal{H}}$,  restricted to the eigenspace of the corresponding vertical eigenvalue $\frac{2z_{3}(2z_{3}+2)}{t_{1}}$. 
				\end{enumerate}
			\end{cor}
			\begin{proof}
				For representations being of the described type we obtain that the computation of $m(2z_{2},z_{2},z_{3})$ becomes significantly more easy: The three summands having $z_{1}-2z_{2}-3$ as the first entries of the partition functions vanish as those entries are negative. The remaining four summands are of the type
				\begin{align*}
					\wp(3z_{2},a_{2})	=\begin{cases}
						0,\quad \text{if } \frac{a_{1}}{3}\notin \N_{0}\text{ or }a_{2}\notin \Z\\
						\max(0,z_{2}+1+\min(0,a_{2}+z_{2})),\quad \text{else}.
					\end{cases}
				\end{align*} 
				A short computation yields the desired result:
				\begin{align*}
					&m(2z_{2},z_{2},z_{3})\\
					&=\wp(3z_{2},-z_{3}-2z_{2}-2)+\wp(3z_{2},z_{3}-z_{2})-\wp(3z_{2},-z_{3}-z_{2}-1)-\wp(3z_{2},z_{3}-2z_{2}-1)\\
					&=\max(0,z_{2}+1+\min(0,-z_{3}-z_{2}-2 ))+\max(0,z_{2}+1+\min(0,z_{3}))\\
					&-\max(0,z_{2}+1+\min(0,-z_{3}-1 ))-\max(0,z_{2}+1+\min(0,z_{3}-z_{2}-1)) \\
					&=0+z_{2}+1-\max(0,z_{2}-z_{3})-\max(0,z_{2}+1+\min(0,z_{3}-z_{2}-1))\\
					& =\begin{cases}
						1,\text{ if } z_{3}\leq z_{2}\\
						0,\text{ if } z_{3}>z_{2}.
					\end{cases}
				\end{align*}
				Next we compute the first eigenvalue of the horizontal Laplacian for fixed $\varrho(z_{3})\in \hat{\SO(3)}$. In order to do this, note that the scalar $c_{g}(\lambda)=g(\lambda,2\rho +\lambda)$ increases whenever it is applied to a weight which is higher than $\lambda$: We consider an arbitrary positive root $\lambda^{+}\in R^{+}$, a dominant integral weight $\lambda$ and conclude
				\begin{align*}
					c_{g}(\lambda+\lambda^{+})&=g(\lambda+\lambda^{+},\lambda+\lambda^{+}+2\rho)\\
					&=g(\lambda+2\rho+\lambda)+\underset{\geq 0}{\underbrace{g(\lambda,\lambda^{+})}}+\underset{\geq 0}{\underbrace{g(\lambda^{+},\lambda)}}+\underset{\geq 0}{\underbrace{g(\lambda^{+},2\rho)}}+\underset{\geq 0}{\underbrace{g(\lambda^{+},\lambda^{+})}}\geq c_{g}(\lambda).
				\end{align*}   In a similar way to the calculation above, one proves that $m(z_{1},z_{2},z_{3})=0$ if $2z_{1}-z_{2}<3z_{3}$ or if $z_{1}+z_{2}<3z_{3}$. Assume that $\varrho(z_{1},z_{2},z_{3})$ is $\U^{\bullet}(2)$-spherical, i.e. $3a:=2z_{1}-z_{2}-3z_{3}\geq 0 $ and $3b:=z_{1}+z_{2}-3z_{3}\geq 0$. Adding $a(\lambda_{1}-\lambda_{2})+b(\lambda_{2}-\lambda_{3})$ to $3z_{3}\lambda_{1}+z_{3}\lambda_{2}+z_{3}\mu$ yields $z_{1}\lambda_{1}+z_{2}\lambda_{2}+z_{3}\mu$. This proves that the highest weight of $\varrho(2z_{3},z_{3},z_{3})$ is lower than or equal to the highest weight of $\varrho(z_{1},z_{2},z_{3})$.
				
			\end{proof}
			With the computation of the spectrum of the Aloff-Wallach manifold, we are now able to complete Urakawa's list of the first eigenvalue of a compact simply connected normal homogeneous space with positive curvature.
			
			\begin{cor}\label{cor: list of first eigenvalue on compact simply connected normal hom pos curvature}
				The only compact simply connected normal homogeneous space which admits a positively curved family of metrics is the Aloff-Wallach manifold $(V_{3},h_{r_{0},r_{1}})$ with $r_{0},r_{1}>0$. The completed list of the first eigenvalue of a compact simply connected normal homogeneous space with positive curvature which is not diffeomorphic to a sphere can be found in Table \ref{Table: first eigenvalues}. Each space $G/K$ which is not the Aloff-Wallach manifold is equipped with the restricted Killing form of $G$. 
			\end{cor}
				\begin{table}[h!]
				\centering
				\renewcommand{\arraystretch}{1.3}
				\begin{tabular}{cc}
					\toprule
					$G/K$ & $\eta_{1}$\\
					\midrule
					$\SU(n+1)/S(\U(n)\times\U(1)) $ & $1$ \\
					%	\textcolor{red}{diese und SP(m) muss raus!} $\SO(m)$ & $\eta_{1}^{B} \geq 2\alpha\delta \cdot 2\frac{n^{2}+2n}{n-1}$ \\
					$\Sp(n+1)/(\Sp(n)\times\Sp(1))$& $\frac{n+1}{n+2}$ \\
					${F_4}/\Spin(9)$ & $\frac{2}{3}$ \\
					$\Sp(2)/\SU(2)$ & $\frac{8}{3}$ \\
					$\SU(5)/(\Sp(2)\times S^{1})$ & $1$\\
					$\Sp(n)/(\Sp(n-1)\times S^{1})$&$\frac{4n}{n+1}$ \\
					$\SU(3)\times\SO(3)/\U^{\bullet}(2)$& $\qquad\quad\:\:\:\qquad \frac{12}{r_{0}},\quad r_{0},r_{1}>0.$\\
					\bottomrule
				\end{tabular}
				\caption{The first eigenvalue of a compact simply connected normal homogeneous space with positive curvature}
				\label{Table: first eigenvalues}
				\vspace{-0.5cm}		
			\end{table}	
			\begin{proof}
				The classification of simply connected normal homogeneous spaces with positive sectional curvature summarized by Urakwa \cite{Urakawa list} misses the normal homogeneous realization of the Aloff-Wallach manifold $W^{1,1}$ which was constructed by Wilking \cite{Wilking} $13$ years later. As we computed in Theorem  \nolinebreak \ref{thm: spectrum Aloff Wallach} the first eigenvalue of $(W^{1,1},g_{t_{0},t_{1}})$ is given by $2\alpha\delta\cdot 12$ with $2\alpha\delta>\delta^{2}$ which translates to $\frac{12}{r_{0}}$ with respect to the normal homogeneous realization $(V_{3},h_{r_{0},r_{1}}),$ $r_{0},r_{1}>0$. The first eigenvalues of the other spaces can be found in Urakawa's list \cite{Urakawa list} except for the cases $\Sp(2)/\SU(2)$ and $\SU(5)/(S^{1}\times\Sp(2))$ for which Urakawa only gives eigenvalue estimates. Each space $G/K$ is equipped with the restricted Killing form of $G$. Urakawa gave the first eigenvalue $\frac{5}{12}$ of $\Sp(2)$ as an lower estimate for the first eigenvalue of $\Sp(2)/\SU(2)$. The embedding of $\SU(2)$ into $\SU(3)$ is not standard: %equipped with the restricted Killing form of $\Sp(2)$. 
				The Lie algebra $\su(2)$ is embedded as the maximal subalgebra of $\su(2)$ in $\sp(2)$. Using the fundamental weights $\omega_{1},\omega_{2}$ of $\sp(2)$ described in Proposition \ref{prop: fundamental weights of sp(2)}, the eigenvalue of the Casimir element of a representation with highest weight $n_{1}\omega_{1}+n_{2}\omega_{2}$ is given by 
				\begin{align*}
					\frac{n_{1}(n_{1}+2)+n_{1}(n_{2}+2)+n_{2}(n_{1}+2)+n_{2}(n_{2}+2)+2n_{2}(n_{2}+2)}{12},
				\end{align*} 
				see \cite[p. 106]{Yamaguchi}. The computation of the branching rules of a maximal subalgebra is well described in \cite{McKay} and can be executed using the software LiE by van Leeuwen. Using that the eigenvalue of the Casimir element is monotonically increasing whenever $n_{1}$ or $n_{2}$ increases, one finds that the spherical representation $\varrho(4\cdot\omega_{1}+0\cdot\omega_{2})$ produces the smallest eigenvalue \nolinebreak $\frac{8}{3}$. Let us proceed with the second case $\SU(5)/(S^{1}\times\Sp(2))$. If $\lambda_{i}$ denotes for $i=1,\dots,5$ the projection onto the $i $-th diagonal element of the diagonal embedded torus of $\su(5)$, a system of simple roots is given by $\lambda_{i}-\lambda_{j}$ for $1\leq \i<j\leq 4$ and the corresponding fundamental weights are given by $\omega_{i}=\sum_{j=1}^{i}\lambda_{j}$, see \cite[Chap. 15.1]{Fulton_Harris}. The Freudenthal formula evaluates for a representation of highest weight $\lambda=\sum_{j=1}^{4}n_{j}\omega_{j}$ to
				\begin{align*}
					%	&\frac{n_{1}p_{1}+n_{2}p_{2}+n_{3}p_{3}+n_{4}p_{4}+2(p_{1}+p_{2}+p_{3}+p_{4})}{10}\\
					%	&=\frac{4n_{1}^{2}+3n_{1}n_{2}+2n_{1}n_{3}+n_{1}n_{4}+3n_{2}n_{1}+6n_{2}^{2}+4n_{2}n_{3}+2n_{2}n_{4}+2n_{3}n_{1}+4n_{3}n_{2}+6n_{3}^{2}+3n_{3}n_{4}+n_{4}n_{1}+2n_{4}n_{2}+3n_{4}n_{3}+4n_{4}^{2}+2(10n_{1}+15n_{2}+15n_{3}+10n_{4})}{150}\\
					&\frac{(4n_{1}+3n_{2}+2n_{3}+n_{4})(2+n_{1})+(3n_{1}+6n_{2}+4n_{3}+2n_{4})(2+n_{2})}{50}\\
					&+\frac{(2n_{1}+4n_{2}+6n_{3}+3n_{4})(2+n_{3})+(n_{1}+2n_{2}+3n_{3}+4n_{4})(2+n_{4})}{50},
				\end{align*}
				see Yamaguchi \cite[p. 106]{Yamaguchi}. Next, we want to compute the branching $\Sp(2)S^{1}\subset  \SU(5)$ using the program LiE. We obtain the corresponding restriction matrix by multiplying the restriction matrix of $\Sp(2)S^{1}\subset S(\U(4)\times S^{1}) $ with the one the complex Grassmannian $S(\U(4)\times S^{1}) \subset  \SU(5)$. That can be done using that the embedding $\Sp(2)\subset \U(4)$ is maximal, $S^{1}= S^{1}$ is the identity and that the branching of the complex Grassmannian can be done using LiE. %Computing the first spherical representations 
				%The program also provides the restriction matrix of this branching. The restriction matrix of the desired branching  is then the product of the restriction matrix of $\Sp(2)S^{1}\rightarrow\ $ $\Sp(2)\rightarrow \U(4)$ is maximal and 
				One finds that the spherical representation $\varrho(\omega_{1}+\omega_{4})$ produces the first eigenvalue $1$ which is the upper bound stated by Urakawa.
			\end{proof}
			\subsection{Comments on the other homogeneous
				$3$-$(\alpha,\delta)$-Sasaki manifolds}
			%------------------------------------------------------------------------
			We shall now explain how the root system relates to spherical representations (for a complete algebraic classification of all admissible highest weights for the classical series, we refer the reader to \cite{AgricolaCaglieroHenkel26}): Let $G$ be compact and simple, $K\subset G$ such that $G/K$ is a $3$-Sasakian manifold with fibre being $H\in\{\SU(2),\SO(3)\}$. As $H$ acts isometrically, freely and transitively on each fibre, each $H$-representation occurs. Choose the maximal torus $\t$ of $\g$ which contains the maximal torus $\t_{\k}$ of $\k$ and the maximal torus $\t_{\su(2)}$ of $\su(2)$. Those yield the orthogonal decomposition $\t=\t_{\k}\oplus\t_{\su(2)}$. The root system of $\su(2)$ can after a right choice of simple roots be associated to the maximal positive root $\frac{\mu_{1}}{2}$ which is contained in the fundamental Weyl chamber \cite{Goertsches Roschig Stecker}. In the case of $W^{1,1}$, this root is given by $\lambda_{1}-\lambda_{2}$. Any $K'$-spherical representation $\varrho_{0}\otimes \varrho_{1}$ corresponds to a $K$-spherical representation $\varrho_{0}$, i.e. the subgroup $K$ acts trivially on at least one vector. Hence, there must exist at least one weight of $\varrho_{0}$ which vanishes on $\t_{\k} $ and the corresponding weight vector is fixed by $K$. As the weight can only act non-trivially on $\t_{\su(2)}$, it is a half-integer multiple of the highest root. At the same time, the highest root describes the action of $\t_{\su(2)}$ on the vector fixed by $K$. The weights which occur in this way are precisely the weights of the $H$-representation $\varrho_{0}^{H,K}$ on the space of fixed vectors $V_{\varrho_{0}}^{K}$, see Theorem \ref{thm: spectrum main result}. It decomposes into the desired $H$ representations $\Phi_{\varrho_{0}}=\{\varrho_{1}\in \hat{H}\mid \varrho_{0}\otimes\varrho_{1}\in \hat{G\times H}_{K'}\}$, whose corresponding highest weights can be found as half-integer multiples of the highest root. \\
			
			\begin{figure}[h!]
				\centering				
				\begin{tikzpicture}[scale=1]
					% Define the range for the axes
					\def\xmin{-4}
					\def\xmax{4}
					\def\ymin{-2}
					\def\ymax{2}
					
					% Draw coordinate axes with arrows at the ends
					\draw[->] (\xmin,0) -- (\xmax+0.5,0) node[anchor=north west] {$\t_{\su(2)}$}; 
					\draw[->] (0,\ymin) -- (0,\ymax+0.5) node[anchor=south east] {$\t_\k$};

					\foreach \x in {-3,-2,-1,0,1,2,3}
					{
						\fill (\x,0) circle (2pt); 
					}

					\node[below] at (3,0) {$\frac{z}{2}\mu_{1}$};
				\end{tikzpicture}
				
				\caption{$H$ weights of a $K'$-spherical representation $\varrho_{0}\otimes \varrho_{1}(z)$.}
				\label{fig:su2_weights}
			\end{figure}
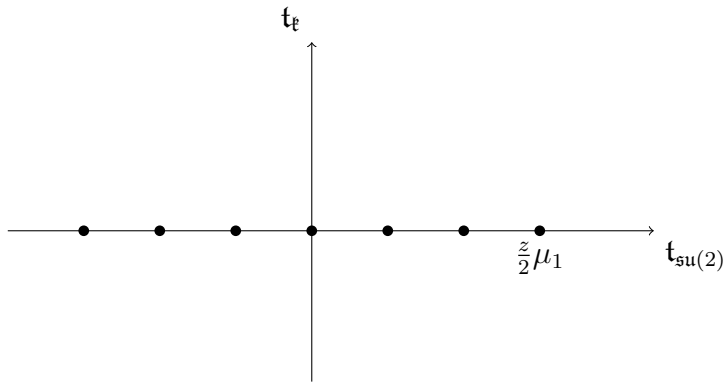
			
			We proved in Corollary \ref{cor: first eigenvalue horizontal} that in the case of the Aloff-Wallach space, the restriction $\varrho(z_{1},z_{2})$ of any $\U^{\bullet}(2)$-spherical representation $\varrho(z_{1},z_{2},z_{3})$ to $\SU(3)$ has highest weight which is higher than or equal to $z_{3}(\lambda_{1}+\lambda_{2})$, where $\lambda_{1}+\lambda_{2}$ is the maximal root of $\su(3)$ with respect to the simple roots $\lambda_{1}-\lambda_{2},\lambda_{2}-\lambda_{3}$ corresponding to the parameterization $z_{1},z_{2},z_{3}$. Hence, the representation $\varrho(2z_{3},z_{3},z_{3})$ yields the first eigenvalue of the horizontal Laplacian with respect to the fibre representation $\varrho(z_{3})$. This is completely in line with our general results described in Remark \ref{rem: horizontal eigenvalue estimates for each z}.\\
			
			Except for the spheres $S^{4n+3}=\Sp(n+1)/\Sp(n)$ which have fibre $\SU(2)$, all positive homogeneous $3$-$(\alpha,\delta)$-Sasaki manifolds have fibre $\SO(3)$, see \cite{Boyer Galicki Mann}. By identifying the corresponding $\su(2)$-representations with multiples of the maximal root, this can be seen in the root data: The Lie algebra $\sp(n)$ is the only real compact simple Lie algebra for which the half of the maximal root is a dominant integral element. In the other cases, $\SU(2)$-representations can not be realized, as positive half integer multiples of the maximal root are not integral elements. However, $\SO(3)$-representations can be realized, as positive integer multiples of the maximal root are for all simple Lie algebras integral elements.

			Let us investigate the relation between the maximal root and spherical representations in the example of $S^{7}=\Sp(2)/\Sp(1)$:
			
			\subsubsection*{The case of $\Sp(2)/\Sp(1)$}
			%--------------------------------------------
			%
			Consider $(\Sp(2)/\Sp(1),g_{t_{0},t_{1}})$ equipped with its $3$-$(\alpha,\delta)$-Sasaki metric. Here, the isotropy group $K$ is the lower right embedded $\Sp(1)$ and $H$ is the upper left embedded $\Sp(1)$. Its Lie algebra $\h$ is the Lie algebra of Reeb vector fields. The spectrum of $\Sp(2)/\Sp(1)$ can be found in \cite{Bettiol}, but in order to understand the relation of $\Sp(1)'$-spherical representations to their root system, we compute the $\Sp(1)'$-spherical representations explicitly. The naturally reductive realization of $\Sp(2)/\Sp(1)$ is given by
			\begin{align*}
				((\Sp(2)\times\Sp(1))/\Sp(1)',h_{r_{0},r_{1}}),\quad \Sp(1)'=\{(AB,B)\mid A\in K,\: B\in H\}.
			\end{align*} 
			Arbitrary maximal tori $\t_{\k}\subset \k$ and $\t_{\h}\subset\h$ give rise to the maximal torus of $\sp(2)\oplus\sp(1)$: $$\t_{\sp(2)\oplus\h}=\t_{\k}\oplus\t_{\h}\oplus\t_{\h}=\spann_{\R}\{E_{0},E_{1},e_{1}\}.$$  We denote by $\nu_{0},\nu_{1},\nu_{1}'$ the same $\sp(1)$ roots, but considered as roots of $\k$, $\h$ and the second copy of $\h$. 
			Following the construction in \cite[p. 102]{Yamaguchi}, combined with \cite[Cha. 8.2]{Hall}, one easily concludes:
			\begin{prop}\label{prop: fundamental weights of sp(2)}
				%---------------
				With respect to the simple roots $\left\{\nu_{0},\frac{1}{2}(\nu_{1}-\nu_{0}),\nu_{1}'\right\}$ of $\sp(2)\oplus\sp(1)$, the fundamental weights are given by $$\omega_{1}=\nu_{1}/2,\quad \omega_{2}=(\nu_{1}+\nu_{0})/2,\quad \omega_{3}=\nu_{1}'/2.$$ They parameterize the highest weights of irreducible representations $$\varrho(n_{1},n_{2},n_{3})\in \hat{\Sp(2)\times\Sp(1)},\quad n_{1},n_{2},n_{3}\in \N_{0}. $$ The sum of positive roots is given by $\rho={\nu_{0}}/{2}+\nu_{1}+\nu_{1}'/2$. The Weyl group of $\Sp(2)$ is the permutation of the indices ${j}\in \{0,1\}$ and an arbitrary number of sign changes.
			\end{prop} Note that the corresponding maximal root of $\mathfrak{sp}(2)$ with respect to those simple roots is given by $2\omega_{1}$. In the notation introduced in Subsection \ref{subsec: spherical representations}, one immediately computes the torus $\t_{\k'}$ and its dual $i\t^{*}_{\k'}$
			\begin{align*}
				\t_{\k'}=a_{0}E_{0}+a_{1}(E_{1}+e_{1}),\quad i\t^{*}_{\k'}=\{a_{0}\nu_{0}^{\bullet}+a_{1}\nu_{1}^{\bullet}\},
			\end{align*}
			the restricted positive roots of $\sp(2)\oplus\sp(1)$ and the positive roots of $\k'$
			\begin{align*}
				R^{+}_{\bullet}=\{\nu_{0}^{\bullet},\nu_{1}^{\bullet},\frac{1}{2}(\nu_{0}^{\bullet}+\nu_{1}^{\bullet})\},\quad R^{+}_{\k'}=\{\nu_{0}^{\bullet},\nu_{1}^{\bullet}\},
			\end{align*}
			the set $\Gamma$ and the partition function $\wp$
			\begin{align*}
				\Gamma=&\{(\nu_{0}^{\bullet}+\nu_{1}^{\bullet})/2, (\nu_{1}^{\bullet}-\nu_{0}^{\bullet})/2, \nu_{1}^{\bullet}\}, \quad m_{\beta}=1,\quad \forall\beta\in \Gamma\\
				\wp(a_{0}\nu_{0}^{\bullet}+a_{1}\nu_{1}^{\bullet})=&\#\{(m_{1},m_{2},m_{3})\in \N_{0}\mid a_{0}=\frac{m_{1}-m_{2}}{2},\quad a_{1}=\frac{m_{1}+m_{2}}{2}+m_{3}\}\\
				=&\begin{cases}
					0,\quad \text{if } a_{1}\notin \N_{0}\text{ or } 2a_{0}\notin \Z\\
					\max(0,a_{1}-a_{0}-\max(0,-2a_{0})+1),\quad \text{else.} 
				\end{cases}
			\end{align*}
			For $n_{1},n_{2},n_{3}\in \N_{0}$ the multiplicity $m(z_{1},z_{2},z_{3})=[\varrho(n_{1},n_{2},n_{3})\mycolon 1_{K'}]$ can be computed by
			$$m(n_{1},n_{2},n_{3})=\sum_{w\in W}\sgn(w)\wp(\r{w(n_{1}\omega_{1}+n_{2}\omega_{2}+n_{3}\omega_{3}+\rho)}-\r{\rho}).$$
			A tedious, but straight-forward  evaluation of this expression leads to:
			\begin{prop}
				For $n_{1},n_{2},n_{3}\in \N_{0}$ the representation $\varrho(n_{1},n_{2},n_{3})$ is $\Sp(1)'$-spherical if and only if $n_{1}=n_{3}$. In particular, a representation $\varrho(n_{1},n_{2})\in \hat{\Sp(2)}$ extends to an $\Sp(1)'$-spherical representation $\varrho(n_{1},n_{2},n_{1})$ if and only if it has highest weight that is higher or equal to the multiple $n_{3}\omega_{1}$ of the maximal root. The restriction of $\varrho(n_{1},n_{2},n_{1})\in \hat{\Sp(2)\times \Sp(1)}_{\Sp(1)'}$ to $\Sp(2)$ produces the horizontal eigenvalue with respect to the fibre representation $\varrho(n_{3})$. 
			\end{prop}
		This result is compatible with \cite[Thm. 3.8]{Bettiol}. Let us now start with the evaluation of the spectrum of the Aloff-Wallach manifold.

			\subsection{Geometric interpretation of the spectrum and electronic supplement}\label{subsection: geometric interpretation of the spectrum}
			The formulas to compute the spectrum of $(W^{1,1},g_{t_{0},t_{1}})$ stated in Theorem \ref{thm: spectrum Aloff Wallach} still need to be evaluated in order to access the spectrum. It is very cumbersome to manually calculate and analyze the spectrum from the arising combinatorial conditions. We show how to handle this problem using a suitable Python script on a computer.% that arise to determine the $\U^{\bullet}(2)$-spherical representations. 
			\subsection*{Electronic supplement}%\label{pageref: electronic supplement}
\begin{figure}[h]
				\centering
				\includegraphics[width=0.8\linewidth]{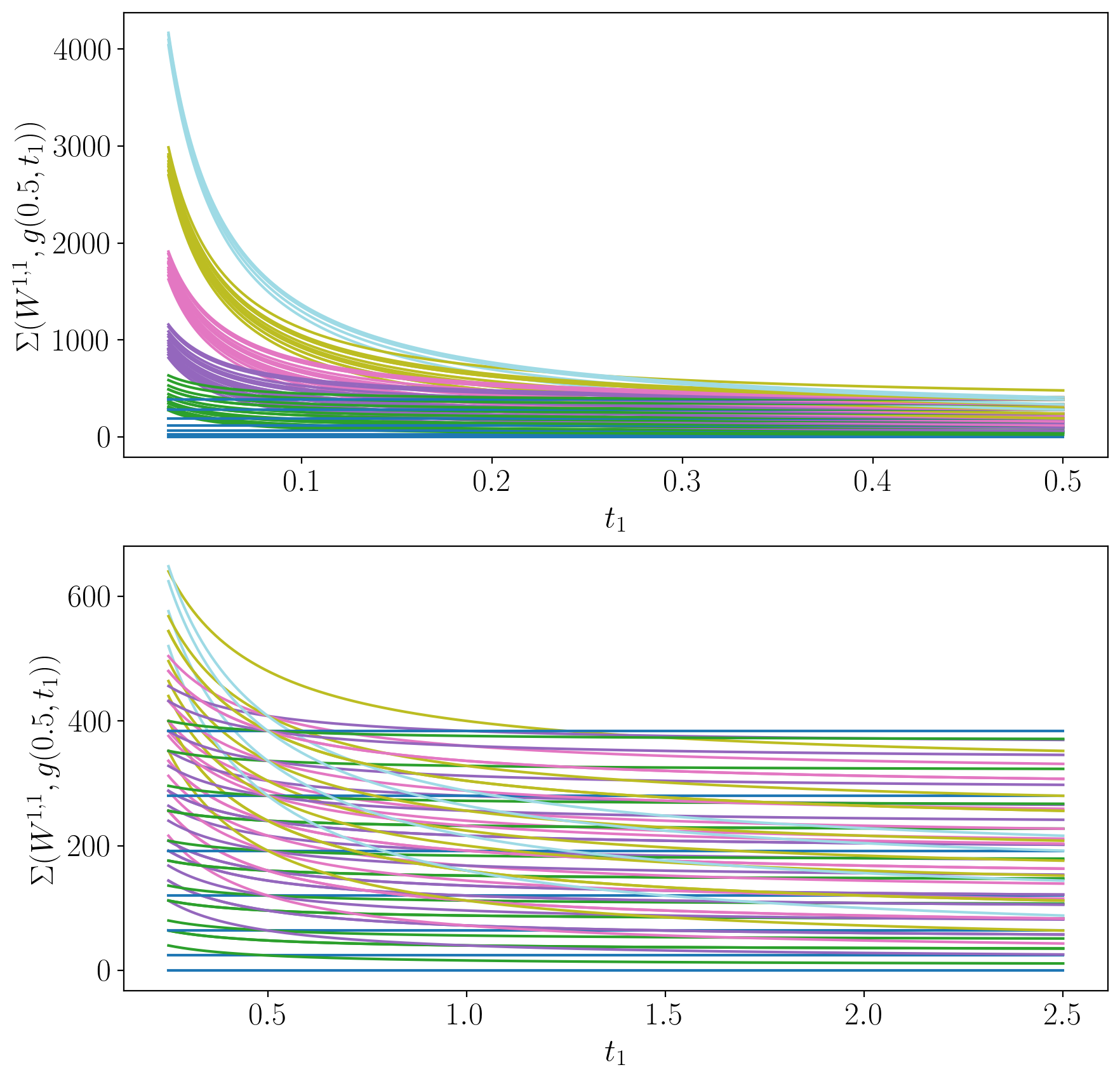}
				\caption{Spectrum with $t_{0}=0.5$, $t_{1}$ varying \cite{Agricola Henkel 24}}
				\label{fig:eigenvaluesunmodified}
			\end{figure}

			\addtocontents{toc}{\protect\setcounter{tocdepth}{1}}%damit wird Electronic supplement nicht im TOC angezeigt
			The explicit computation of eigenvalues and multiplicities is quite challenging due to the complexity of the partition function which determines the $\U^{\bullet}(2)$-spherical representations. We implemented the partition function, the multiplicity which arises from it and the eigenvalue formula into Python. The results can be found in the electronic supplement \cite{Agricola Henkel 24} which contains an overview of the first $\U^{\bullet}(2)$-spherical representations which are parameterized by $z_{1},z_{2},z_{3}\in \nolinebreak\{0,\dots,20\}$, along with their corresponding spectral data. The associated Python code is well described and can be modified to obtain the spectrum for an arbitrary finite range for $z_{1},z_{2},z_{3}\in \N_{0}$ when running it in a suitable environment, such as Jupyter Notebook. For convenience of the reader, we extracted a subset of these data points in Table \ref{table: first eigenvalues, reps and mult}.
			The electronic supplement additionally contains graphs of three different variations of the spectrum: \begin{itemize}
				\item[a)]  $t_{0}$ is fixed and $t_{1}$ is varying, see Figure \ref{fig:eigenvaluesunmodified} 
				\item[b)]  $t_{0}$ varies in dependence of $t_{1}$ such that the volume is constant, see Subsection \ref{subsection: spectrum with constant volume} 
				\item[c)] $t_{0}$ is fixed and $t_{1}$ is varying, the spectrum is divided by the total curvature 
			\end{itemize}
			The graphs can be modified by editing the Python code. Additionally, the electronic supplement contains the evaluation of a Python code which confirms that the basic eigenvalues correspond to those of $\mathbb{CP}^{2}$, see Subsection \ref{subsection: interpretation as canonical variation}. It also contains a verification that the spectrum computed by us and counted with multiplicities converges in the limiting case $t_{1}\rightarrow t_{0}=2$ to the spectrum computed by Urakawa, see the explanation below. \\
			\addtocontents{toc}{\protect\setcounter{tocdepth}{2}}%damit wird Electronic supplement nicht im TOC angezeigt
			
			Urakawa \cite{Urakawa} considered the Killing form $g_{2,2}$ of $\SU(3)$ to compute the spectrum on $W^{1,1}$. A glance at the formula for $\eta(z_{1},z_{2},z_{3},t_{0},t_{1})$ stated in \eqref{item: spectrum formula}, Theorem \ref{thm: spectrum Aloff Wallach} shows that the limiting case $t_{1}\rightarrow t_{0}=2$ forces all eigenvalues generated by $\U^{\bullet}(2)$-spherical representation with the same numbers $z_{1},z_{2}$ to converge to the same eigenvalue $\eta(z_{1},z_{2},2,2)$. The resulting formula for $\eta(z_{1},z_{2},2,2)$ is precisely the Freudenthal formula for standard normal homogeneous spaces which has been used by Urakawa. We need to dive into the computation of spherical representations in order to verify the compatibility of the results.  By Urakawa's computation a necessary and sufficient condition that $\varrho(z_{1},z_{2})\in \hat{SU(3)}$ is $ S^{1} $-spherical is $z_{1}+z_{2}\equiv 0 \mod (3)$, see \eqref{item: S1 spherical}, Theorem \ref{thm: spectrum Aloff Wallach}. Implementing Urakawa's multiplicity function for $S^{1}$-spherical representations into Python yields\\
			\begin{align*}
				\begin{array}{l|@{\hspace{0.8mm}}cccccccc}
					z_1 &  0  & 2 & 3 & 3 & 4 & 5 & 5 & 6 \\
					z_2 &  0  & 1 & 0 & 3 & 2 & 1 & 4 & 0 \\
					\mult & 1 & 32& 30 & 30 & 243 & 280 & 280 & 140, \\
				\end{array}
			\end{align*}${}$\\
			where the first $ \U^{\bullet}(2) $-spherical representations we computed are given by\\
			
			\begin{center}
				\emph{$ \begin{array}{l|@{\hspace{0.8mm}}ccccccccccccc}
						z_1 &  0  & 2 & 2 & 3 & 3 & 4 & 4 & 4 & 5 & 5& 5 & 5 & 6 \\
						z_2 &  0  & 1 & 1 & 0 & 3 & 2 & 2 & 2 & 1 & 1& 4 & 4 & 0 \\
						z_3 & 0 & 0 & 1 & 1 & 1 & 0 & 1 & 2 & 1 & 2 & 1 & 2 & 2 \\
						\mult& 1 & 8& 24 & 30 & 30 & 27 & 81 & 135 & 105 & 175& 105 & 175 & 140\\
					\end{array} $}.
			\end{center} ${}$\\
			We indeed observe that each $ S^{1} $-spherical representation displayed in the table decomposes into $ \U^{\bullet}(2) $-spherical representations and the multiplicity splits. We shall now explain what the additional information $z_{3}$ encodes. Recall that an $S^{1}$-spherical representation $\varrho(z_{1},z_{2})$ induces an $\SO(3)$-representation $\varrho(z_{1},z_{2})^{\SO(3),S^{1}}$ which acts on the representation space $(V_{\varrho(z_{1},z_{2})}^{S^{1}})^{*}$. This representation might not be irreducible and by Remark \nolinebreak  \ref{rem: GxH action} and Theorem \nolinebreak  \ref{thm: spectrum main result},  it decomposes into the desired irreducible $\SO(3)$-representations $\Phi_{\varrho(z_{1},z_{2})}=\{\varrho(z_{3})\in \hat{\SO(3)}\mid \varrho(z_{1},z_{2})\otimes\varrho(z_{3})\text{ is $\U^{\bullet}(2)$-spherical }\}$. If we leave the domain of $\SU(3)$-normal homogeneous metrics, the Laplacian no longer acts as a scalar on $(V_{\varrho(z_{1},z_{2})}^{S^{1}})^{*}$. However, it still operates scalarly  on its $\SO(3)$-irreducible subspaces with eigenvalue $\eta(z_{1},z_{2},z_{3},t_{0},t_{1})$. Hence, the multiplicities of all $\U^{\bullet}(2)$-spherical representations which contain the same  $S^{1}$-spherical representation $\varrho(z_{1},z_{2})$ have to add up and yield the multiplicity $\mult(z_{1},z_{2})$. The fact that all $S^{1}$-spherical representations extend to $\U^{\bullet}(2)$-spherical representations and the corresponding multiplicities add up could be verified for $z_{1},z_{2}\in \{0,\dots,150\}$. Depending on the computer performance, the range can be enlarged by executing the Python code described in \cite{Agricola Henkel 24}. We confirm that the spectrum counted with multiplicities converges for $t_{0}=2$ and $t_{1}\rightarrow 2$ to the spectrum computed by Urakawa. \\

\begin{figure}[h]
				\centering
				\includegraphics[width=0.7\linewidth]{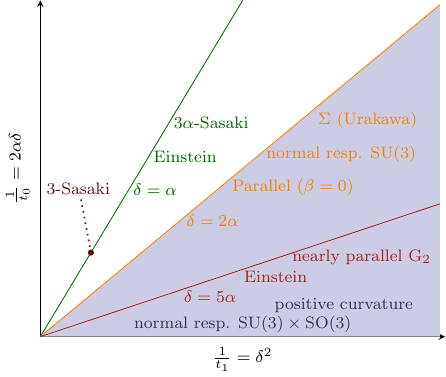}
				\caption{$\SU(3)\times \SO(3)$-naturally reductive metrics of $W^{1,1}$.}
				\label{fig:3-ad-sasaki-aw}
			\end{figure}

			This observation is reflected in Figure \ref{fig:eigenvaluesunmodified}: We clearly see the spectrum of the $\SU(3)$-normal homogeneous metric which is the vertical line $t_{1}=t_{0}=0.5$. Each eigenvalue has highest multiplicity which reflects that this metric is highly symmetric; making it possible for Urakawa \cite{Urakawa} to access the spectrum. From this point, the spectrum splits and it becomes obvious how the enlargement of the group from $\SU(3)$ to $\SU(3)\times \SO(3)$ makes it possible to capture more geometry. Each splitting corresponds to an irreducible $\SO(3)$-representation. As Figure \ref{fig:3-ad-sasaki-aw} clearly displays, the metric has positive sectional curvature for $t_{1}<t_{0}$ and this is precisely the area on which Wilking gave the $\SU(3)\times\SO(3)$-normal homogeneous realization. The nearly parallel $G_{2}$-metric, i.e. $t_{1}=0.2$, is not a special point in the spectrum. Basing on the positive curvature of the metric, this is also the area on which the Lichnerowicz-Obata estimate can be applied, see Proposition \ref{prop: Lichnerowicz} and Figure \ref{fig:picture-estimate} which yields an eigenvalue estimate from below. This estimate becomes more trivial, the further one enters the realm of naturally reductive metrics which are not normal homogeneous, i.e. $t_{1}>t_{0}$. Here the $3$-$ \alpha$-Sasaki metrics can be found. The $3$-Sasaki metric, i.e. $t_{1}=1$, is not a special point in Figure \ref{fig:eigenvaluesunmodified}.\\
			
			We remark that the spectrum on $W^{1,1}$ has been computed for the $3$-$\alpha$-Sasaki metric with $\alpha=\delta=2$ by the physicist Termonia \cite{Termonia} and his formula to compute the eigenvalues coincides with ours. However, his computations are only made for a particular metric and are difficult for a mathematician to understand.

			\subsection{Interpretation as a canonical variation of metrics over $\mathbb{CP}^2$}\label{subsection: interpretation as canonical variation}
			Consider the Riemannian submersion with totally geodesic fibres
			\begin{align*}
				(\SO(3),g_{t_{1}})\rightarrow (\SU(3)/S^{1},g_{t_{0},t_{1}})\rightarrow ( \SU(3)/\U(2),g_{t_{0}})=(\mathbb{CP}^{2},g_{t_{0}}).
			\end{align*}
			By Theorem \ref{thm: spectrum main result}, the spectrum of $( \mathbb{CP}^{2},g_{t_{0}})$ is obtained as the subset of $ \spec (\SU(3)/S^{1},g_{t_{0},t_{1}})$ in which only those $\U^{\bullet}(2)$-spherical representations $\varrho_{1}(z_{1},z_{2})\otimes\varrho_{2}(z_{3})$ contribute on which $\SO(3)$ restricts trivially, i.e. $z_{3}=0$. This can be verified for finite many eigenvalues \cite{Agricola Henkel 24}.
			Moreover, this Riemannian submersion is the reason for the bundle-like structure of the eigenvalues under the deformation of $ t_{1} $. Two eigenvalues are in the same "bundle" if they come from the same eigenvalue of the fibre $\SO(3)$ which is uniquely determined by the value of $z_{3}$.
			\subsubsection{Comparison to \cite{Semmelmann Nagy}}\label{NaSe}
			%----------------------------------------------------
			Nagy and Semmelmann \cite{Semmelmann Nagy} give eigenvalue estimates of the first eigenvalue of the canonical variation of a compact, $4n+3$-dimensional $3$-Sasaki manifold with non-constant sectional curvature, i.e. $t_{0}=0.5$ and $t_{1}$ varies. They compute the first $(\eta_{1}^{\mathcal{H}}=8)$, and in the case of a $7$-dimensional manifold with $\SO(3)$-symmetry the second $(\eta_{2}^{\mathcal{H}}=16)$ eigenvalue of the horizontal Laplacian explicitly \cite[Thm. 1.1, Thm. 1.3]{Semmelmann Nagy}. By Corollary \ref{cor: first eigenvalue horizontal} we can confirm both eigenvalues $\eta_{1}^{\mathcal{H}}=\eta^{\mathcal{H}}(2,1,1)=8$, $\eta_{2}^{\mathcal{H}}=\nolinebreak\eta^{\mathcal{H}}(4,2,2)=\nolinebreak 16$. The first horizontal eigenvalue provides an estimate of the first basic eigenvalue  $\eta_{1}^{B}\geq\nolinebreak 8(n+\nolinebreak1)$ \cite[Cor 3.10]{Semmelmann Nagy} which we improved in the homogeneous case, see Corollary \ref{cor: estimates basic eigenvalue}. In the particular situation of $G\in \{ \SU(m), \SO(2m), E_{6}, E_{7}, E_{8} \}$, the inequality $\eta_{1}^{B}\geq\nolinebreak 8(n+2)$ holds, with equality for $n=1$ on ($W^{1,1},g_{0.5,t_{1}})$. The estimate of the first eigenvalue \cite[Thm. 3.17 ]{Semmelmann Nagy} contains a mistake, as they do not take the estimate of the first basic eigenvalue into account. This was confirmed in a discussion with the authors. Keeping this in mind, the corrected estimate yields $f_{1}(t,n)\leq \eta_{1}(g_{0.5,t_{1}})$ where
			\begin{align*}
				f_{1}(t_{1},n)=\begin{cases}
					8(n+t_{1}^{-1}),\quad t_{1}\geq 1
					\\
					8(n+1),\quad t_{1}\leq 1.
				\end{cases}
			\end{align*}
			Corollary \ref{cor: first eigenvalue horizontal} proves that $z_{3}\in \N_{0}$ generates the first eigenvalue $4\cdot2z_{3}$ of the horizontal Laplacian of $(W^{1,1},g_{0.5,t_{1}})$ restricted to the eigenspace of the vertical Laplacian with respect to the representation $\varrho(z_{3})$. Those values equal our estimates, see Remark \ref{rem: horizontal eigenvalue estimates for each z} and are proven to be attained \cite[Prop. 3.4, 3.15]{Semmelmann Nagy}.

\begin{figure}
				\centering
				\includegraphics[width=0.68\linewidth]{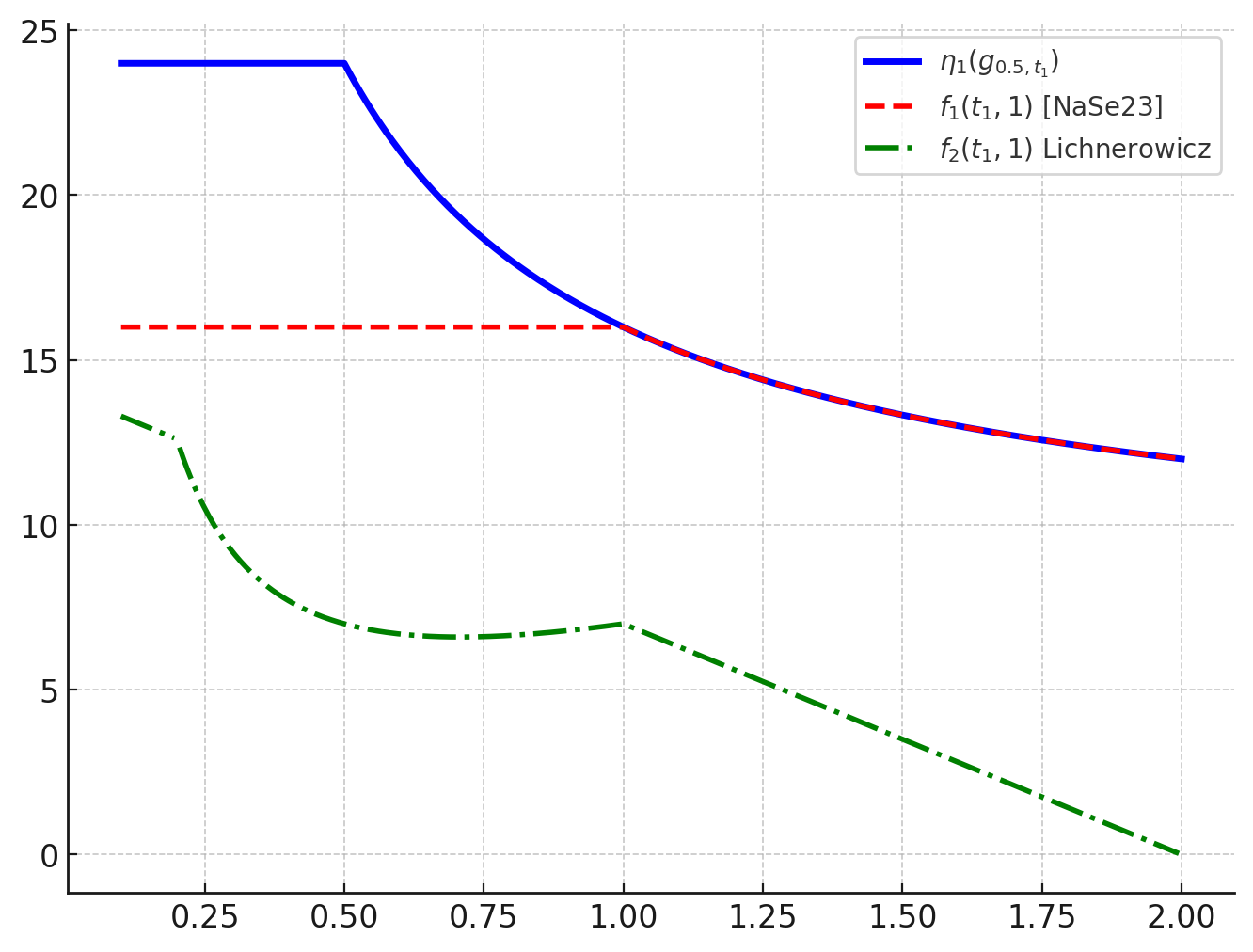}
				\caption{First eigenvalue $\eta_{1}(g_{{0.5},t_{1}})$ depending on $t_{1}$ and its estimates. Nagy, Semmelmann: $f_{1}(t_{1},1)$, Lichnerowicz, Obata: $f_{2}(t_{1},1)$}
				\label{fig:picture-estimate}
			\end{figure}

			The $\SU(3)\times\SO(3) $-normal homogeneous metrics of $(W^{1,1},g_{t_{0},t_{1}})$ have strictly positive sectional curvature. Moreover, it turns out that any homogeneous $3$-$(\alpha,\delta)$-Sasaki manifold contains metrics of positive Ricci curvature. The Lichnerowicz-Obata theorem \cite[p. 82]{Chavel} relates this curvature property to the first eigenvalue of the manifold.
			\begin{prop}\label{prop: Lichnerowicz}
				Let $(M,g_{t_{0},t_{1}})$ be a positive homogeneous $3$-$(\alpha,\delta)$-Sasaki manifold of dimension $4n+3$ where $t_{0}=\frac{1}{2}$ is fixed. The function
				\begin{align*}
					f_{2}(t_{1},n)=\begin{cases}
						\frac{\left(4 n + 3\right) \left(2 n - 3 t_{1} + 4\right)}{2 n + 1},\quad \frac{2n+4}{3}>t_{1}>1 \text{ or } \frac{1}{2n+3}>t_{1}
						\\\frac{2 n t_{1}^{2} \left(4 n + 3\right) + 4 n + 3}{t_{1} \left(2 n + 1\right)},\quad 1>t_{1}>\frac{1}{2n+3}
					\end{cases}
				\end{align*} gives a lower bound of the Ricci curvature $\Ric\geq (4n+2)\frac{f_{2}(t_{1},n)}{4n+3}>0$. By the Lichnerowicz-Obata theorem it yields a lower estimate of the first eigenvalue: $\eta_{1}(g_{0.5,t_{1}})\geq f_{2}(t_{1},n)$ with equality if and only if $M\cong S^{4n+3}$. The first basic eigenvalue is bounded from below by $\lim_{t_{1}\rightarrow 0}f_{2}(t_{1},n)=\frac{(2 (2 + n) (4n+3))}{( 2 n+1)}$.\\
				This results for the $7$-dimensional case in
				\begin{align*}
					f_{2}(t_{1},1)=\begin{cases}
						14-7t_{1},\quad 2>t_{1}>1 \text{ or } \frac{1}{5}>t_{1}
						\\\frac{7}{3t_{1}}+\frac{14t_{1}}{3},\quad 1>t_{1}>\frac{1}{5}
					\end{cases}
				\end{align*}
				with lower bound for the first basic eigenvalue $\lim_{t_{1}\rightarrow 0}f_{2}(t_{1},n)=14$.
			\end{prop}
			\begin{proof}%\ref{prop: Lichnerowicz}
				The Lichnerowicz-Obata theorem \cite[p. 82]{Chavel} states that if $\Ric\geq C(4n+3-1)$ for a $C>0$ then $\eta_{1}\geq C(4n+3) $ with equality if and only if the manifold is isometric to a round sphere of radius $1/\sqrt{C}$. The function $f_{2}(t_{1},n)$ relates to $C(t_{1},n)$ by $f_{2}=(4n+3)C$. The formula of the Ricci curvature of a $3$-$(\alpha,\delta)$-Sasaki manifold can be found in \cite[Prop. 2.3.3]{Agricola Dileo 20} and is given by
				\begin{align*}
					\Ric(X,Y)=2\alpha(2\delta(n+2)-3\alpha)g(X,Y)+2(\alpha-\delta)((2n+3)\alpha-\delta)\sum_{i=1}^{3}\eta_{i}(X)\eta_{i}(Y).
				\end{align*}
				We set $||X||=1$ and get for $\eta_{i}(X_{i}):=a_{i}$:
				\begin{align*}
					\Ric(X,X)=2\alpha(2\delta(n+2)-3\alpha)+2(\alpha-\delta)((2n+3)\alpha-\delta)\sum_{i=1}^{3}a_{i}^{2}.
				\end{align*}
				In our concrete situation, we have $2=t_{0}^{-1}=2\alpha\delta $, i.e. $\alpha=1/\delta$. Plugging this into the Ricci curvature, yields
				\begin{align*}
					\Ric(X,X)=4(n+2)-\frac{6}{\delta^{2}}+2\left(\frac{1-\delta^{2}}{\delta}\right)\left(\frac{2n+3-\delta^{2}}{\delta}\right)\sum_{i=1}^{3}a_{i}^{2}.			
				\end{align*}
				If the left hand side $4(n+2)-\frac{6}{\delta^{2}}$ is negative, i.e. $t_{1}>\frac{2n+4}{3}$ and $X$ is horizontal, the Ricci curvature becomes negative. So, assume that $t_{1}<\frac{2n+4}{3}$. If the summand on the right hand side is greater than $0$, i.e. $t_{1}^{-1}=\delta^{2}<1$ or if  $t_{1}^{-1}=\delta^{2}>2n+3$, we get the following sharp lower bound which is realized by horizontal vector fields:
				\begin{align*}
					\Ric(X,X)\geq 4(n+2)-\frac{6}{\delta^{2}}=C(4n+2).	
				\end{align*}
				From this we compute that $f_{2}=\frac{\left(4 n + 3\right) \left(2 n - 3 t_{1} + 4\right)}{2 n + 1}$. This covers the first case. For the second case $1>t_{1}>\frac{1}{2n+3}$, i.e. $2n+3>\delta^{2}>1$, the summand to the right hand side becomes negative. Assuming that $X$ is vertical, the Ricci curvature simplifies to:
				\begin{align*}
					\Ric(X,X)\geq 4(n+2)-\frac{6}{\delta^{2}}+2\left(\frac{1-\delta^{2}}{\delta}\right)\left(\frac{(2n+3)-\delta^{2}}{\delta}\right)=C (4n+2)
					%	&\Rightarrow (4n+3)C=(4(n+2)-\frac{6}{\delta^{2}}+2\left(\frac{1-\delta^{2}}{\delta}\right)\left(\frac{(2n+3)-\delta^{2}}{\delta}\right))\frac{4n+3}{4n+2}\\
					%	&=\frac{4 \delta^{4} n + 3 \delta^{4} + 8 n^{2} + 6 n}{\delta^{2} \left(2 n + 1\right)}=\frac{2 n t_{1}^{2} \left(4 n + 3\right) + 4 n + 3}{t_{1} \left(2 n + 1\right)}
				\end{align*}
				From this we compute that $f_{2}=\frac{2 n t_{1}^{2} \left(4 n + 3\right) + 4 n + 3}{t_{1} \left(2 n + 1\right)}$.
			\end{proof}
			The estimate of the first eigenvalue obtained in \cite{Semmelmann Nagy} is in all cases better than the Lichnerowicz-Obata estimate. Our improvement of the estimate of the first basic eigenvalue of \cite{Semmelmann Nagy}, see Corollary \ref{cor: estimates basic eigenvalue}, yields precisely the line of the first eigenvalue, see Figure \ref{fig:picture-estimate}.

			\subsection{Spectrum with constant volume}\label{subsection: spectrum with constant volume}
\begin{figure}[h!]
				\centering
				\includegraphics[width=0.6\linewidth]{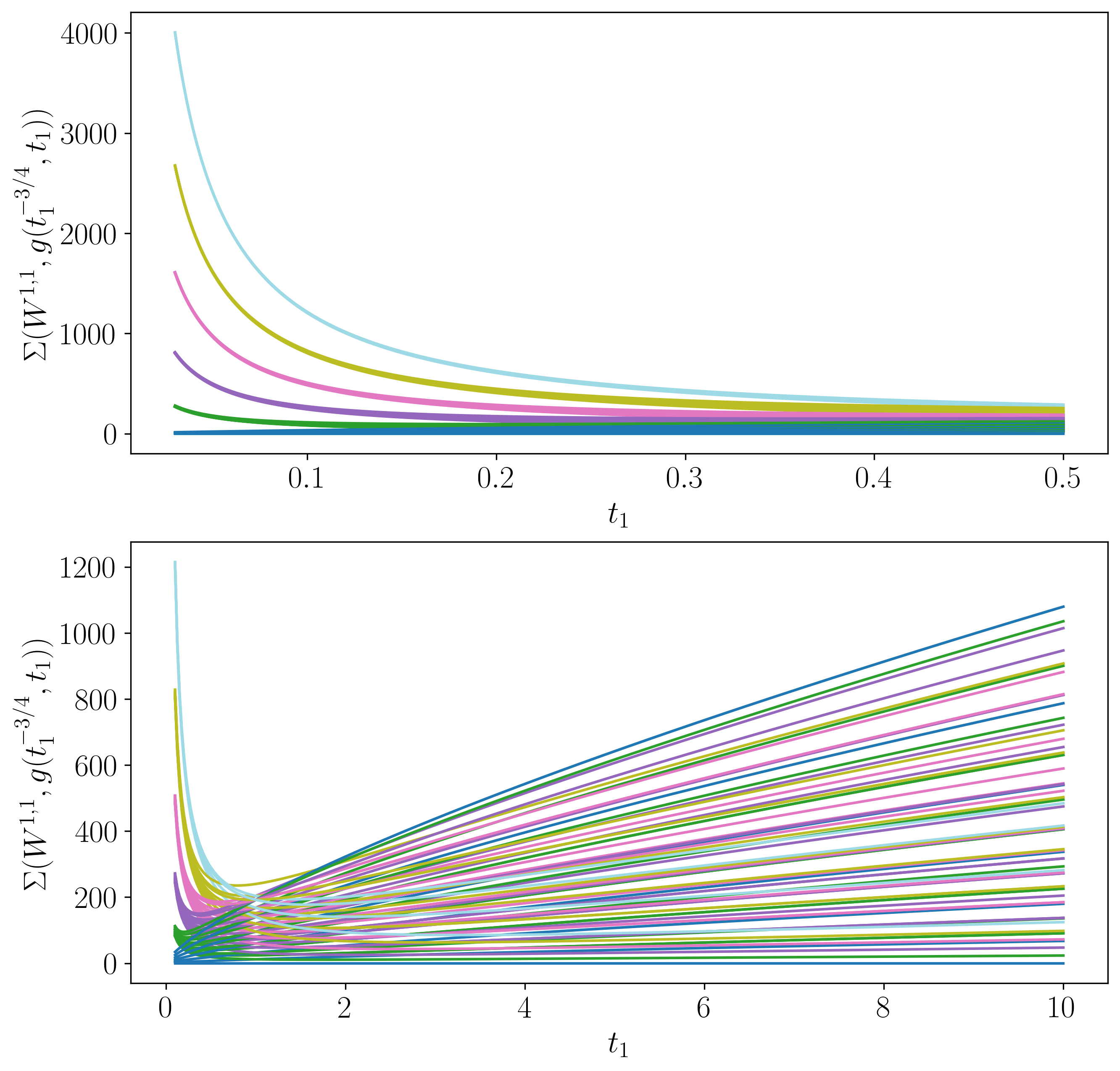}
				\caption{Spectrum with constant volume \cite{Agricola Henkel 24}}
				\label{fig:eigenvaluesvolconst}
			\end{figure}

			The spectrum naturally depends on conformal changes of the metric and therefore on the volume: If the metric is scaled by a constant $c>0$, the spectrum is inversely proportional to $c$, where the volume is proportional to $c^{n/2}$. Berger \cite{Berger first eigenvalue} considered constant conformal changes by posing the problem of determining for which compact manifolds $M$ it is possible to find a constant $C(M)$ such that for any Riemannian metric $g$ of $M$ the first eigenvalue $\eta_{1}(g)$ satisfies
			\begin{align}
				\eta_{1}(g)\vol(M,g)^{2/n}\leq C(M).\label{eq: bounded first eigenvalue}
			\end{align}
			Berger was inspired by Hersch \cite{Hersch}, who proved that $M=S^{2}$ satisfies this inequality for  $C=8\pi$. Urakawa \cite{Urakawa on the least eigenvalue} proved that such an upper bound does not exist if $M$ is a Lie group which is not a torus. A positive answer was given by El Soufi and Ilias \cite{El Soufi Ilias} who proved that \eqref{eq: bounded first eigenvalue} holds for any manifold if $C=C(M,[g])$ additionally depends on the conformal class $[g]$ of $g$.  We will see that there does not exist a constant $C=C(W^{1,1})$ which does not depend on the conformal classes $[g_{t_{0},t_{1}}]$ such that \eqref{eq: bounded first eigenvalue} holds and the same argument works for any other $3$-$(\alpha,\delta)$-Sasaki manifold. Integrating the volume form of $\SU(3)$ along each fibre yields the formula \begin{align*}
				\vol(\SU(3),g_{t_{0},t_{1}})=\vol(W^{1,1},g_{t_{0},t_{1}})\cdot\vol(S^{1},g_{t_{0},t_{0}}).
			\end{align*} The volume form of $(\SU(3),g_{t_{0},t_{1}})$ is left-invariant and hence uniquely determined by its value at the origin, where it is proportional to $\sqrt{t_{0}^{5}t_{1}^{3}}$. Analogously, the volume of each fibre $(S^{1},g_{t_{0},t_{0}})$ is proportional to $\sqrt{t_{0}}$. We conclude that $\vol(W^{1,1},g_{t_{0},t_{1}})  $ is proportional to $t_{0}^{2}\sqrt{t_{1}^{3}}$ so it is constant if $t_{0}=t_{1}^{-3/4}$. In this case, each eigenvalue is proportional to $t_{1}^{3/4}$ or is the sum of a horizontal eigenvalue which is proportional to $t_{1}^{3/4}$ and a vertical eigenvalue which is proportional to $t_{1}^{-1}$. We obtain that
			\begin{align*}
				\lim_{t_{1}\rightarrow \infty}	\eta_{1}(g_{t_{1}^{-3/4},t_{1}})\cdot \vol(W^{1,1},g_{t_{1}^{-3/4},t_{1}})^{2/7}=\infty,\quad \lim_{t_{1}\rightarrow 0}	\eta_{1}(g_{t_{1}^{-3/4},t_{1}})\cdot \vol(W^{1,1},g_{t_{1}^{-3/4},t_{1}})^{2/7}=0,
			\end{align*}
			giving another class of counterexamples to Berger's question. A plot of the spectrum under this deformation can be found in Figure \ref{fig:eigenvaluesvolconst}.			
			\\\\

	\appendix
	\section{Table of eigenvalues of $(W^{1,1},g_{t_{0},t_{1}})$ for small $z_{1},z_{2},z_{3}$}
	As explained in the beginning of Subsection \ref{subsection: geometric interpretation of the spectrum}, the electronic supplement \cite{Agricola Henkel 24} displays the first eigenvalues and multiplicities generated by the first $\U^{\bullet}(2)$-spherical representations $\varrho(z_{1},z_{2},z_{3})$, where $z_{1},z_{2},z_{3}\in \{0,\ldots,20\}$ and provides a suitable Python code which can be adapted to extend the calculations to an arbitrary, finite range for $z_{1},z_{2},z_{3}\in \N_0$. For the reader's convenience we extract the first $78$ entries of this list in Table \ref{table: first eigenvalues, reps and mult}.

	{\renewcommand{\frac}[2]{#1/#2}
		\begin{table}[h!]
			\centering
			\renewcommand{\arraystretch}{1.2}
			\normalsize
			\begin{tabular}{cc}
				\begin{minipage}[t]{0.50\linewidth}
					\centering
					\begin{tabular}{c c c c c}
						\toprule
						\textbf{$z_1$} & \textbf{$z_2$} & \textbf{$z_3$} & \textbf{$\eta$} & mult \\
						\midrule
						0 & 0 & 0 & $0$ & $1$ \\
						2 & 1 & 0 & $\frac{12}{t_{0}}$ & $8$ \\
						2 & 1 & 1 & $\frac{8}{t_{1}} + \frac{4}{t_{0}}$ & $24$ \\
						3 & 0 & 1 & $\frac{8}{t_{1}} + \frac{16}{t_{0}}$ & $30$ \\
						3 & 3 & 1 & $\frac{8}{t_{1}} + \frac{16}{t_{0}}$ & $30$ \\
						4 & 2 & 0 & $\frac{32}{t_{0}}$ & $27$ \\
						4 & 2 & 1 & $\frac{8}{t_{1}} + \frac{24}{t_{0}}$ & $81$ \\
						4 & 2 & 2 & $\frac{24}{t_{1}} + \frac{8}{t_{0}}$ & $135$ \\
						5 & 1 & 1 & $\frac{8}{t_{1}} + \frac{40}{t_{0}}$ & $105$ \\
						5 & 1 & 2 & $\frac{24}{t_{1}} + \frac{24}{t_{0}}$ & $175$ \\
						5 & 4 & 1 & $\frac{8}{t_{1}} + \frac{40}{t_{0}}$ & $105$ \\
						5 & 4 & 2 & $\frac{24}{t_{1}} + \frac{24}{t_{0}}$ & $175$ \\
						6 & 0 & 2 & $\frac{24}{t_{1}} + \frac{48}{t_{0}}$ & $140$ \\
						6 & 3 & 0 & $\frac{60}{t_{0}}$ & $64$ \\
						6 & 3 & 1 & $\frac{8}{t_{1}} + \frac{52}{t_{0}}$ & $192$ \\
						6 & 3 & 2 & $\frac{24}{t_{1}} + \frac{36}{t_{0}}$ & $320$ \\
						6 & 3 & 3 & $\frac{48}{t_{1}} + \frac{12}{t_{0}}$ & $448$ \\
						6 & 6 & 2 & $\frac{24}{t_{1}} + \frac{48}{t_{0}}$ & $140$ \\
						7 & 2 & 1 & $\frac{8}{t_{1}} + \frac{72}{t_{0}}$ & $243$ \\
						7 & 2 & 2 & $\frac{24}{t_{1}} + \frac{56}{t_{0}}$ & $405$ \\
						7 & 2 & 3 & $\frac{48}{t_{1}} + \frac{32}{t_{0}}$ & $567$ \\
						7 & 5 & 1 & $\frac{8}{t_{1}} + \frac{72}{t_{0}}$ & $243$ \\
						7 & 5 & 2 & $\frac{24}{t_{1}} + \frac{56}{t_{0}}$ & $405$ \\
						7 & 5 & 3 & $\frac{48}{t_{1}} + \frac{32}{t_{0}}$ & $567$ \\
						8 & 1 & 2 & $\frac{24}{t_{1}} + \frac{84}{t_{0}}$ & $400$ \\
						8 & 1 & 3 & $\frac{48}{t_{1}} + \frac{60}{t_{0}}$ & $560$ \\
						8 & 4 & 0 & $\frac{96}{t_{0}}$ & $125$ \\
						8 & 4 & 1 & $\frac{8}{t_{1}} + \frac{88}{t_{0}}$ & $375$ \\
						8 & 4 & 2 & $\frac{24}{t_{1}} + \frac{72}{t_{0}}$ & $625$ \\
						8 & 4 & 3 & $\frac{48}{t_{1}} + \frac{48}{t_{0}}$ & $875$ \\
						8 & 4 & 4 & $\frac{80}{t_{1}} + \frac{16}{t_{0}}$ & $1125$ \\
						8 & 7 & 2 & $\frac{24}{t_{1}} + \frac{84}{t_{0}}$ & $400$ \\
						8 & 7 & 3 & $\frac{48}{t_{1}} + \frac{60}{t_{0}}$ & $560$ \\
						9 & 0 & 3 & $\frac{48}{t_{1}} + \frac{96}{t_{0}}$ & $385$ \\
						9 & 3 & 1 & $\frac{8}{t_{1}} + \frac{112}{t_{0}}$ & $462$ \\
						9 & 3 & 2 & $\frac{24}{t_{1}} + \frac{96}{t_{0}}$ & $770$ \\
						9 & 3 & 3 & $\frac{48}{t_{1}} + \frac{72}{t_{0}}$ & $1078$ \\
						9 & 3 & 4 & $\frac{80}{t_{1}} + \frac{40}{t_{0}}$ & $1386$ \\
						9 & 6 & 1 & $\frac{8}{t_{1}} + \frac{112}{t_{0}}$ & $462$ \\
						\bottomrule
					\end{tabular}
					
				\end{minipage} &
				\begin{minipage}[t]{0.50\linewidth}
					\centering
					\begin{tabular}{c c c c c}
						\toprule
						\textbf{$z_1$} & \textbf{$z_2$} & \textbf{$z_3$} & \textbf{$\eta$} & mult \\
						\midrule
						9 & 6 & 2 & $\frac{24}{t_{1}} + \frac{96}{t_{0}}$ & $770$ \\
						9 & 6 & 3 & $\frac{48}{t_{1}} + \frac{72}{t_{0}}$ & $1078$ \\
						9 & 6 & 4 & $\frac{80}{t_{1}} + \frac{40}{t_{0}}$ & $1386$ \\
						9 & 9 & 3 & $\frac{48}{t_{1}} + \frac{96}{t_{0}}$ & $385$ \\
						10 & 2 & 2 & $\frac{24}{t_{1}} + \frac{128}{t_{0}}$ & $810$ \\
						10 & 2 & 3 & $\frac{48}{t_{1}} + \frac{104}{t_{0}}$ & $1134$ \\
						10 & 2 & 4 & $\frac{80}{t_{1}} + \frac{72}{t_{0}}$ & $1458$ \\
						10 & 5 & 0 & $\frac{140}{t_{0}}$ & $216$ \\
						10 & 5 & 1 & $\frac{8}{t_{1}} + \frac{132}{t_{0}}$ & $648$ \\
						10 & 5 & 2 & $\frac{24}{t_{1}} + \frac{116}{t_{0}}$ & $1080$ \\
						10 & 5 & 3 & $\frac{48}{t_{1}} + \frac{92}{t_{0}}$ & $1512$ \\
						10 & 5 & 4 & $\frac{80}{t_{1}} + \frac{60}{t_{0}}$ & $1944$ \\
						10 & 5 & 5 & $\frac{120}{t_{1}} + \frac{20}{t_{0}}$ & $2376$ \\
						10 & 8 & 2 & $\frac{24}{t_{1}} + \frac{128}{t_{0}}$ & $810$ \\
						10 & 8 & 3 & $\frac{48}{t_{1}} + \frac{104}{t_{0}}$ & $1134$ \\
						10 & 8 & 4 & $\frac{80}{t_{1}} + \frac{72}{t_{0}}$ & $1458$ \\
						11 & 1 & 3 & $\frac{48}{t_{1}} + \frac{144}{t_{0}}$ & $1001$ \\
						11 & 1 & 4 & $\frac{80}{t_{1}} + \frac{112}{t_{0}}$ & $1287$ \\
						11 & 4 & 1 & $\frac{8}{t_{1}} + \frac{160}{t_{0}}$ & $780$ \\
						11 & 4 & 2 & $\frac{24}{t_{1}} + \frac{144}{t_{0}}$ & $1300$ \\
						11 & 4 & 3 & $\frac{48}{t_{1}} + \frac{120}{t_{0}}$ & $1820$ \\
						11 & 4 & 4 & $\frac{80}{t_{1}} + \frac{88}{t_{0}}$ & $2340$ \\
						11 & 4 & 5 & $\frac{120}{t_{1}} + \frac{48}{t_{0}}$ & $2860$ \\
						11 & 7 & 1 & $\frac{8}{t_{1}} + \frac{160}{t_{0}}$ & $780$ \\
						11 & 7 & 2 & $\frac{24}{t_{1}} + \frac{144}{t_{0}}$ & $1300$ \\
						11 & 7 & 3 & $\frac{48}{t_{1}} + \frac{120}{t_{0}}$ & $1820$ \\
						11 & 7 & 4 & $\frac{80}{t_{1}} + \frac{88}{t_{0}}$ & $2340$ \\
						11 & 7 & 5 & $\frac{120}{t_{1}} + \frac{48}{t_{0}}$ & $2860$ \\
						11 & 10 & 3 & $\frac{48}{t_{1}} + \frac{144}{t_{0}}$ & $1001$ \\
						11 & 10 & 4 & $\frac{80}{t_{1}} + \frac{112}{t_{0}}$ & $1287$ \\
						12 & 0 & 4 & $\frac{80}{t_{1}} + \frac{160}{t_{0}}$ & $819$ \\
						12 & 3 & 2 & $\frac{24}{t_{1}} + \frac{180}{t_{0}}$ & $1400$ \\
						12 & 3 & 3 & $\frac{48}{t_{1}} + \frac{156}{t_{0}}$ & $1960$ \\
						12 & 3 & 4 & $\frac{80}{t_{1}} + \frac{124}{t_{0}}$ & $2520$ \\
						12 & 3 & 5 & $\frac{120}{t_{1}} + \frac{84}{t_{0}}$ & $3080$ \\
						12 & 6 & 0 & $\frac{192}{t_{0}}$ & $343$ \\
						12 & 6 & 1 & $\frac{8}{t_{1}} + \frac{184}{t_{0}}$ & $1029$ \\
						12 & 6 & 2 & $\frac{24}{t_{1}} + \frac{168}{t_{0}}$ & $1715$ \\
						12 & 6 & 3 & $\frac{48}{t_{1}} + \frac{144}{t_{0}}$ & $2401$ \\
						\bottomrule
					\end{tabular}
					
				\end{minipage} \\
			\end{tabular}
			\caption{First eigenvalues of $(W^{1,1},g_{t_{0},t_{1}})$ and their multiplicities}
			\label{table: first eigenvalues, reps and mult}
	\end{table}}

\end{document}